\documentclass{article}

\usepackage{amsmath}

\numberwithin{equation}{section} 

\counterwithin{figure}{section}
\counterwithin{table}{section}

\usepackage{amsxtra}

\usepackage{mathtools}

\usepackage{amssymb}

\usepackage[all,arrow]{xy}

\usepackage{geometry}
\usepackage{threeparttable}
\usepackage{diagbox}
\usepackage{paralist,bbding,pifont}

\usepackage[thmmarks]{ntheorem}
{
	\theoremstyle{nonumberplain}
	\theoremheaderfont{\bfseries}
	\theorembodyfont{\normalfont}
	\theoremsymbol{\mbox{$\Box$}}
	\newtheorem{pf}{Proof}
}
{
	\theoremstyle{plain}
	\theoremheaderfont{\bfseries}
	\theorembodyfont{\normalfont}
	\newtheorem{rmk}{Remark}[section]
}
{
	\theoremstyle{plain}
	\theoremheaderfont{\bfseries}
	\theorembodyfont{\normalfont}
	
}
{
	\theoremstyle{plain}
	\theoremheaderfont{\bfseries}
	\theorembodyfont{\normalfont}
	
}
{
	\theoremstyle{plain}
	\theoremheaderfont{\bfseries}
	\theorembodyfont{\normalfont}
	\newtheorem*{ackn}{Acknowledgements.}
}
\newtheorem{thm}{Theorem}[section]

\newtheorem{prop}{Proposition}[section]
\newtheorem{lem}{Lemma}[section]
{
	\theoremstyle{plain}
	\theoremheaderfont{\bfseries}
	\theorembodyfont{\normalfont}
	\newtheorem{defn}{Definition}[section]
}
{
	\theoremstyle{plain}
	\theoremheaderfont{\bfseries}
	\theorembodyfont{\normalfont}
	
}
\newtheorem{Q}{Problem}

\usepackage{nicematrix,tikz}
\usepackage{dsfont}
\usetikzlibrary{matrix}
\usetikzlibrary{decorations.pathmorphing} 
\usetikzlibrary{positioning} 

\usepackage{tabularx,amsmath}
\usepackage{float}
\usepackage{booktabs}	
\usepackage{caption}   

\usepackage{hyperref}
\hypersetup{hypertex=true,
	colorlinks=true,
	linkcolor=red,
	anchorcolor=green,
	citecolor=blue}

\title{On $7$-manifolds with $b_{2}=2$: diffeomorphism classification and nonconnected moduli spaces of positive Ricci curvature metrics}
\author{Fupeng Xu}
\date{}

\begin{document}
	
	\maketitle
	
	\begin{abstract}
		We derive the $s$-invariants of certain simply connected $7$-manifolds whose second homology groups are isomorphic to $\mathbb{Z}^{2}$. We apply the $s$-invariants to give a partial classification of simply connected total spaces of circle bundles over $\left(\mathbb{C}P^{1}\times\mathbb{C}P^{2}\right)\#\mathbb{C}P^{3}$ up to diffeomorphism. As an application, we show that there is a simply connected $7$-manifold whose space and moduli space of positive Ricci curvature metrics both have infinitely many path components. We also determine bordism groups $\Omega_{8}^{Spin}\left(K_{2}\right)$ and $\Omega_{8}^{Spin}\left(K_{2};\mathrm{pr}_{1}^{*}\gamma^{1}\right)$ that are required in the deduction of $s$-invariants.
	\end{abstract}
	
	
	\section{Introduction}

If no otherwise is stated, manifolds are smooth, connected and closed. Homeomorphisms and diffeomorphims between oriented manifolds preserve orientation. Lie groups are connected and compact. Fiber bundles over oriented manifolds with oriented fiber manifold are oriented.

Manifold classification is a fundamental problem in topology and has wide applications in many areas of mathematics, while a complete classification becomes increasingly complicated as dimension increases. Except for dimension $4$, homeomorphism and diffeomorphism classification of simply connected manifolds are known up to dimension $6$ (\cite{Smale62Spin5mfd,Barden65simplycnt5mfd} for dimension $5$, \cite{Wall1966Classification,Jupp1973,Zhubr00} for dimension $6$). When $n\geqslant7$, a complete classification of simply connected $n$-manifolds remains unclear, and classification of manifolds with extra structures are frequently considered. We ask the following problem:

\begin{Q}\label{Problem: Main}
	Assume we have a family of simply connected $7$-manifolds that are the total spaces of circle bundles over a fixed manifold. Give their classification up to diffeomorphism.
\end{Q} 

Solution to Problem \ref{Problem: Main} can be related to differential geometry. In differential geometry we ask the following basic problem:

\begin{Q}
	Given a smooth manifold $M$, determine whether $M$ admits a Riemannian metric that satisfies certain curvature condition $P$ (abbreviated as $P$ metric for convenience). When $M$ admits a $P$ metric, determine ``how many'' $P$ metrics $M$ can support. For example study whether the space and moduli space of $P$ metrics on $M$ has trivial, finite or finitely generated homotopy groups.
\end{Q}

There are some classical conditions $P$, for example positive (non-negative, negative, non-positive) scalar (Ricci, sectional) curvature, Ricci flat metrics and Einstein metrics. For positive scalar curvature, the existence problem (\cite{HITCHIN1974,SchoenYau79,GromovLawson80,Rosenberg1986I,Rosenberg1986II,ROSENBERG1986319,Stolz92,STOLZ1994}) and counting problem (\cite{HITCHIN1974,Carr88,KreckStolz1993,CS13,BER17}) have been extensively studied. While for positive Ricci curvature and positive sectional curvature these problems remain open.

Homogeneous spaces and principal bundles are rich sources of manifolds with positive Ricci curvature. By \cite[Proposition 3.4]{Nash1979} if $G$ is a compact Lie group and $H$ is a closed subgroup of $G$, then $G\big/H$ admits a metric with positive Ricci curvature if and only if its fundamental group $\pi_{1}\left(G\big/H\right)$ is finite. By \cite[Theorem 0.4]{PBJPWT98PositiveRicciPrinBdl} if $\left(Y,g_{Y}\right)$ is a Riemannian manifold with positive Ricci curvature and $P\to Y$ is a principal $H$-bundle, such that $H$ is a Lie group and $P$ has finite fundamental group, then $P$ admits an $H$-equivariant metric $g_{P}$ such that $\left(P,g_{P}\right)$ has positive Ricci curvature and $\left(P,g_{P}\right)\to\left(Y,g_{Y}\right)$ is a Riemmanian submersion. 

We focus on the case $H=S^{1}$. \cite{WZ1986} considers free circle actions on $S^{3}\times S^{5}$ and studies the corresponding orbit spaces, showing that there are infinitely many simply connected Einstein manifolds $\left(M_{m,n},g_{M_{m,n}}\right)$ with positive cosmology constants. Here $m$, $n$ are coprime integers. In particular they have positive Ricci curvature. 
	
The manifold $M_{m,n}$ can be equivalently described as the total spaces of circle bundles over $\mathbb{C}P^{1}\times\mathbb{C}P^{2}$. \cite{KS88} uses modified surgery theory and gives the homeomorphism and diffeomorphism classification of $M_{m,n}$. In particular $M_{m,n}$ and $M_{\overline{m},\overline{n}}$ are diffeomorphic if and only if
	\begin{eqnarray}\label{Equation: M_{m,n} classification}
		\begin{cases}
			\overline{n}=\pm n,\\
			\overline{m}\equiv m\mod 2^{\lambda_{2}(n)}7^{\lambda_{7}(n)}n^{2},
		\end{cases}
	\end{eqnarray}
	where
	\begin{eqnarray}
		\lambda_{2}(n)=\begin{cases}
			0,&n\equiv2,6\mod 8,\\
			1,&n\equiv1,7\mod 8,\\
			2,&n\equiv3,5\mod 8,\\
			3,&n\equiv0,4\mod 8;
		\end{cases}\label{Equation: def of lambda_2}\\
		\lambda_{7}(n)=\begin{cases}
			0, &n\equiv1,2,5,6\mod7,\\
			1, &n\equiv0,3,4\mod7.
		\end{cases}\label{Equation: def of lambda_7}
	\end{eqnarray}
	
\cite{KreckStolz1993} computes the $s$-invariants of Riemannian manifolds $\left(M_{m,n},g_{M_{m,n}}\right)$:
	$$s\left(M_{m,n},g_{M_{m,n}}\right)=-\frac{3m\left(n^{2}+3\right)\left(n^{2}-1\right)}{896n^{2}}\in\mathbb{Q}.$$
This $s$-invariant is defined for a Riemannian manifold $(M,g)$ of dimension $(4q-1)$ $(q\geqslant2)$, such that $(M,g)$ has positive scalar curvature and all its rational Pontryagin classes vanish. It is an invariant on the path components of $\mathfrak{R}_{scal>0}(M)$, the space of positive scalar curvature metrics on $M$, and when $H^{1}\left(M;\mathbb{Z}\big/2\right)=0$ its absolute value is an invariant on the path components of $\mathfrak{M}_{scal>0}(M)$, the moduli space of positive scalar curvature metrics. Here given an oriented manifold $M$, we denote by $\mathrm{Diff}(M)$ the group of orientation-preserving self-diffeomorphism on $M$; when $\mathfrak{R}_{scal>0}(M)$ is non-empty $\mathrm{Diff}(M)$ acts on $\mathfrak{R}_{scal>0}(M)$ by $\varphi\cdot g:=\left(\varphi^{-1}\right)^{*}g$ and the corresponding orbit space is defined as the moduli space $\mathfrak{M}_{scal>0}(M)$. \cite{KreckStolz1993} combines the diffeomorphism classification and shows that there is a manifold $M_{m,n}$ of which $\mathfrak{R}_{Ric>0}(M)$ and $\mathfrak{M}_{Ric>0}(M)$, the space and moduli space of positive Ricci curvature metrics on $M$, have infinitely many path components. See \cite[Chapter 1]{TW15} for more details of topology on $\mathfrak{R}_{scal>0}(M)$, $\mathrm{Diff}(M)$ and $\mathfrak{M}_{scal>0}(M)$.

Here is another example. Diagonal embeddings of $U(1)$ into $SU(3)$ induce free circle actions on $SU(3)$. \cite{Aloff1975AnIF} studies the corresponding orbit spaces and shows that there are infinitely many simply connected Riemannian manifolds $\left(N_{m,n},g_{N_{m,n}}\right)$ admitting positive sectional curvature with $m$ and $n$ coprime.
	
The manifolds $N_{m,n}$ can be equivalently characterized as the total spaces of circle bundles over the complete flag manifold $SU(3)\big/T=Fl(3)$, where $T\cong S^{1}\times S^{1}$ is a maximal torus of $SU(3)$. \cite{Kreck1991SomeNH,Kreck1991SomeNHCorrection} use modified surgery theory and classify $N_{m,n}$ up to homeomorphism and diffeomorphism. \cite{KreckStolz1993} computes the $s$-invariants of Riemannian manifolds $\left(N_{m,n},g_{N_{m,n}}\right)$, combines the diffeomorphism classification and shows that there is a manifold $N_{m,n}$ of which $\mathfrak{R}_{sec>0}(M)$ and $\mathfrak{M}_{sec>0}(M)$, the space and moduli space of positive sectional curvature metrics, are disconnected.

There are another approaches constructing manifolds with positive Ricci curvature systematically. \cite{BURDICK2019} extends the result of \cite{Perelman1997} and introduces the \textbf{core metric}. This is a kind of positive Ricci curvature metric, and projective spaces $\mathbb{C}P^{n}$, $\mathbb{H}P^{n}$, $\mathbb{O}P^{2}$ admit core metrics (\cite[Theorem C]{BURDICK2019}). The property of admitting a core metric behaves well under connected sum: if $n$-manifolds $M_{1}$ and $M_{2}$ both admit core metrics, so does their connected sum $M_{1}\#M_{2}$; in particular $M_{1}\#M_{2}$ admits a positive Ricci curvature metric (\cite[Theorem B]{BURDICK2019}, \cite[Corollary 3.18]{burdick2020}). Meanwhile assume we are given a linear $S^{p}$-bundle $E\to B$ over $q$-manifold $B$ such that $p+q\geqslant6$ and $\min\left\{p,q\right\}\geqslant2$,
$E$ admits a core metric if $B$ does; in particular $E$ admits a positive Ricci curvature metric (\cite[Theorem 1.5, Theorem C]{Reiser_2023}).

There is a more general result on enumeration of positive Ricci curvature metrics. Let $M$ be a spin manifold of dimension $4q-1$ $(q\geqslant2)$ such that it admits a \textbf{socket metric}. This is a kind of positive Ricci curvature metric and is weaker than the core metric, and a core metric is automatically a socket metric. Then $\mathfrak{R}_{Ric>0}(M)$ admits infinitely many components (\cite[Theorem B]{burdick2020}). If moreover $H^{1}\left(M;\mathbb{Z}\big/2\right)=0$ and all its rational Pontryagin classes vanish, then $\mathfrak{M}_{Ric>0}(M)$ also has infinitely many components (\cite[Theorem D]{burdick2020})

Building upon the successful classification and geometric analysis of families $M_{m,n}$ and $N_{m,n}$, we investigate a natural yet more complicated generalization: the total spaces of circle bundles over the connected sum $N=\left(\mathbb{C}P^{1}\times\mathbb{C}P^{2}\right)\#\mathbb{C}P^{3}$ which are simply connected. They can be parametrized as follows. Note that $H^{2}\left(N\right)=\mathbb{Z}\left\{\alpha,\beta,\gamma\right\}\cong\mathbb{Z}^{3}$, where $\alpha$, $\beta$ and $\gamma$ correspond to the Euler classes of Hopf bundles over $\mathbb{C}P^{1}$, $\mathbb{C}P^{2}$ and $\mathbb{C}P^{3}$ respectively. Let $M_{m,n,l}$ be the total space of circle bundle over $N$ with Euler class $m\alpha+n\beta+l\gamma$. It can be shown that $M_{m,n,l}$ is simply connected if and only if the greatest common divisor of $m$, $n$ and $l$, denoted by $\gcd(m,n,l)$ and assumed to be positive, is $1$ (Lemma \ref{Lemma: E 1cnt}).

\begin{Q}\label{Problem: Classify M_{m,n,l}}
	Classify simply connected $M_{m,n,l}$ up to diffeomorphism. Given integers $m$, $n$, $l$, $\overline{m}$, $\overline{n}$, $\overline{l}$ such that $M_{m,n,l}$ and $M_{\overline{m},\overline{n},\overline{l}}$ are simply connected, determine when $M_{m,n,l}$ and $M_{\overline{m},\overline{n},\overline{l}}$ are diffeomorphic. Express the results using algebraic equation in $m$, $n$, $l$ and $\overline{m}$, $\overline{n}$, $\overline{l}$.
\end{Q}

$M_{m,n,l}$ admits positive Ricci curvature metrics, and the reason is as follows. $\mathbb{C}P^{2}$ and $\mathbb{C}P^{3}$ admits core metrics (\cite[Theorem C]{BURDICK2019}). Hence $\mathbb{C}P^{1}\times\mathbb{C}P^{2}=S^{2}\times\mathbb{C}P^{2}$ admits a core metric (\cite[Theorem C]{Reiser_2023}) and the connected sum $N=\left(\mathbb{C}P^{1}\times\mathbb{C}P^{2}\right)\#\mathbb{C}P^{3}$ also admits a core metric (\cite[Theorem C]{BURDICK2019}). In particular $N$ admits a positive Ricci curvature metric. Since $M_{m,n,l}$ is the total space of a circle bundle over $N$, it also admits a positive Ricci curvature metric (\cite[Theorem 0.4]{PBJPWT98PositiveRicciPrinBdl}). 

It is then natural to ask the following problem:

\begin{Q}\label{Problem: (moduli) space of positive Ricci curvature metrics of M_{m,n,l}}
	Determine whether the space or moduli space of positive Ricci curvature metrics on $M_{m,n,l}$ is contractible or has infinitely many components.
\end{Q}
The powerful general results of \cite[Theorems B and D]{burdick2020} on the infiniteness of moduli space components are not directly applicable to $M_{m,n,l}$
, since we do not know whether $M_{m,n,l}$ admits a socket metric from the existing knowledge. 

In this paper we give a partial diffeomorphism classification of simply connected $M_{m,n,l}$ and answer Problem \ref{Problem: Classify M_{m,n,l}} partially. We also show that there are manifolds $M_{m,n,l}$ such that $\mathfrak{R}_{Ric>0}\left(M_{m,n,l}\right)$ and $\mathfrak{M}_{Ric>0}\left(M_{m,n,l}\right)$ have infinitely many path components, answering Problem \ref{Problem: (moduli) space of positive Ricci curvature metrics of M_{m,n,l}} for certain manifolds $M_{m,n,l}$. Our main results are as follows. First, we provide a coarse classification (Theorem \ref{Theorem: coarse classify}), showing that the simply connected 7-manifolds $M_{m,n,l}$ naturally fall into five distinct families $\mathcal{M}_{1}\sim \mathcal{M}_{5}$ based on their cohomology and characteristic classes. We then achieve a complete diffeomorphism classification for the families $\mathcal{M}_{3}$ (Theorem \ref{Theorem: Classify mfds in M_3}) and $\mathcal{M}_{2}$ (Theorem \ref{Theorem: Classify mfds in M_2}). For the most intricate family $\mathcal{M}_{1}$, which is our primary focus, we obtain a partial classification (Theorem \ref{Theorem: Partial classification of mfds in M_1}). Finally, leveraging this classification, we answer Problem \ref{Problem: (moduli) space of positive Ricci curvature metrics of M_{m,n,l}} in the affirmative for infinitely many such manifolds (Theorem \ref{Theorem: mfd whose space & moduli space of Ric+ metrics have infinitely many components}).

\begin{thm}\label{Theorem: coarse classify}
	Consider the following $5$ families of manifolds:
	\begin{compactenum}
		\item $\mathcal{M}_{1}=\left\{M_{m,n,l}\middle|m,n,l\in\mathbb{Z};n,l\neq0;\gcd(m,n,l)=1\right\}$,
		\item $\mathcal{M}_{2}=\left\{M_{m,n,0}\middle|m,n\in\mathbb{Z};n\neq0;\gcd(m,n)=1\right\}$,
		\item $\mathcal{M}_{3}=\left\{M_{m,0,l}\middle|m,l\in\mathbb{Z};m,l\neq0;\gcd(m,l)=1\right\}$,
		\item $\mathcal{M}_{4}=\left\{M_{m,0,0}\middle|m=\pm1\right\}$,
		\item $\mathcal{M}_{5}=\left\{M_{0,0,l}\middle|l=\pm1\right\}$.
	\end{compactenum}
	If $M\in\mathcal{M}_{i}$, $\overline{M}\in\mathcal{M}_{j}$ and $1\leqslant i<j\leqslant5$, then $M$ and $\overline{M}$ are not homeomorphic.
\end{thm}

\begin{rmk}
	It can be shown that $M_{m,n,l}$ and $M_{-m,-n,-l}$ are orientation-reversing diffeomorphic. Hence the family $\mathcal{M}_{4}$ contains only one diffeomorphism class if orientation is not considered, and the same is also true for $\mathcal{M}_{5}$. It suffices to focus on classification of manifolds in $\mathcal{M}_{1}$, $\mathcal{M}_{2}$ and $\mathcal{M}_{3}$.
\end{rmk}



\begin{thm}[Complete classification of manifolds in $\mathcal{M}_{3}$]\label{Theorem: Classify mfds in M_3}
	Assume we have nonzero integers $m$, $l$, $\overline{m}$, $\overline{l}$ such that $\gcd(m,l)=\gcd\left(\overline{m},\overline{l}\right)=1$. Then $M_{m,0,l}$ and $M_{\overline{m},0,\overline{l}}$ are diffeomorphic if and only if one of the following is true:
		\begin{compactenum}
			\item $|m|,\ |l|\geqslant3$ and $\left(\overline{m},\overline{l}\right)=(m,l)$;
			\item $|m|\leqslant2$, $|l|\geqslant3$ and $\overline{m}=\pm m$, $\overline{l}=l$;
			\item $|m|\geqslant3$, $|l|\leqslant2$ and $\overline{m}=m$, $\overline{l}=\pm l$;
			\item $|m|,\ |l|\leqslant2$ and $\overline{m}=\pm m$, $\overline{l}=\pm l$.
		\end{compactenum}
\end{thm}

\begin{thm}[Complete classification of manifolds in $\mathcal{M}_{2}$]\label{Theorem: Classify mfds in M_2}
	Assume we have integers $m$, $\overline{m}$ and nonzero integers $n$, $\overline{n}$ such that $\gcd(m,n)=\gcd\left(\overline{m},\overline{n}\right)=1$. Then $M_{m,n,0}$ and $M_{\overline{m},\overline{n},0}$ are diffeomorphic if and only if
		\begin{equation}
			\left\{
			\begin{alignedat}{4}
				\overline{n}&=&\pm n, &\\
				\overline{m}&\equiv&m\quad &\mathrm{mod}\ 2^{\widetilde{\lambda}_{2}(n)}n^{2},
			\end{alignedat}
			\right.
		\end{equation}
		where
		\begin{equation}
			\widetilde{\lambda}_{2}(n)=\left\{
			\begin{alignedat}{6}
				0,\ & n & \equiv & 2,6\quad & \mathrm{mod}\ & 8,\\
				1,\ & n & \equiv & 1,7\quad & \mathrm{mod}\ & 8,\\
				2,\ & n & \equiv & 0,3,4,5\quad & \mathrm{mod}\ & 8.
			\end{alignedat}
			\right.
		\end{equation}
\end{thm}

\begin{thm}[Partial classification of manifolds in $\mathcal{M}_{1}$]\label{Theorem: Partial classification of mfds in M_1}
	Assume we have nonzero integers $n$, $l$, $\overline{n}$, $\overline{l}$.
	\begin{compactenum}
		\item Suppose $\gcd(n,l)=\gcd\left(\overline{n},\overline{l}\right)=1$. Manifolds $M_{0,n,l}$ and $M_{0,\overline{n},\overline{l}}$ are diffeomorphic if and only if 
		\begin{equation}
			\left\{
			\begin{alignedat}{2}
				\overline{n}&=\pm n,\\
				\overline{l}&=l.
			\end{alignedat}
			\right.
		\end{equation}
		\item Suppose $m$ is a non-zero integer such that $\gcd(m,n,l)=\gcd\left(\overline{n},\overline{l}\right)=1$ and $(m,n,l)\equiv\left(\overline{m},\overline{n},\overline{l}\right)\equiv(0,1,0)\ \mathrm{mod}\ 2$. Manifolds $M_{m,n,l}$ and $M_{0,\overline{n},\overline{l}}$ are diffeomorphic if and only if 
		\begin{equation}
			\left\{
				\begin{alignedat}{2}
					\overline{l}&=l,\\
					\overline{n}&=\pm n,\\
					{m}&\equiv 0\ \mathrm{mod}\ \frac{2^{\lambda_{2}(n)}7^{\lambda_{7}(n)}n^{2}}{3^{\mu_{3}(n)}},
				\end{alignedat}
			\right.
		\end{equation}
		where
		$$\mu_{3}(n)=\begin{cases}
			0, &n\equiv1,2\mod3,\\
			1, &n\equiv0\mod3.
		\end{cases}$$
		\item Suppose $m$ and $\overline{m}$ are nonzero integers such that $\gcd(m,n)=\gcd\left(\overline{m},\overline{n}\right)=1$ and $(m,n,l)\equiv\left(\overline{m},\overline{n},\overline{l}\right)\equiv(0,1,0)\ \mathrm{mod}\ 2$. Manifolds $M_{m,n,l}$ and $M_{\overline{m},\overline{n},\overline{l}}$ are diffeomorphic if and only if 
		\begin{equation}
			\left\{
				\begin{alignedat}{2}
					\overline{l}&=l,\\
					\overline{n}&=\pm n,\\
					\overline{m}&\equiv m\ \mathrm{mod}\ \frac{2^{\lambda_{2}(n)}7^{\lambda_{7}(n)}n^{2}}{3^{\mu_{3}(n)}}.
				\end{alignedat}
			\right.
		\end{equation}
	\end{compactenum}
\end{thm}

\begin{thm}\label{Theorem: mfd whose space & moduli space of Ric+ metrics have infinitely many components}
	There are integers $m$, $n$, $l$ with $\gcd(m,n,l)=1$, $n$, $l\neq0$ and $(m,n,l)\equiv(0,1,0)\ \mathrm{mod}\ 2$, such that $\mathfrak{R}_{Ric>0}\left(M_{m,n,l}\right)$ and $\mathfrak{M}_{Ric>0}\left(M_{m,n,l}\right)$ both have infinitely many path components.
\end{thm}

Although $M_{m,n,l}$ and $M_{m,n}$, $N_{m,n}$ are constructed similar, and modified surgery theory is the common technique applied to classify these manifolds, the classification of $M_{m,n,l}$ presents a substantial new challenge compared to the cases of $M_{m,n}$ and $N_{m,n}$. The reason is as follows. The $s$-invariants of manifolds derived from modified surgery depend on choices of ordered bases of its second cohomology and can be viewed as polynomial functions of even degrees on the ordered bases. The second cohomology groups of manifolds $M_{m,n}$ and $N_{m,n}$ studied in \cite{KS88,Kreck1991SomeNH,Kreck1991SomeNHCorrection} are $\mathbb{Z}$, whose generators differ by a sign. Hence $s$-invariants for $M_{m,n}$ and $N_{m,n}$ are authentic diffeomorphism invariants. However, the second cohomology groups of manifolds $M_{m,n,l}$ we concern here are $\mathbb{Z}^{2}$, different ordered bases of second cohoology group are differed by an element in $GL(2,\mathbb{Z})$ and might lead different value of $s$-invariants. While differentiating $M_{m,n,l}$ we shall consider all plausible ordered bases of $H^{2}\left(M_{m,n,l}\right)$, which transform under the action of the full group $GL(2,\mathbb{Z})$. Taming this action and extracting genuine diffeomorphism classification requires a delicate and novel analysis, which is undertaken in Section \ref{Section: 7 Partial classification of manifolds in M_1}.

The paper is organized as follows. In Section \ref{Section: 3 A first classification} we first list the result of cohomology and characteristic classes of $M_{m,n,l}$, classify these manifolds into five families $\mathcal{M}_{i}(1\leqslant i\leqslant5)$ and prove Theorem \ref{Theorem: coarse classify}. Then we give more detailed results for manifolds in $\mathcal{M}_{1}$. In Section \ref{Section: 4 Computation of bordism groups} we compute necessary bordism groups and derive the $s$-invariant of spin manifolds in $\mathcal{M}_{1}$. In Section \ref{Section: 7 Partial classification of manifolds in M_1} we compute the explicit expressions of $s$-invariants for certain spin manifolds in $\mathcal{M}_{1}$ and give a partial classification, proving Theorem \ref{Theorem: Partial classification of mfds in M_1} for spin manifolds. In Section \ref{Section: 8 Application to positive Ricci curvature} we apply the classification to differential geometry and show that there are spin manifolds in $\mathcal{M}_{1}$ whose spaces and moduli spaces of positive Ricci curvature metrics have infinitely many components, proving Theorem \ref{Theorem: mfd whose space & moduli space of Ric+ metrics have infinitely many components}. 
Throughout this paper we mainly focus on the spin manifolds in $\mathcal{M}_{1}$ and do not expand computation for other manifolds.

\begin{ackn}
	The author thanks Professor Wilderich Tuschmann and Professor Yi Jiang sincerely, whose insightful comments and discussions inspired the initial direction of this research. The author is also deeply grateful to Professor Yang Su and Professor Yi Jiang for their invaluable guidance and support throughout this work. The author's research was supported NSFC 11801298 and NSFC 12471069.
\end{ackn}
	
	\section{A first classification of $M_{m,n,l}$}\label{Section: 3 A first classification}

In this section we use cohomology ring and characteristic classes to classify simply connected $M_{m,n,l}$ into five families. In Section \ref{Section: Classification into five families} we fix the orientation convention and list the results of cohomology ring and characteristic classes of $M_{m,n,l}$. Their proofs are basic algebraic topology and are omitted. Then we apply these results and classify simply connected manifolds $M_{m,n,l}$ into five families by comparing fourth Betti numbers and first rational Pontryagin classes, proving Theorem \ref{Theorem: coarse classify}. In Section \ref{Section: coh ring & char classes of mfds in M_1} we discuss manifolds in $\mathcal{M}_{1}$ more carefully. We specify bases of second cohomology groups and generators of fourth cohomology groups. Then we compare cohomology ring structures and further classify some manifolds in $\mathcal{M}_{1}$.

\subsection{Classification into five families $\mathcal{M}_{1}\sim\mathcal{M}_{5}$}\label{Section: Classification into five families}

Recall that $N=\left(\mathbb{C}P^{1}\times\mathbb{C}P^{2}\right)\#\mathbb{C}P^{3}$ and $H^{2}(N)=\mathbb{Z}\left\{\alpha,\beta,\gamma\right\}\cong\mathbb{Z}^{3}$. Here $\alpha$, $\beta$ and $\gamma$ correspond to the Euler classes of Hopf bundles over $\mathbb{C}P^{1}$, $\mathbb{C}P^{2}$ and $\mathbb{C}P^{3}$ respectively. Let $\xi=\xi_{m,n,l}$ denote the complex line bundle over $N$ with first Chern class $m\alpha+n\beta+l\gamma$. Let $M=M_{m,n,l}$ and $V=V_{m,n,l}$ denote the total spaces of the unit disc bundle and circle bundle associated to $\xi_{m,n,l}$ respectively. Let $M_{m,n,l}\xrightarrow{p}N$, $V_{m,n,l}\xrightarrow{\pi}N$ denote the bundle projections and let $M_{m,n,l}\xrightarrow{\iota}V_{m,n,l}$ denote the boundary inclusion. Then $\pi$ is a homotopy equivalence and $p=\pi\circ\iota$.
These manifolds are oriented as follows.
\begin{compactenum}
	\item We orient $\mathbb{C}P^{n}$ so that the hyperplane $\left\{\left[z_{0}:\cdots:z_{n}\right]\in\mathbb{C}P^{n}\middle|z_{n}=0\right\}\cong\mathbb{C}P^{n-1}$ is Poincar\'{e} dual to the Euler class of Hopf bundle over $\mathbb{C}P^{n}$. 
	\item Given two oriented manifolds $M_{1}$ and $M_{2}$, their product manifold $M_{1}\times M_{2}$ is canonically oriented. If $M_{1}$ and $M_{2}$ have the same dimension, then their connected sum $M_{1}\#M_{2}$ is also canonically oriented. In particular $N=\left(\mathbb{C}P^{1}\times\mathbb{C}P^{2}\right)\#\mathbb{C}P^{3}$ is oriented.
	\item If $M$ is an oriented manifold and $\xi:F\to E\to M$ is an oriented fiber bundle over $M$ with oriented fiber manifold $F$, then the orientations of $\xi$, $F$ and $M$ determine an orientation on the total space manifold $E$. In particular $M_{m,n,l}$ and $V_{m,n,l}$ are oriented.
\end{compactenum}

\begin{lem}\label{Lemma: E 1cnt}
	$M_{m,n,l}$ is simply connected if and only if $\mathrm{gcd}(m,n,l)=1$. In this case $H_{2}\left(M_{m,n,l}\right)\cong\mathbb{Z}^{2}$ is torsion-free.
\end{lem}

The following are results of cohomology rings and characteristic classes.

\begin{prop}\label{Proposition: 4th cohomology gps of E}
	Suppose $\mathrm{gcd}(m,n,l)=1$. Then $H^{2}\left(M_{m,n,l}\right)\cong\mathbb{Z}^{2}$ and  $H^{4}\left(M_{m,n,l}\right)$ is generated by pairwise products of elements in $H^{2}\left(M_{m,n,l}\right)$. More precisely:
	\begin{compactenum}
		\item when $n,l\neq0$, $H^{4}\left(M_{m,n,l}\right)$ is a finite abelian group and $b_{4}=0$;
		\item when $n\neq0$ and $l=0$, we can find $z\in H^{2}(N)$ such that $\left(p^{*}z,p^{*}\gamma\right)$ is a basis of $H^{2}\left(M_{m,0,l}\right)$ and
		$$H^{4}\left(M_{m,n,0}\right)=\mathbb{Z}\big/n^{2}\left\{\left(p^{*}z\right)^{2}\right\}\oplus\mathbb{Z}\left\{\left(p^{*}\gamma\right)^{2}\right\},\ b_{4}\left(M_{m,n,0}\right)=1;$$
		\item when $m,l\neq0$ and $n=0$, we can find $z\in H^{2}(N)$ such that $\left(p^{*}\beta,p^{*}z\right)$ is a basis of $H^{2}\left(M_{m,0,l}\right)$ and
		$$	H^{4}\left(M_{m,0,l}\right)=\mathbb{Z}\left\{\left(p^{*}\beta\right)^{2}\right\}\oplus\mathbb{Z}\big/|m|\left\{\left(p^{*}\beta\right)\left(p^{*}z\right)\right\}\oplus\mathbb{Z}\big/|l|\left\{\left(p^{*}z\right)^{2}\right\},\ b_{4}\left(M_{m,0,l}\right)=1;$$
		\item when $n=l=0$ and $m=\pm1$ we have
		\begin{align*}
			H^{2}\left(M_{m,0,0}\right)&=\mathbb{Z}\left\{p^{*}\beta,p^{*}\gamma\right\},\\
			H^{4}\left(M_{m,0,0}\right)&=\mathbb{Z}\left\{\left(p^{*}\beta\right)^{2},\left(p^{*}\gamma\right)^{2}\right\}\cong\mathbb{Z}^{2},\ b_{4}\left(M_{m,0,0}\right)=2;
		\end{align*}
		\item when $m=n=0$ and $l=\pm1$ we have
		\begin{align*}
			H^{2}\left(M_{0,0,l}\right)&=\mathbb{Z}\left\{p^{*}\beta,p^{*}\alpha\right\},\\
			H^{4}\left(M_{0,0,l}\right)&=\mathbb{Z}\left\{\left(p^{*}\alpha\right)\left(p^{*}\beta\right),\left(p^{*}\beta\right)^{2}\right\}\cong\mathbb{Z}^{2}, b_{4}\left(M_{0,0,l}\right)=2.
		\end{align*}
	\end{compactenum}
\end{prop}

\begin{prop}\label{Proposition: char classes of E}
	The characteristic classes $w_{2}$ and $p_{1}$ of $M_{m,n,l}$ and $V_{m,n,l}$ are given as follows.
	\begin{compactenum}
		\item $w_{2}\left(V_{m,n,l}\right)=\pi^{*}\rho_{2}(m\alpha+(n+1)\beta+l\gamma)$, $p_{1}\left(V_{m,n,l}\right)=\pi^{*}\left(2mn\alpha\beta+\left(n^{2}+3\right)\beta^{2}+\left(l^{2}+4\right)\gamma^{2}\right)$. $V_{m,n,l}$ is spin if and only if $(m,n,l)\equiv(0,1,0)\ \mathrm{mod}\ 2$. 
		\item $w_{2}\left(M_{m,n,l}\right)=p^{*}\rho_{2}(\beta)$, $p_{1}\left(M_{m,n,l}\right)=p^{*}\left(3\beta^{2}+4\gamma^{2}\right)$. $M_{m,n,l}$ is spin if and only if $(m,n,l)\equiv(0,1,0)\ \mathrm{mod}\ 2$. 
		
		In terms of description in Proposition \ref{Proposition: 4th cohomology gps of E} these characteristic classes also have the following expressions:
		\begin{compactenum}
			\item[(1)] when $n,l\neq0$ we have $p_{1}^{\mathbb{Q}}\left(M_{m,n,l}\right)=0$;
			\item[(2)] when $n\neq0$ and $l=0$ we have $w_{2}\left(M_{m,n,0}\right)=m\rho_{2}\left(p^{*}z\right)$, $p_{1}\left(M_{m,n,0}\right)=3m^{2}\left(p^{*}z\right)^{2}+4\left(p^{*}\gamma\right)^{2}$, $p_{1}^{\mathbb{Q}}\left(M_{m,n,0}\right)=4\left(p^{*}\gamma\right)^{2}$;
			\item[(3)] when $m,l\neq0$ and $n=0$ we have $p_{1}\left(M_{m,0,l}\right)=3\left(p^{*}\beta\right)^2+4m^{2}\left(p^{*}z\right)^{2}$, $p_{1}^{\mathbb{Q}}\left(M_{m,0,l}\right)=3\left(p^{*}\beta\right)^2$;
			\item[(4)] when $n=l=0$ and $m=\pm1$ we have $p_{1}\left(M_{m,0,0}\right)=3\left(p^{*}\beta\right)^{2}+4\left(p^{*}\gamma\right)^{2}$;
			\item[(5)] when $m=n=0$ and $l=\pm1$ we have $p_{1}\left(M_{0,0,l}\right)=3\left(p^{*}\beta\right)^{2}$.
		\end{compactenum}
	\end{compactenum}
\end{prop}

The following lemma is helpful to reduce the classification when orientation is not mainly concerned.
\begin{lem}\label{Lemma: E_{m,n,l}=E_{-m,-n,-l}}
	The manifolds $M_{m,n,l}$ and $M_{-m,-n,-l}$ are orientation-reversing diffeomorphic.
\end{lem}
Hence $\mathcal{M}_{4}$ contains only one diffeomorphism class, and the same conclusion is also true for the family $\mathcal{M}_{5}$. 

We also have the following formula of signature of $V_{m,n,l}$.
\begin{prop}\label{Proposition: signature of (D,E)}
	The signature of $\left(V_{m,n,l},M_{m,n,l}\right)$ is 
	\begin{equation}\label{Equation: signature of (D,E)}
		\sigma\left(V_{m,n,l},M_{m,n,l}\right)
		=\sigma\left(\begin{pmatrix}
			0 & n & 0 \\ n & m & 0 \\ 0 & 0 & l
		\end{pmatrix}\right)=\mathrm{sgn}(l).
	\end{equation}		
\end{prop}

With these conclusions above, we are ready to prove Theorem \ref{Theorem: coarse classify}.

\begin{pf}[of Theorem \ref{Theorem: coarse classify}]
	By Proposition \ref{Proposition: 4th cohomology gps of E} we see that $b_{4}$ can be $0$, $1$ and $2$ depending on different values of $m$, $n$ and $l$. Now we apply Proposition \ref{Proposition: char classes of E} and use the first rational Pontryagin class to distinguish fanilies $\mathcal{M}_{2}$ and $\mathcal{M}_{3}$, fanilies $\mathcal{M}_{4}$ and $\mathcal{M}_{5}$. In the following we write $M=M_{m,n,l}$ and $V=V_{m,n,l}$ for simplicity. 
	
	When $n\neq0$ and $l=0$ we have $p_{1}^{\mathbb{Q}}(M)=4p^{*}\left(\gamma^{2}\right)$, which has divisibility $4$ in the integral lattice $\mathbb{Z}\left\{p^{*}\left(\gamma^{2}\right)\right\}$. When $m\neq0$, $n=0$ and $l\neq0$ we have $p_{1}^{\mathbb{Q}}(M)=3p^{*}\left(\beta^{2}\right)$, which has divisibility $3$ in the integral lattice $\mathbb{Z}\left\{p^{*}\left(\beta^{2}\right)\right\}$. Since the rational Pontryagin class is a homeomorhism invariant, the calculation above implies that if $M\in\mathcal{M}_{2}$ and $\overline{M}\in\mathcal{M}_{3}$, then $M$ and $\overline{M}$ are not homeomorphic.
	
	When $m\pm1$ and $n=l=0$ we have $p_{1}^{\mathbb{Q}}(M)=3p^{*}\left(\beta^{2}\right)+4p^{*}\left(\gamma^{2}\right)$, which is primitive in the lattice $\mathbb{Z}\left\{p^{*}\left(\beta^{2}\right),p^{*}\left(\gamma^{2}\right)\right\}$. When $m=n=0$ and $l=\pm1$ we have $p_{1}^{\mathbb{Q}}(M)=3p^{*}\left(\beta^{2}\right)$, which is not primitive in the lattice $\mathbb{Z}\left\{p^{*}\left(\alpha\beta\right),p^{*}\left(\beta^{2}\right)\right\}$ and has divisibility 3. Since the rational Pontryagin class is a homeomorhism invariant, the calculation above implies that if $M\in\mathcal{M}_{4}$ and $\overline{M}\in\mathcal{M}_{5}$, then $M$ and $\overline{M}$ are not homeomorphic.
\end{pf}

\subsection{More information of manifolds in $\mathcal{M}_{1}$}\label{Section: coh ring & char classes of mfds in M_1}

In the remaining of this section we further analyze $H^{4}(M)$ and the ring structure on $H^{*}(M)$ for manifolds in $\mathcal{M}_{1}$. We will equip $H^{2}(M)$ with bases and equip $H^{4}(M)$ with generators, and we intend to further express the generators of $H^{4}(M)$ as pairwise products of the chosed bases of $H^{2}(M)$ if possible. Then we express the characteristic classes in terms of these bases and generators. In Section \ref{Section: Classification into five families} the result is not so detailed as the families $\mathcal{M}_{i}$, $2\leqslant i\leqslant5$, since whether $m=0$ or not vastly affects the calculation.

First we compute the exponent of $H^{4}(M)$, give bases of $H^{2}(M)$ and generators of $H^{4}(M)$. The \textbf{exponent} of a finite abelian group is the minimal positive integer that annihilates all elements. We begin with some necessary notations. 
\begin{compactenum}
	\item The least common divisor of integers $P$ and $Q$ are denoted by $\mathrm{lcm}(P,Q)$ and is assumed to be positive.
	\item When $\gcd(n,l)=1$, by B\'{e}zout's theorem there are integers $p$ and $q$ such that $pn+ql=1$.
	\item When $m\neq0$ we denote $d=\gcd(m,n)$ and write $m=dm_{1}$, $n=dn_{1}$. Since $\gcd(m,n,l)=1$ we have $\gcd(d,l)=1$, and by B\'{e}zout's theorem there are integers $u$, $v$, $u'$, $v'$, $s$ and $t$ such that $um_{1}+vn_{1}=1$, $u'm_{1}+v'n_{1}=1$ and $sd+tl=1$.
\end{compactenum}

\begin{prop}\label{Proposition: M1 mfds, H^4, coh ring & char classes}
	Let $m\in\mathbb{Z}$, $n,l\in\mathbb{Z}\backslash\{0\}$ such that $\gcd(m,n,l)=1$.
	\begin{compactenum}
		\item When $m=0$, order and exponent of $H^{4}(M)$ are $n^{2}|l|$ and $|nl|$ respectively. We can find $z=-q\beta+p\gamma\in H^{2}(N)$ such that 
		\begin{compactenum}
			\item $\left(p^{*}z,p^{*}\alpha\right)$ is a basis of $H^{2}(M)$;
			\item $H^{4}(M)=\mathbb{Z}\big/|n|\left\{\left(p^{*}\alpha \right)\left(p^{*} z\right)\right\}\oplus\mathbb{Z}\big/|n|l\left\{\left(p^{*}z\right)^{2}\right\}\cong\mathbb{Z}\big/|n|\oplus\mathbb{Z}\big/|n|l$;
			\item $w_{2}(M)=\rho_{2}\left(lp^{*}z\right)$;
			\item  $p_{1}(M)=\left(3l^{2}+4n^{2}\right)\left(p^{*}z\right)^{2}$.
		\end{compactenum} 
		\item When $m\neq0$, order and exponent of $H^{4}(M)$ are $n^{2}|l|$ and $d\cdot\mathrm{lcm}\left(n_{1}^{2},l\right)$ respectively. We can find
		$f=-v\alpha+u\beta$, $g=-m_{1}t\alpha-n_{1}t\beta+s\gamma\in H^{2}(N)$ and $\omega=m_{1}\alpha\beta+n_{1}\beta^{2}$, $\rho=-v'\alpha\beta+u'\beta^{2}\in H^{4}(N)$
		such that 
		\begin{compactenum}
			\item $\left(p^{*}f,p^{*}g\right)$ is a basis of $H^{2}(M)$;
			\item $H^{4}(M)=\mathbb{Z}\big/d\left\{p^{*}\omega\right\}\oplus\mathbb{Z}\big/n_{1}^{2}\left\{p^{*}\rho\right\}\oplus\mathbb{Z}\big/|l|\left\{\left(p^{*}\gamma\right)^{2}\right\}$;
			\item  $\left(p^{*}\right)^{2}$, $\left(p^{*}f\right)\left(p^{*}g\right)$ and $\left(p^{*}g\right)^{2}$ also generate $H^{4}(M)$;
			\item $w_{2}(M)=\rho_{2}\left(m_{1}p^{*}f-lvp^{*}g\right)$;
			\item  $p_{1}(M)=3m_{1}^{2}\left(p^{*}f\right)^{2}-6vm_{1}l\left(p^{*}f\right)\left(p^{*}g\right)+\left(3v^{2}l^{2}+4d^{2}\right)\left(p^{*}g\right)^{2}
			=3v'p^{*}\omega+3m_{1}p^{*}\rho+4\left(p^{*}\gamma\right)^{2}$.
		\end{compactenum} 
		
		\item When $m\neq0$ and $\gcd(n,l)=1$, order and exponent of $H^{4}(M)$ are $n^{2}|l|$ and $dn_{1}^{2}|l|$ respectively. We can find
			$z=-q\beta+p\gamma\in H^{2}(N)$ and 
			$\omega=m_{1}\alpha\beta+n_{1}\beta^{2}$, $\rho=-v'\alpha\beta+u'\beta^{2}\in H^{4}(N)$
		such that 
		\begin{compactenum}
			\item $\left(p^{*}z,p^{*}\alpha\right)$ is a basis of $H^{2}(M)$;
			\item $H^{4}(M)=\mathbb{Z}\big/d\left\{p^{*}\omega\right\}\oplus\mathbb{Z}\big/dn_{1}^{2}\left\{p^{*}\rho\right\}\oplus\mathbb{Z}\big/|l|\left\{\left(p^{*}\gamma\right)^{2}\right\}$;
			\item $\left(p^{*}z\right)^{2}$, $\left(p^{*}z\right)\left(p^{*}\alpha\right)$ also generate $H^{4}(M)$ and $\left(p^{*}\alpha\right)^{2}=0$;
			\item $w_{2}(M)=\rho_{2}\left(-lp^{*}z-pm_{1}dp^{*}\alpha\right)$;
			\item $p_{1}(M)=\left(3l^{2}+4n_{1}^{2}d^{2}\right)\left(p^{*}z\right)^{2}+\left(6pm_{1}dl-8qm_{1}n_{1}d^{2}\right)\left(p^{*}z\right)\left(p^{*}\alpha\right)=3v'p^{*}\omega+3m_{1}p^{*}\rho+4\left(p^{*}\gamma\right)^{2}$.
		\end{compactenum} 
		\item When $m\neq0$ and $d=1$, order and exponent of $H^{4}(M)$ are $n^{2}|l|$ and $\mathrm{lcm}\left(n^{2},l\right)$ respectively. We can find $z=-v\alpha+u\beta\in H^{2}(N)$ such that
		\begin{compactenum}
			\item $\left(p^{*}z,p^{*}\gamma\right)$ is a basis of $H^{2}(M)$;
			\item $H^{4}(M)=\mathbb{Z}\big/n^{2}\left\{\left(p^{*}z\right)^{2}\right\}\oplus\mathbb{Z}\big/|l|\left\{\left(p^{*}\gamma\right)^{2}\right\}$;
			\item $w_{2}(M)=\rho_{2}\left(mp^{*}z-vlp^{*}\gamma\right)$;
			\item $p_{1}(M)=3m^{2}\left(p^{*}z\right)^{2}+4\left(p^{*}\gamma\right)^{2}$.
		\end{compactenum} 
	\end{compactenum}
\end{prop}

\begin{rmk}
	Following Proposition \ref{Proposition: M1 mfds, H^4, coh ring & char classes}, Statement 2 it is routine to compute the transformation relation between two sets of generators of $H^{4}(M)$:
	$$\begin{pmatrix}
		\left(p^{*}f\right)^{2} & \left(p^{*}f\right)\left(p^{*}g\right) & \left(p^{*}g\right)^{2}
	\end{pmatrix}=\begin{pmatrix}
		p^{*}\omega & p^{*}\rho & \left(p^{*}\gamma\right)^{2}
	\end{pmatrix}\begin{pmatrix}
		-2uvu'+u^{2}v' & vu'tn_{1}-ut & u't^{2}m_{1}n_{1}+t^{2}n_{1}\\
		uvn_{1}+u & vtn_{1}^{2} & -t^{2}m_{1}n_{1}^{2}\\
		0 & 0 & s^{2}
	\end{pmatrix}.$$
	Since $u'=u+\lambda n_{1}$, $v'=v-\lambda m_{1}$ for some $\lambda\in\mathbb{Z}$ and $p^{*}\omega$ has order $d$, we assume that $\lambda$ is divided by $d$ from now on and obtain 
	\begin{eqnarray}\label{Equation: generator transf when m neq 0}
		\begin{pmatrix}
			\left(p^{*}f\right)^{2} & \left(p^{*}f\right)\left(p^{*}g\right) & \left(p^{*}g\right)^{2}
		\end{pmatrix}=\begin{pmatrix}
			p^{*}\omega & p^{*}\rho & \left(p^{*}\gamma\right)^{2}
		\end{pmatrix}\begin{pmatrix}
			-u^{2}v & -u^{2}tm_{1} & ut^{2}m_{1}n_{1}+t^{2}n_{1}\\
			uvn_{1}+u & vtn_{1}^{2} & -t^{2}m_{1}n_{1}^{2}\\
			0 & 0 & s^{2}
		\end{pmatrix}.
	\end{eqnarray}
\end{rmk}

\begin{rmk}\label{Remark: M1 mfd coh ring, m neq0, n & l coprime, assume d divides u'}
	Following Proposition \ref{Proposition: M1 mfds, H^4, coh ring & char classes}, Statement 3 it is routine to compute the transformation relation between two sets of generators of $H^{4}(M)$:
	\begin{align*}
		\left\{
		\begin{alignedat}{2}
			\left(p^{*}z\right)^{2}&=q^{2}v'p^{*}\omega+q^{2}m_{1}p^{*}\rho+p^{2}\left(p^{*}\gamma\right)^{2},\\
			\left(p^{*}z\right)\left(p^{*}\alpha\right)&=-qu'p^{*}\omega+qn_{1}p^{*}\rho.
		\end{alignedat}
		\right.
	\end{align*}
	When $d$ and $n_{1}$ are also coprime, from B\'{e}zout's identities $u'm_{1}+v'n_{1}=1$ and $pdn_{1}+ql=1$ we may assume that $u'$ is divided by $d$ so that
	\begin{align}
		\left\{
		\begin{alignedat}{2}
			\left(p^{*}z\right)^{2}&=q^{2}v'p^{*}\omega+q^{2}m_{1}p^{*}\rho+p^{2}\left(p^{*}\gamma\right)^{2},\\
			\left(p^{*}z\right)\left(p^{*}\alpha\right)&=qn_{1}p^{*}\rho.
		\end{alignedat}
		\right.
	\end{align}
	Recall that $p^{*}\omega$ has order $d$, $p^{*}\rho$ has order $dn_{1}^{2}$ and $\left(p^{*}\gamma\right)^{2}$ has order $|l|$. From B\'{e}zout's identies and the assumption that $d$ divides $u'$ we see that $q$ and $v'$ are both coprime to $d$, $q$ is coprime to $n_{1}$ and $d$, and $p$ is coprime to $l$. Hence $\left(p^{*}z\right)^{2}$ has order $dn_{1}^{2}|l|$ and $\left(p^{*}z\right)\left(p^{*}\alpha\right)$ has order $d\left|n_{1}\right|$.
\end{rmk}

\begin{rmk}
	When $m\neq0$ the computation is vastly different from that of the case $m=0$. Hence comparison and classification are divided into the following $3$ parts:
	\begin{compactenum}
		\item Determine when $M_{0,n,l}$ and $M_{0,\overline{n},\overline{l}}$ are diffeomorphic.
		\item Determine when $M_{m,n,l}$ and $M_{0,\overline{n},\overline{l}}$ are diffeomorphic. Here $m\neq0$.
		\item Determine when $M_{m,n,l}$ and $M_{\overline{m},\overline{n},\overline{l}}$ are diffeomorphic. Here $m\neq0$, $\overline{m}\neq0$.
	\end{compactenum}
\end{rmk}

After comparing fourth cohomology groups and cup product structures, we obtain the following proposition:

\begin{prop}\label{Proposition: M1 mfds, first classification comparing H^4}
	Let $n$, $l$, $\overline{n}$, and $\overline{l}$ be non-zero integers.
	\begin{compactenum}
		\item Suppose $\gcd(n,l)=\gcd\left(\overline{n},\overline{l}\right)=1$. Then $M=M_{0,n,l}$ and $\overline{M}=M_{0,\overline{n},\overline{l}}$ are diffeomorphic only if $\overline{n}=\pm n$ and $\overline{l}=\pm l$.
		\item Suppose $m$ is another non-zero integer and $\gcd(m,n,l)=\gcd\left(\overline{n},\overline{l}\right)=1$. Then $M=M_{m,n,l}$ and $\overline{M}=M_{0,\overline{n},\overline{l}}$ are diffeomorphic only if $n$ divides $m$, $\overline{n}=\pm n$ and $\overline{l}=\pm l$.
		\item Suppose $m$ and $\overline{m}$ are non-zero integers and $\gcd(m,n)=\gcd\left(\overline{m},\overline{n}\right)=1$. Then $M=M_{m,n,l}$ and $\overline{M}=M_{\overline{m},\overline{n},\overline{l}}$ are diffeomorphic only if $\overline{n}=\pm n$ and $\overline{l}=\pm l$.
	\end{compactenum}
\end{prop}

Now we prove these propositions.

\begin{pf}[of Proposition \ref{Proposition: M1 mfds, H^4, coh ring & char classes}]
	Recall the matrix representation of homomorphism $H^{2}(N)\xrightarrow{e\cup\cdot}H^{4}(N)$:
	$$e\cup\begin{pmatrix}
		\alpha & \beta & \gamma
	\end{pmatrix}=\begin{pmatrix}
		\alpha\beta & \beta^{2} & \gamma^{2}
	\end{pmatrix}\begin{pmatrix}
		n & m & 0\\
		0 & n & 0\\
		0 & 0 & l
	\end{pmatrix}.$$
	When $m=0$ the matrix is diagonal and the cokernel is easier to analyze. While the computation becomes more complicated when $m\neq0$. Hence we treat these two cases separately.
	
	\noindent\textbf{Proof of Statement 1}
	
	We start with the case that $m=0$. It is straightforward to see that
	\begin{align*}
		H^{4}(M)&=\mathbb{Z}\big/|n|\left\{\left(p^{*}\alpha\right)\left(p^{*}\beta\right)\right\}\oplus\mathbb{Z}\big/|n|\left\{\left(p^{*}\beta\right)^{2}\right\}\oplus\mathbb{Z}\big/|l|\left\{\left(p^{*}\gamma\right)^{2}\right\}\\
		&\cong\mathbb{Z}\big/|n|\oplus\mathbb{Z}\big/|n|\oplus\mathbb{Z}\big/|l|,
	\end{align*}
	and we immediately obtain that $H^{4}(M)$ has order $n^{2}|l|$ and exponent $|nl|$.
	
	Next we equip with $H^{2}(M)$ a basis and express generators of $H^{4}(M)$. By assumption $\gcd(n,l)=1$. By B\'{e}zout's theorem there are integers $p$ and $q$ such that $pn+ql=1$. Then we have
	$$\begin{pmatrix}
		\alpha & e & z
	\end{pmatrix}=\begin{pmatrix}
		\alpha & \beta & \gamma
	\end{pmatrix}\begin{pmatrix}
		1 & 0 & 0\\0 & n & -q\\0 & l & p
	\end{pmatrix},\quad\begin{pmatrix}
		\alpha & \beta & \gamma
	\end{pmatrix}=\begin{pmatrix}
		\alpha & e & z
	\end{pmatrix}\begin{pmatrix}
		1 & 0 & 0\\0 & p & q\\0 & -l & n
	\end{pmatrix},$$
	and we obtain
	\begin{align*}
		H^{2}(M)&=\mathbb{Z}\left\{p^{*}z,p^{*}\alpha\right\}\cong\mathbb{Z}^{2};\\
		p^{*}\begin{pmatrix}
			\alpha\beta & \beta^{2} & \gamma^{2}
		\end{pmatrix}&=\begin{pmatrix}
			p^{*}(\alpha z) & p^{*}\left(z^{2}\right)
		\end{pmatrix}\begin{pmatrix}
			-l & 0 & 0\\0 & l^{2} & n^{2}
		\end{pmatrix},\\
		H^{4}(M)&=\mathbb{Z}\big/|n|\left\{\left(p^{*}\alpha\right)\left(p^{*} z\right)\right\}\oplus\mathbb{Z}\big/|nl|\left\{\left(p^{*}z\right)^{2}\right\}\cong\mathbb{Z}\big/|n|\oplus\mathbb{Z}\big/|nl|.
	\end{align*}
	The value of $p^{*}z\in H^{2}(M)$ does not depend on choices of integers $p$ and $q$ which satisty B\'{e}zout's identity. 
	
	It is routine to compute characteristic classes, and Proposition \ref{Proposition: M1 mfds, H^4, coh ring & char classes}, Statement 1 is proved.
	
	\noindent\textbf{Proof of Statement 2}
	
	Now suppose $m\neq0$. Let $d$, $m_{1}$, $n_{1}$ and $u$, $v$, $u'$, $v'$, $s$, $t$ be introduced as before. Set
	$$\begin{pmatrix}
		\omega & \rho
	\end{pmatrix}=\begin{pmatrix}
		\alpha\beta & \beta^{2}
	\end{pmatrix}\begin{pmatrix}
		m_{1} & -v'\\ n_{1} & u'
	\end{pmatrix},\quad\begin{pmatrix}
		\alpha\beta & \beta^{2}
	\end{pmatrix}=\begin{pmatrix}
		\omega & \rho
	\end{pmatrix}\begin{pmatrix}
		u' & v' \\ -n_{1} & m_{1}
	\end{pmatrix},$$
	and we have
	\begin{align*}
		\mathrm{im}\left(H^{2}(N)\xrightarrow{e\cup\cdot}H^{4}(N)\right)&=\mathbb{Z}\left\{n_{1}d\alpha\beta,m_{1}d\alpha\beta+n_{1}d\beta^{2},l\gamma^{2}\right\}\\
		&=\mathbb{Z}\left\{n_{1}^{2}d\rho,d\omega,l\gamma^{2}\right\},\\
		H^{4}(M)&=\mathbb{Z}\big/n_{1}^{2}d\left\{p^{*}\rho\right\}\oplus\mathbb{Z}\big/d\left\{p^{*}\omega\right\}\oplus\mathbb{Z}\big/|l|\left\{\left(p^{*}\gamma\right)^{2}\right\}\\
		&\cong\mathbb{Z}\big/n_{1}^{2}d\oplus\mathbb{Z}\big/d\oplus\mathbb{Z}\big/|l|
	\end{align*}
	We immediately obtain that $H^{4}(M)$ has order $n^{2}|l|$ and exponent $d\cdot\mathrm{lcm}\left(n_{1}^{2},l\right)$. 
	
	In this case $p^{*}\rho$ depends on choices of integers $u'$ and $v'$ that satisfy B\'{e}zout's identity unless $d=1$. Suppose $u''m_{1}+v''n_{1}=1$ for another pair of integers $\left(u'',v''\right)$, and there is a unique integer $\lambda$ such that $u''=u'+\lambda n_{1}$ and $v''=v'-\lambda m_{1}$. Denote $\rho_{\lambda}=-v''\alpha\beta+u''\beta^{2}$, and we have $\rho_{\lambda}=\rho+\lambda\omega$. Then $p^{*}\rho_{\lambda}=p^{*}\rho+\lambda p^{*}\omega$ depends on the value of $\lambda\ \mathrm{mod}\ d$ unless $d=1$ so that $p^{*}\omega=p^{*}(e\beta)=0$.
	
	Now we equip $H^{2}(M)$ with a basis when $d$ is an arbitrary positive integer. We have
	\begin{align*}
		\begin{pmatrix}
			e & f & g
		\end{pmatrix}=\begin{pmatrix}
			\alpha & \beta & \gamma
		\end{pmatrix}\begin{pmatrix}
			m_{1}d & -v & -m_{1}t \\
			n_{1}d & u & -n_{1}t \\
			l & 0 & s
		\end{pmatrix},\quad
		\begin{pmatrix}
			\alpha & \beta & \gamma
		\end{pmatrix}=\begin{pmatrix}
			e & f & g
		\end{pmatrix}\begin{pmatrix}
			su & sv & t \\
			-n_{1} & m_{1} & 0 \\
			-ul & -vl & d
		\end{pmatrix};
	\end{align*}
	and we obtain
	$$H^{2}(M)=\mathbb{Z}\left\{p^{*}f,p^{*}g\right\}.$$
	Values of $p^{*}f$ and $p^{*}g$ in $H^{4}(M)$ depend on choices of $u$, $v$, $s$ and $t$. Pairs $(u,v)$ and $\left(u',v'\right)$ can be different.
	
	It is routine to compute characteristic classes,
	and Proposition \ref{Proposition: M1 mfds, H^4, coh ring & char classes}, Statement 2 is justified.
	
	\noindent\textbf{Proof of Statement 3}
	
	Now suppose $m\neq0$ and $\gcd(n,l)=1$. Let $d$, $m_{1}$, $n_{1}$, $u'$, $v'$, $p$ and $q$ be introduced as before. It follows from the proof of Statement 2 that we still have
	\begin{align*}
		H^{4}(M)&=\mathbb{Z}\big/n_{1}^{2}d\left\{p^{*}\rho\right\}\oplus\mathbb{Z}\big/d\left\{p^{*}\omega\right\}\oplus\mathbb{Z}\big/|l|\left\{\left(p^{*}\gamma\right)^{2}\right\}\\
		&\cong\mathbb{Z}\big/n_{1}^{2}d\oplus\mathbb{Z}\big/d\oplus\mathbb{Z}\big/|l|.
	\end{align*}
	By assumption $n$ and $l$ are coprime, hence we equip $H^{2}(M)$ with another basis. We have
	\begin{align*}
		\begin{pmatrix}
			\alpha & e & z
		\end{pmatrix}=\begin{pmatrix}
			\alpha & \beta & \gamma
		\end{pmatrix}\begin{pmatrix}
			1 & dm_{1} & 0 \\
			0 & dn_{1} & -q \\
			0 & l & p
		\end{pmatrix},\quad
		\begin{pmatrix}
			\alpha & \beta & \gamma
		\end{pmatrix}=\begin{pmatrix}
			\alpha & e & z
		\end{pmatrix}\begin{pmatrix}
			1 & -pdm_{1} & -qdm_{1} \\
			0 & p & q \\
			0 & -l & dn_{1}
		\end{pmatrix};
	\end{align*}
	and we obtain 
	$H^{2}(M)=\mathbb{Z}\left\{p^{*}z,p^{*}\alpha\right\}\cong\mathbb{Z}^{2}$.
	Here the value of $p^{*}z$ depends on choice of integers $p$ and $q$. Suppose $p'n+q'l=1$ for another pair of integers $\left(p',q'\right)$, and there is a unique integer $\lambda$ such that $p'=p+\lambda l$ and $q'=q-\lambda n$. Denote $z_{\lambda}=-q'\beta+p'\gamma$, and we have $z_{\lambda}=z+\lambda\left(e-m\alpha\right)$, $p^{*}z_{\lambda}=p^{*}z-\lambda mp^{*}\alpha$.
	
	It is routine to compute characteristic classes,
	and Proposition \ref{Proposition: M1 mfds, H^4, coh ring & char classes}, Statement 3 is justified.
	
	\noindent\textbf{Proof of Statement 4}
	
	Finally we assume $m\neq0$ and $d=1$. In this case $m_{1}=m$ and $n_{1}=n$. Now we have
	\begin{align*}
		H^{4}(M)&=\mathbb{Z}\big/n^{2}\left\{p^{*}\rho\right\}\oplus\mathbb{Z}\big/|l|\left\{\left(p^{*}\gamma\right)^{2}\right\}\\
		&\cong\mathbb{Z}\big/n^{2}\oplus\mathbb{Z}\big/|l|,
	\end{align*}
	and we immediately obtain that $H^{4}(M)$ has order $n^{2}|l|$ and exponent $\mathrm{lcm}\left(n^{2},l\right)$.
	
	Since $d=1$, in the relation $sd+tl=1$ we take $s=1$, $t=0$ and obtain \begin{align*}
		\begin{pmatrix}
			e & z & \gamma
		\end{pmatrix}=\begin{pmatrix}
			\alpha & \beta & \gamma
		\end{pmatrix}\begin{pmatrix}
			m & -v & 0 \\
			n & u & 0 \\
			l & 0 & 1
		\end{pmatrix},\quad
		\begin{pmatrix}
			\alpha & \beta & \gamma
		\end{pmatrix}=\begin{pmatrix}
			e & z & \gamma
		\end{pmatrix}\begin{pmatrix}
			u & v & 0 \\
			-n & m & 0 \\
			-ul & -vl & 1
		\end{pmatrix}.
	\end{align*}
	We also have $H^{2}(M)=\mathbb{Z}\left\{p^{*}z,p^{*}\gamma\right\}\cong\mathbb{Z}^{2}$.
	The value of $p^{*}z$ depends on choices of integers $u$ and $v$. Suppose $u''m+v''n=1$ for another pair of integers $\left(u'',v''\right)$, and there is a unique integer $\lambda$ such that $u''=u+\lambda n$ and $v''=v-\lambda m$. Denote $z_{\lambda}=-v''\alpha+u''\beta$, and we have $z_{\lambda}=z+\lambda(e-l\gamma)$, $p^{*}z_{\lambda}=p^{*}z-\lambda lp^{*}\gamma$. 
	
	In this case
	\begin{align*}
		p^{*}\rho&=-v'\left(p^{*}\alpha\right)\left(p^{*}\beta\right)+u'\left(p^{*}\beta\right)^{2}\\
		&=-v'\left(-np^{*}z-ulp^{*}\gamma\right)\left(mp^{*}z-vlp^{*}\gamma\right)+u'\left(mp^{*}z-vlp^{*}\gamma\right)^{2}\\
		&=\left(u'm^{2}+v'mn\right)\left(p^{*}z\right)^{2}=m\left(p^{*}z\right)^{2}
	\end{align*}
	since $z\gamma=0$, $\left(p^{*}\gamma\right)^{2}$ has order $|l|$ and $u'm+v'n=1$.
	Now we have
	\begin{align*}
		H^{4}(M)&=\mathbb{Z}\big/n^{2}\left\{m\left(p^{*}z\right)^{2}\right\}\oplus\mathbb{Z}\big/|l|\left\{\left(p^{*}\gamma\right)^{2}\right\}\\
		&=\mathbb{Z}\big/n^{2}\left\{\left(p^{*}z\right)^{2}\right\}\oplus\mathbb{Z}\big/|l|\left\{\left(p^{*}\gamma\right)^{2}\right\}\\
		&\cong\mathbb{Z}\big/n^{2}\oplus\mathbb{Z}\big/|l|
	\end{align*}
	as $m$ and $n$ are coprime.
	It is routine to compute the characteristic classes,
	and Proposition \ref{Proposition: M1 mfds, H^4, coh ring & char classes}, Statement 4 is justified. Now we complete the proof of Proposition \ref{Proposition: M1 mfds, H^4, coh ring & char classes}.
\end{pf}

\begin{pf}[of Proposition \ref{Proposition: M1 mfds, first classification comparing H^4}]
	We give a brief sketch before providing details. 
	\begin{compactenum}
		\item To prove the first statement it suffices to compare exponents and orders of fourth cohomology groups.
		\item To prove the second statement, we first apply the ordered basis of $H^{2}\left(M_{m,n,l}\right)$ given in the second part of Proposition \ref{Proposition: M1 mfds, H^4, coh ring & char classes}, and we obtain that $\gcd\left(d,n_{1}\right)=\gcd\left(n_{1},l\right)=1$ according to the cup product structure $H^{2}(M)\times H^{2}(M)\xrightarrow{\cdot\cup\cdot}H^{4}(M)$. Next we apply the ordered basis of $H^{2}\left(M_{m,n,l}\right)$ given in the third part of Proposition \ref{Proposition: M1 mfds, H^4, coh ring & char classes} and deduce that $n_{1}=\pm1$ again from the comparison of cup product structure. 
		\item To prove the third statement, first we compare exponents and orders of fourth cohomology groups to deduce that $\gcd\left(n^{2},l\right)=\gcd\left(\overline{n}^{2},\overline{l}\right)$, then we apply the ordered bases of second cohomology groups given in the third part of Proposition \ref{Proposition: M1 mfds, H^4, coh ring & char classes} and deduce that $\overline{n}=\pm n$ and $\overline{l}=\pm l$.
	\end{compactenum} 
	
	\noindent\textbf{Proof of Statement 1}
	
	We begin with the first statement. Let $n$, $l$ and $\overline{n}$, $\overline{l}$ be given as above. The orders of $H^{4}(M)$ and $H^{4}\left(\overline{M}\right)$ are $n^{2}|l|$ and $\overline{n}^{2}\left|\overline{l}\right|$ respectively. The exponents of $H^{4}(M)$ and $H^{4}\left(\overline{M}\right)$ are $|nl|$ and $\left|\overline{n}\overline{l}\right|$ respectively. Hence we obtain $\overline{n}=\pm n$, $\overline{l}=\pm l$ and complete proof of Statement 1.
	
	\noindent\textbf{Proof of Statement 2}
	
	Next we consider the second statement. Let $m$, $n$, $l$ and $\overline{n}$, $\overline{l}$ be given as above. From the comparison of orders and exponents of fourth cohomology groups we obtain $d^{2}n_{1}^{2}|l|=\overline{n}^{2}\left|\overline{l}\right|$ and $d\cdot\mathrm{lcm}\left(n_{1}^{2},l\right)=\left|\overline{n}\overline{l}\right|$. We denote $\delta=\gcd\left(n_{1}^{2},l\right)\geqslant1$ and write $n_{1}^{2}=\delta n_{2}$, $l=\delta l_{1}$ with $n_{2}>0$ and $\mathrm{sgn}(l)=\mathrm{sgn}\left(l_{2}\right)$. Let the exponents and orders of fourth cohomology groups be equal, and we obtain
	\begin{align*}
		\left|\overline{n}\right|=d\delta,\quad \left|\overline{l}\right|=n_{2}\left|l_{1}\right|.
	\end{align*}
	Since $\gcd\left(d,l\right)=\gcd\left(\overline{n},\overline{l}\right)=1$, integers $d$, $\delta$, $n_{2}$ and $l_{1}$ are pairwise coprime. Since $n_{1}^{2}=\delta n_{2}$, there are coprime integers $\delta_{1}$ and $n_{3}$ such that $\delta=\delta_{1}^{2}$, $n_{2}=n_{3}^{2}$. We assume $\delta_{1}>0$ and $\mathrm{sgn}\left(n_{3}\right)=\mathrm{sgn}\left(n_{1}\right)$ so that $n_{1}=\delta_{1}n_{3}$. Now the tuples of integers can be rewritten as
	\begin{align*}
		(m,n,l)=\left(dm_{1},d\delta_{1}n_{3},\delta_{1}^{2}l_{1}\right),\quad
		\left(0,\overline{n},\overline{l}\right)=\left(0,\pm d\delta_{1}^{2},\pm n_{3}^{2}l_{1}\right).
	\end{align*}
	
	Then we $\delta_{1}=1$ by comparing cohomology rings. Let ${M}\xrightarrow{\varphi} \overline{M}$ be a diffeomorphism. Then $H^{*}\left(\overline{M}\right)\xrightarrow{\varphi^{*}}H^{*}\left({M}\right)$ is an isomorphism. Set
	$$\varphi^{*}\begin{pmatrix}
		\overline{p}^{*}\overline{z} & \overline{p}^{*}\alpha
	\end{pmatrix}=\begin{pmatrix}
		p^{*}f & p^{*}g
	\end{pmatrix}\begin{pmatrix}
		A & C \\ B & D
	\end{pmatrix},\ \begin{pmatrix}
	A & C \\ B & D
	\end{pmatrix}\in GL(2,\mathbb{Z}),$$
	and $\delta_{1}=1$ will follow from $\varphi^{*}\left(\left(\overline{p}^{*}\alpha\right)^{2}\right)=0$. Note that $\varphi^{*}\left(\left(\overline{p}^{*}\alpha\right)^{2}\right)=0$ since $\alpha^{2}=0$. Meanwhile
	\begin{align*}
		\varphi^{*}\left(\left(\overline{p}^{*}\alpha\right)^{2}\right)&=\left(\varphi^{*}\left(\overline{p}^{*}\alpha\right)\right)^{2}\\
		&=C^{2}\left(p^{*}f\right)^{2}+2CD\left(p^{*}f\right)\left(p^{*}g\right)+D^{2}\left(p^{*}g\right)^{2}\\
		&=\left(-C^{2}uv-2CDu^{2}tm_{1}+D^{2}\left(ut^{2}m_{1}n_{1}+t^{2}n_{1}\right)\right)p^{*}\omega\\
		&\ +\left(C^{2}\left(uvn_{1}+u\right)+2CDvtn_{1}^{2}-D^{2}t^{2}m_{1}n_{1}^{2}\right)p^{*}\rho
		+D^{2}s^{2}\left(p^{*}\gamma\right)^{2}.
	\end{align*}
	Recall that $p^{*}\omega$ has order $d$, $p^{*}\rho$ has order $dn_{1}^{2}=d\delta_{1}^{2}n_{3}^{2}$ and $\left(p^{*}\gamma\right)^{2}$ has order $|l|=\delta_{1}^{2}\left|l_{1}\right|$. Since $\psi^{*}\left(\left(\overline{p}^{*}\alpha\right)^{2}\right)=0$ we must have
	$$
		\left\{
			\begin{alignedat}{2}
				C^{2}\left(uvn_{1}+u\right)+2CDvtn_{1}^{2}-D^{2}t^{2}m_{1}n_{1}^{2}&\equiv0\ \mathrm{mod}\ dn_{3}^{2}\delta_{1}^{2},\\
				D^{2}s^{2}&\equiv0\ \mathrm{mod}\ \delta_{1}^{2}l_{1},
			\end{alignedat}
		\right.
	$$
	and the $\delta_{1}$-primary part reads
	$$
		\left\{
			\begin{alignedat}{2}
				C^{2}\left(uv\delta_{1}n_{3}+u\right)&\equiv0\ \mathrm{mod}\ \delta_{1}^{2},\\
				D^{2}s^{2}&\equiv0\ \mathrm{mod}\ \delta_{1}^{2}.
			\end{alignedat}
		\right.
	$$

	We begin with the second congruence equation. It follows from $1=sd+tl=sd+t\delta_{1}^{2}l_{1}$ that $s$ and $\delta_{1}$ are coprime. Hence the second equation reads $D^{2}\equiv0\ \mathrm{mod}\ \delta_{1}^{2}$, namely $D\equiv0\ \mathrm{mod}\ \delta_{1}$. 
	Next we consider the first congruence equation and start with the modulus $\delta_{1}$, obtaining $C^{2}u\equiv0\ \mathrm{mod}\ \delta_{1}$. It follows from $1=um_{1}+vn_{1}=um_{1}+v\delta_{1}n_{3}$ that $u$ and $\delta_{1}$ are coprime, hence $C^{2}\equiv0\ \mathrm{mod}\ \delta_{1}$. Now we return to the first equation with modulus $\delta_{1}^{2}$, and we obtain $C^{2}u\equiv0\ \mathrm{mod}\ \delta_{1}^{2}$. Since $u$ and $\delta_{1}$ are coprime we have $C^{2}\equiv0\ \mathrm{mod}\ \delta_{1}^{2}$, namely $C\equiv0\ \mathrm{mod}\ \delta_{1}$.
	Now both $C$ and $D$ are divided by $\delta_{1}$. Since $AD-BC=\pm1$ we see $C$ and $D$ are coprime, and $\delta_{1}$ must be $1$. As a result, $n_{3}=n_{1}$, $l_{1}=l$, integers $d$, $n_{1}$ and $l$ are pairwise coprime, $\overline{n}=\pm d$ and $\overline{l}=\pm n_{1}^{2}l$. 
	
	It remains to show $n_{1}=\pm1$. Since $d$, $n_{1}$ and $l$ are pairwise coprime, we apply the ordered basis $\left(p^{*}z,p^{*}\alpha\right)$ of $H^{2}(M)$ given in Proposition \ref{Proposition: M1 mfds, H^4, coh ring & char classes}, Statement 3 and the assumption that $d$ divides $u'$ from Remark \ref{Remark: M1 mfd coh ring, m neq0, n & l coprime, assume d divides u'}. Set
	$$\varphi^{*}\begin{pmatrix}
		\overline{p}^{*}\overline{z} & \overline{p}^{*}\alpha
	\end{pmatrix}=\begin{pmatrix}
		p^{*}z & p^{*}\alpha
	\end{pmatrix}\begin{pmatrix}
		\widetilde{A} & \widetilde{C} \\ \widetilde{B} & \widetilde{D}
	\end{pmatrix},\ \begin{pmatrix}
	\widetilde{A} & \widetilde{C} \\ \widetilde{B} & \widetilde{D}
	\end{pmatrix}\in GL(2,\mathbb{Z}).$$
	We will deduce $n_{1}=\pm1$ from the analysis of $\varphi^{*}\left(\left(\overline{p}^{*}\alpha\right)^{2}\right)$ and $\varphi^{*}\left(\left(\overline{p}^{*}\overline{z}\right)\left(\overline{p}^{*}\alpha\right)\right)$. 
	
	It is straightforward to compute that
	$$
		\left\{
			\begin{alignedat}{2}
				\varphi^{*}\left(\left(\overline{p}^{*}\overline{z}\right)^{2}\right)&=\widetilde{A}^{2}\left(p^{*}z\right)^{2}+2\widetilde{A}\widetilde{B}\left(p^{*}z\right)\left(p^{*}\alpha\right),\\
				\varphi^{*}\left(\left(\overline{p}^{*}\overline{z}\right)\left(\overline{p}^{*}\alpha\right)\right)&=\widetilde{A}\widetilde{C}\left(p^{*}z\right)^{2}+\left(\widetilde{A}\widetilde{D}+\widetilde{B}\widetilde{C}\right)\left(p^{*}z\right)\left(p^{*}\alpha\right),\\
				\varphi^{*}\left(\left(\overline{p}^{*}\alpha\right)^{2}\right)&=\widetilde{C}^{2}\left(p^{*}z\right)^{2}+2\widetilde{C}\widetilde{D}\left(p^{*}z\right)\left(p^{*}\alpha\right).
			\end{alignedat}
		\right.
	$$
	Recall that $\left(\overline{p}^{*}\overline{z}\right)^{2}$ has order $dn_{1}^{2}|l|$, $\left(\overline{p}^{*}\overline{z}\right)\left(\overline{p}^{*}\alpha\right)$ has order $d$, $\left(\overline{p}^{*}\alpha\right)^{2}=0$ and when $u'$ is divided by $d$ we have
	\begin{align*}
		\left\{
		\begin{alignedat}{2}
			\left(p^{*}z\right)^{2}&=q^{2}v'p^{*}\omega+q^{2}m_{1}p^{*}\rho+p^{2}\left(p^{*}\gamma\right)^{2},\\
			\left(p^{*}z\right)\left(p^{*}\alpha\right)&=qn_{1}p^{*}\rho,\\
			\left(p^{*}\alpha\right)^{2}&=0,
		\end{alignedat}
		\right.
	\end{align*}
	and $\left(p^{*}z\right)^{2}$ has order $dn_{1}^{2}|l|$, $\left(p^{*}z\right)\left(p^{*}\alpha\right)$ has order $d\left|n_{1}\right|$.
	
	We begin with $\varphi^{*}\left(\left(\overline{p}^{*}\alpha\right)^{2}\right)$. Since $\left(\overline{p}^{*}\alpha\right)^{2}=0$, we must have
	$$
		\left\{
			\begin{alignedat}{2}
				\widetilde{C}^{2}&\equiv0\ \mathrm{mod}\ dn_{1}^{2}|l|,\\
				2\widetilde{C}\widetilde{D}&\equiv0\ \mathrm{mod}\ d\left|n_{1}\right|.
			\end{alignedat}
		\right.
	$$
	Since $d$, $n_{1}$ and $l$ are pairwise coprime, from the first equation we obtain that $\widetilde{C}^{2}\equiv0\ \mathrm{mod}\ d$, $\widetilde{C}\equiv0\ \mathrm{mod}\ \left|n_{1}\right|$ and $\widetilde{C}^{2}\equiv0\ \mathrm{mod}\ |l|$. Then it follows from $\widetilde{A}\widetilde{D}-\widetilde{B}\widetilde{C}=\pm1$ that both $\widetilde{A}$ and $\widetilde{D}$ are coprime to $d$, $n_{1}$ and $l$.
	
	Next we consider $\varphi^{*}\left(\left(\overline{p}^{*}\overline{z}\right)\left(\overline{p}^{*}\alpha\right)\right)$. Since $\left(\overline{p}^{*}\alpha\right)\left(\overline{p}^{*}\alpha\right)$ has order $d$, we must have
	$$
		\left\{
			\begin{alignedat}{2}
				\widetilde{A}\widetilde{C}d&\equiv0\ \mathrm{mod}\ dn_{1}^{2}|l|,\\
				\left(\widetilde{A}\widetilde{D}+\widetilde{B}\widetilde{C}\right)d&\equiv0\ \mathrm{mod}\ d\left|n_{1}\right|.
			\end{alignedat}
		\right.
	$$
	From the analysis of $\varphi^{*}\left(\left(\overline{p}^{*}\alpha\right)^{2}\right)$ we already obtain that $\widetilde{C}\equiv0\ \mathrm{mod}\ \left|n_{1}\right|$. Hence the second equation implies that $\widetilde{A}\widetilde{D}\equiv0\ \mathrm{mod}\ \left|n_{1}\right|$. While $\widetilde{A}$ and $\widetilde{D}$ are both coprime to $n_{1}$, we must have $n_{1}=\pm1$.
	
	As a result, $m=dm_{1}$, $n=dn_{1}=\pm d$ and $\overline{n}=\pm n$, $\overline{l}=\pm l$. 
	This completes the proof of Statement 2.
	
	\noindent\textbf{Proof of Statement 3}
	
	Finally we consider the last statement. Let $m$, $n$, $l$ and $\overline{m}$, $\overline{n}$, $\overline{l}$ be given as above. From the comparison of order and exponent of fourth cohomology groups we obtain $n^{2}|l|=\overline{n}^{2}\left|\overline{l}\right|$ and $\mathrm{lcm}\left(n^{2},l\right)=\mathrm{lcm}\left(\overline{n}^{2},\overline{l}\right)$. We denote $\delta=\gcd\left(n^{2},l\right)$, $n^{2}=\delta n_{2}$ and $l=\delta l_{1}$ as in the proof of Statement 2. Then $\gcd\left(n_{1},l_{1}\right)=1$. We also introduce the bar counterpart and obtain $\overline{\delta}=\delta$, $n_{2}\left|l_{1}\right|=\overline{n_{2}}\left|\overline{l_{1}}\right|$.
	
	We claim $\overline{n_{2}}=n_{2}$ and $\overline{l_{1}}=\pm l_{1}$. To prove this claim we compare cohomology rings. 
	Let $M\xrightarrow{\varphi}\overline{M}$ be a diffeomorphism and set $\psi=\varphi^{-1}$. Then $H^{*}\left(\overline{M}\right)\xrightarrow{\varphi^{*}}H^{*}(M)$ and $H^{*}\left(M\right)\xrightarrow{\psi^{*}}H^{*}\left(\overline{M}\right)$ are isomorphisms. Set
	\begin{align*}
		\varphi^{*}\begin{pmatrix}
			\overline{p}^{*}\overline{z} & \overline{p}^{*}\gamma
		\end{pmatrix}=\begin{pmatrix}
			p^{*}z & p^{*}\gamma
		\end{pmatrix}\begin{pmatrix}
			A & C \\ B & D
		\end{pmatrix},\quad \psi^{*}\begin{pmatrix}
			p^{*}z & p^{*}\gamma
		\end{pmatrix}=\begin{pmatrix}
			\overline{p}^{*}\overline{z} & \overline{p}^{*}\gamma
		\end{pmatrix}\begin{pmatrix}
			\varepsilon D & -\varepsilon C \\ -\varepsilon B & \varepsilon A
		\end{pmatrix}.
	\end{align*}
	Here $\begin{pmatrix}
		A & C \\ B & D
	\end{pmatrix}\in GL(2,\mathbb{Z})$ and $\varepsilon=AD-BC=\pm1$. Then it is straightforward to compute
	\begin{align*}
		\left\{
		\begin{alignedat}{2}
			\varphi^{*}\left(\left(\overline{p}^{*}\overline{z}\right)^{2}\right)&=A^{2}\left(p^{*}z\right)^{2}+B^{2}\left(p^{*}\gamma\right)^{2},\\
			\varphi^{*}\left(\left(\overline{p}^{*}\overline{z}\right)\left(\overline{p}^{*}\gamma\right)\right)&=AC\left(p^{*}z\right)^{2}+BD\left(p^{*}\gamma\right)^{2},\\
			\varphi^{*}\left(\left(\overline{p}^{*}\gamma\right)^{2}\right)&=C^{2}\left(p^{*}z\right)^{2}+D^{2}\left(p^{*}\gamma\right)^{2},
		\end{alignedat}
		\right.
	\end{align*}
	Recall that
	\begin{align*}
		H^{4}(M)&=\mathbb{Z}\big/\delta n_{2}\left\{\left(p^{*}z\right)^{2}\right\}\oplus\mathbb{Z}\big/\delta \left|l_{1}\right|\left\{\left(p^{*}\gamma\right)^{2}\right\},\\
		H^{4}\left(\overline{M}\right)&=\mathbb{Z}\big/\delta \overline{n_{2}}\left\{\left(\overline{p}^{*}\overline{z}\right)^{2}\right\}\oplus\mathbb{Z}\big/\delta\left|\overline{l_{1}}\right|\left\{\left(\overline{p}^{*}\gamma\right)^{2}\right\}.
	\end{align*}
	Hence $\delta\overline{n_{2}}\varphi^{*}\left(\overline{p}^{*}z\right)^{2}=\varphi^{*}\left(\overline{p}^{*}\overline{z}\right)\left(\left(\overline{p}^{*}\gamma\right)\right)=\delta\overline{l_{1}}\varphi^{*}\left(\overline{p}^{*}\gamma\right)^{2}=0$, and we have
	\begin{equation*}
		\left\{\begin{alignedat}{2}
			A^{2}\overline{n_{2}}\equiv0\ \mathrm{mod}\ n_{2},\quad &B^{2}\overline{l_{1}}\equiv0\ \mathrm{mod}\ \left|l_{1}\right|;\\
			AC\equiv0\ \mathrm{mod}\ \delta n_{2},\quad &BD\equiv0\ \mathrm{mod}\ \delta \left|l_{1}\right|;\\
			C^{2}\overline{n_{2}}\equiv0\ \mathrm{mod}\ n_{2},\quad &D^{2}\overline{l_{1}}\equiv0\ \mathrm{mod}\ \left|l_{1}\right|.
		\end{alignedat}\right.
	\end{equation*}
	It follows from $AD-BC=\pm1$ that $\gcd(A,C)=\gcd\left(A^{2},C^{2}\right)=1$ and $\gcd(B,D)=\gcd\left(B^{2},D^{2}\right)=1$. Hence
	\begin{equation*}
		\left\{\begin{alignedat}{2}
			\overline{n_{2}}\equiv0\ \mathrm{mod}\ n_{2},\quad &\overline{l_1}\equiv0\ \mathrm{mod}\ \left|l_1\right|;\\
			AC\equiv0\ \mathrm{mod}\ \delta n_{2},\quad &BD\equiv0\ \mathrm{mod}\ \delta \left|l_1\right|.
		\end{alignedat}\right.
	\end{equation*}
	We also consider the equations from action of $\psi$. By symmetry it suffices to interchange $\left(\overline{p}^{*}\overline{z},\overline{p}^{*}\alpha\right)$, $\left({p}^{*}{z},{p}^{*}\alpha\right)$ and replace $A$, $B$, $C$, $D$ by $\varepsilon D$, $-\varepsilon B$, $-\varepsilon C$ $\varepsilon A$. Then we obtain
	\begin{equation*}
		\left\{\begin{alignedat}{2}
			n_{2}\equiv0\ \mathrm{mod}\ \overline{n_{2}},\quad &l_1\equiv0\ \mathrm{mod}\ \left|\overline{l_1}\right|,\\
			CD\equiv0\ \mathrm{mod}\ \delta\overline{n_{2}},\quad &AB\equiv0\ \mathrm{mod}\ \delta\left|\overline{l_1}\right|.
		\end{alignedat}\right.
	\end{equation*}
	Hence $\overline{l_1}=\pm l_1$ and $\overline{n_{2}}=n_{2}$.
	As a result, $\overline{n}=\pm n$ and $\overline{l}=\pm l$. This completes the proof of Statement 3.
\end{pf}
	
	\section{The $s$-invariant of spin manifolds in $\mathcal{M}_{1}$}\label{Section: 4 Computation of bordism groups}

Now we apply modified surgery theory to deduce more delicate invariants and give a more complete classification. We focus on spin manifolds in $\mathcal{M}_{1}$. In Section \ref{Section: spin bordism gp} we determine the bordism group $\Omega^{Spin}_{8}\left(K_{2}\right)$ that is required in deduction of $s$-invariant, and we also list the result of $\Omega^{Spin}_{8}\left(K_{2};\mathrm{pr}_{1}^{*}\gamma^{1}\right)$. In Section \ref{Section:  s inv for E_1^+ mfds} we extend spin manifolds in $\mathcal{M}_{1}$ to a broader class of manifolds in view of modified surgery theory, namely the $\mathcal{E}_{1}^{+}$-manifolds, and derive their $s$-invariants.

We shall mention that the $s$-invariants of non-spin manifolds in $\mathcal{M}_{1}$ and of manifolds in $\mathcal{M}_{i}$, $2\leqslant i\leqslant5$ can be derived similarly. 

\subsection{The bordism groups }\label{Section: spin bordism gp}

In this section we study the bordism group 
$$\Omega^{Spin}_{8}\left(K_{2}\right)=\left\{[M;x,y]\middle| M\text{ is a closed spin $8$-manifold};x,y\in H^{2}(M)\right\}.$$
We will show that $\Omega^{Spin}_{8}\left(K_{2}\right)$ is both isomorphic to $\mathbb{Z}^{10}$, establish explicit isomorphisms and evaluate specific characteristic numbers on the generators. We also list our results of the group 
$$\Omega^{Spin}_{8}\left(K_{2};\mathrm{pr}_{1}^{*}\gamma^{1}\right)=\left\{[M,x,y]\middle|M\text{ is a closed oriented $8$-manifold};x,y\in H^{2}(M),\rho_{2}(x)=w_{2}(M)\right\}.$$

First we list some results on divisibility of characteristic numbers. They follow from basic algebraic topology and (generalized) Rokhlin's theorem and are helpful to establish homomorphisms from bordism groups to $\mathbb{Z}^{10}$.

\begin{lem}\label{Lemma: Char num are even}
	Let $M$ be a topological closed spin $8$-manifolds and $x,\ y\in H^{2}(M)$.
	Then $\left<x^{2}y^{2},[M]\right>\in\mathbb{Z}$ is even. In particular $\left<x^{4},[M]\right>\in\mathbb{Z}$ is also even.
\end{lem}

\begin{lem}\label{Lemma: Char num are divided by 24, 48}
	Let $M$ be a smooth closed spin $8$-manifold and $x, y\in H^{2}(M)$.
	Then $M$ admits oriented $4$-dimensional embedding submanifolds $M(x,x)$ and $M(x,y)$ which are Poincar\'{e} dual to $x^{2}$ and $xy$ respectively, and the following characteristic numbers of $M$ are integers:
		\begin{align*}
			\frac{\sigma(x,x)}{16}&:=\frac{1}{48}\left<x^{2}\left(p_{1}(M)-2x^{2}\right),[M]\right>,\\
			\frac{\sigma(x,y)-(x-y)^2}{8}&:=\frac{1}{24}\left<xy\left(p_{1}(M)-4\left(x^{2}+y^{2}\right)+6xy\right),[M]\right>.
		\end{align*}
\end{lem}

Now we introduce some $8$-dimensional closed manifolds and their distinguished second cohomology classes. Generators of the bordism groups that we concern can be represented by these manifolds.

\begin{compactenum}
	\item Let $Bott$ be the following manifold. Plumb $8$ copies of the total space unit disc bundle associated to the tangent bundle of $S^{4}$, and we obtain a smooth compact $8$-manifold $W\left(E_{8}\right)$ whose boundary represents a generator of the group of homotopy $7$-spheres $\Theta_{7}\cong\mathbb{Z}\big/28$. Form the boundary connected sum of $28$ copies of $W\left(E_{8}\right)$, and we obtain a smooth compact $8$-manifold whose boundary is now diffeomorphic to $S^{7}$ so that we attach an $8$-dimensional unit disc $D^{8}$ along the boundary via a standard diffeomorphism, obtaining the smooth closed manifold $Bott$. In particular $Bott$ is $3$-connected.
	\item Let $z_{n}\in H^{2}\left(\mathbb{C}P^{n}\right)$ be the Poincar\'{e} dual of the embedding hyperplane $\mathbb{C}P^{n-1}$. Let $V(2)$ be the hypersurface in $\mathbb{C}P^{5}$ of degree $2$ and let $\widehat{z_{5}}$ denote the restriction of $z_{5}$ to $V(2)$.
	\item Consider $\left(S^{2}\right)^{4}$, the product of $4$ copies of $S^{2}$. Let $x_{i}$ be the second cohomology class corresponding to the $i$-th component $S^{2}$ and denote $\Delta=x_{1}+x_{2}+x_{3}+x_{4}$.
	\item Consider $\left(\mathbb{C}P^{2}\right)^{2}$, the product of $2$ copies of $\mathbb{C}P^{2}$. Let $z_{2,i}$ be the second cohomology class corresponding to the $i$-th component $\mathbb{C}P^{2}$ and denote $\Delta'=z_{2,1}+z_{2,2}$.
	\item Given integers $j\geqslant i\geqslant1$, the Milnor hypersurface $H_{i,j}$ is defined as 
	 $$H_{i,j}:=\left\{\left(\left[z_{0}:\cdots:z_{i}\right],\left[w_{0}:\cdots:w_{j}\right]\right)\in\mathbb{C}P^{i}\times\mathbb{C}P^{j}\middle|\sum_{u=0}^{i}z_{u}w_{u}=0\right\}.$$
	Its Poincar\'{e} dual in $\mathbb{C}P^{i}\times\mathbb{C}P^{j}$ is $z_{i}+z_{j}$, and we denote by $\overline{z_{i}}, \overline{z_{j}}$ the restrictions of $z_{i},z_{j}\in H^{2}\left(\mathbb{C}P^{i}\times\mathbb{C}P^{j}\right)$ to $H_{i,j}$ respectively.
\end{compactenum}

Now we state our results of $\Omega^{Spin}_{8}\left(K_{2}\right)$ and $\Omega^{Spin}_{8}\left(K_{2};\mathrm{pr}_{1}^{*}\gamma^{1}\right)$.

\begin{thm}\label{Theorem: results of spin bordism gp}
	We have $\Omega^{Spin}_{8}\left(K_{2}\right)\cong\mathbb{Z}^{10}$. The homomorphism
	\begin{alignat*}{2}
		\Omega^{Spin}_{8}\left(K_{2}\right)&\to\mathbb{Z}^{10},\\
		[M,x,y]&\mapsto \left[\sigma(M),\widehat{A}(M),\frac{1}{2}x^{4},x^{3}y,\frac{1}{2}x^{2}y^{2},xy^{3},\frac{1}{2}y^{4},\frac{\sigma(x,x)}{16},\frac{\sigma(x,y)-(x-y)^2}{8},\frac{\sigma(y,y)}{16}\right]
	\end{alignat*}
	is an isomorphism. A basis of $\Omega^{Spin}_{8}\left(K_{2}\right)$ is
	\begin{align*}
		&\left[\mathbb{H}P^{2},0,0\right],\ \left[Bott,0,0\right],\\ &\left[V(2),\widehat{z_{5}},0\right],\ \left[\mathbb{C}P^{3}\times\mathbb{C}P^{1},z_{3},z_{1}\right],\ \left[V(2),\widehat{z_{5}},\widehat{z_{5}}\right],
		\left[\mathbb{C}P^{3}\times\mathbb{C}P^{1},z_{1},z_{3}\right],\ \left[V(2),0,\widehat{z_{5}}\right],\\
		&\left[\left(S^{2}\right)^{4},\Delta,0\right],\ \left[\left(S^{2}\right)^{4},x_{1}+x_{2},x_{3}+x_{4}\right],\ \left[\left(S^{2}\right)^{4},0,\Delta\right].
	\end{align*}
	 Evaluating certain characteristic numbers on this basis, we have Table \ref{Table: char nums of a basis for spin bordism gp}.
	\begin{table}[ht]
		\centering
		\caption{Certain characteristic numbers evaluating on the basis of $\Omega^{Spin}_{8}\left(K_{2}\right)$}
		\label{Table: char nums of a basis for spin bordism gp}
		\begin{tabular}{ccccccccccc}
			\toprule[1.5pt]
			$\Omega^{Spin}_{8}\left(K_{2}\right)$ & $\sigma$ & $p_{1}^{2}$ & $x^{2}p_{1}$ & $xyp_{1}$ & $y^{2}p_{1}$ & $x^{4}$ & $x^{3}y$ & $x^{2}y^{2}$ & $xy^{3}$ & $y^{4}$ \\
			\midrule[1pt]
			$\left[\mathbb{H}P^{2},0,0\right]$ & 1 & 4 & 0 & 0 & 0 & 0 & 0 & 0 & 0 & 0\\
			$\left[Bott,0,0\right]$ & 224 & 0 & 0 & 0 & 0 & 0 & 0 & 0 & 0 & 0\\
			$\left[V(2),\widehat{z_{5}},0\right]$ & 2 & 8 & 4 & 0 & 0 & 2 & 0 & 0 & 0 & 0\\
			$\left[\mathbb{C}P^{3}\times\mathbb{C}P^{1},z_{3},z_{1}\right]$ & 0 & 0 & 0 & 4 & 0 & 0 & 1 & 0 & 0 & 0\\
			$\left[V(2),\widehat{z_{5}},\widehat{z_{5}}\right]$ & 2 & 8 & 4 & 4 & 4 & 2 & 2 & 2 & 2 & 2\\
			$\left[\mathbb{C}P^{3}\times\mathbb{C}P^{1},z_{1},z_{3}\right]$ & 0 & 0 & 0 & 4 & 0 & 0 & 0 & 0 & 1 & 0\\
			$\left[V(2),0,\widehat{z_{5}}\right]$ & 2 & 8 & 0 & 0 & 4 & 0 & 0 & 0 & 0 & 2\\
			$\left[\left(S^{2}\right)^{4},\Delta,0\right]$ & 0 & 0 & 0 & 0 & 0 & 24 & 0 & 0 & 0 & 0\\
			$\left[\left(S^{2}\right)^{4},x_{1}+x_{2},x_{3}+x_{4}\right]$ & 0 & 0 & 0 & 0 & 0 & 0 & 0 & 4 & 0 & 0\\
			$\left[\left(S^{2}\right)^{4},0,\Delta\right]$ & 0 & 0 & 0 & 0 & 0 & 0 & 0 & 0 & 0 & 24\\
			\bottomrule[1.5pt]
		\end{tabular}
	\end{table}
\end{thm}

\begin{thm}\label{Theorem: results of nonspin bordism gp}
	We have $\Omega^{Spin}_{8}\left(K_{2};\mathrm{pr}_{1}^{*}\gamma_{1}\right)\cong\mathbb{Z}^{10}$. The homomorphism
	\begin{align*}
		\Omega^{Spin}_{8}\left(K_{2};\mathrm{pr}_{1}^{*}\gamma^{1}\right)&\to\mathbb{Z}^{10},\\
		[M,x,y]&\mapsto\left[\sigma(M),e^{\frac{x}{2}}\widehat{A}(M), x^{4},\frac{x^{3}y+x^{2}y^{2}}{2},\frac{x^{3}y-x^{2}y^{2}}{2},\frac{xy^{3}+y^{4}}{2},\frac{xy^{3}-y^{4}}{2},\right.\\
		&\quad \left.\frac{\sigma(x,x)-x^2}{8},\frac{\sigma(x+y,y)}{16},\frac{\sigma(y,y)-x^2}{8}\right]
	\end{align*}
	is an isomorphism. A basis of $\Omega^{Spin}_{8}\left(K_{2};\mathrm{pr}_{1}^{*}\gamma^{1}\right)$ is
	\begin{align*}
		&\left[\mathbb{H}P^{2},0,0\right],\ \left[Bott,0,0\right],\\
		&\left[\mathbb{C}P^{4},z_{4},0\right],\ \left[H_{2,3},\overline{z_{3}},\overline{z_{2}}\right],\ \left[H_{2,3},\overline{z_{3}},-\overline{z_{2}}\right],\ 
		\left[H_{1,4},\overline{z_{1}},\overline{z_{4}}\right],\ \left[H_{1,4},\overline{z_{1}},-\overline{z_{4}}\right],\\ &\left[\left(\mathbb{C}P^{2}\right)^{2},\Delta',0\right],\ \left[H_{1,4},\overline{z_{1}},-\overline{z_{1}}+2\overline{z_{4}}\right],\ \left[H_{2,3},\overline{z_{3}},\overline{z_{2}}+\overline{z_{3}}\right].
	\end{align*}
\end{thm}

\begin{rmk}
	Since $\mathbb{H}P^{2}$ and $Bott$ are $3$-connected, their second cohomology groups are trivial. Hence we denote the bordism classes $\left[\mathbb{H}P^{2},0,0\right]$ and $\left[Bott,0,0\right]$ by $\mathbb{H}P^{2}$ and $Bott$ respectively for convenience.
\end{rmk}

\begin{pf}[of Theorem \ref{Theorem: results of spin bordism gp}]
	Let us first sketch the proof. According to Lemmata \ref{Lemma: Char num are even} and \ref{Lemma: Char num are divided by 24, 48}, the characteristic numbers mentioned in the theorem are integers when evaluating on $8$-dimensional $\left(\mathcal{B},\mathcal{F}\right)$-manifolds, and they define homomorphisms from the bordism group to $\mathbb{Z}$.
	The proof can be summarized into three steps:
	\begin{compactenum}
		\item First apply Atiyah-Hirzebruch spectral sequence (AHSS for short). We obtain that $\Omega^{Spin}_{8}\left(K_{2}\right)$ is a finitely generated abelian group that has rank $10$ and contains no odd torsion. 
		\item In the filtration associated to AHSS there is a nontrivial extension problem. Hence we further apply Adams spectral sequence (ASS for short) to determine the isomorphism type of $\Omega^{Spin}_{8}\left(K_{2}\right)$. We have seen from AHSS that $\Omega^{Spin}_{8}\left(K_{2}\right)$ has no odd torsion, hence it suffices to consider the $2$-primary part. The result implies that $\Omega^{Spin}_{8}\left(K_{2}\right)$ is torsion-free and thus isomorphic to $\mathbb{Z}^{10}$.
		\item Finally we use characteristic numbers to construct an explicit isomorphism $\Omega^{Spin}_{8}\left(K_{2}\right)\to\mathbb{Z}^{10}$. 
	\end{compactenum}
	
	First we apply AHSS. The $E^{2}$-page concentrates in the first quadrant, as shown in Figure \ref{Figure: E^{2} page of AHSS}.
	\begin{figure}[h]
		\begin{center}
			\begin{tikzpicture}
				\matrix (m) [matrix of math nodes,
				nodes in empty cells,nodes={minimum width=5ex,
					minimum height=5ex,inner sep=3pt, outer sep=0pt},
				column sep=1ex,row sep=-1ex]{
					q & & & & & & & & & & & \\
					8 & \mathbb{Z}^{2} & 0 & & & & & & & & & \\
					7 & 0 & 0 & 0 & & & & & & & &\\
					6 & 0 & 0 & 0 & 0 & & & & & & &\\
					5 & 0 & 0 & 0 & 0 & 0 & & & & & &\\
					4 & \mathbb{Z} & 0 & \mathbb{Z}^{2} & 0 & \mathbb{Z}^{3} & 0 & & & & &\\
					3 & 0 & 0 & 0 & 0 & 0 & 0 & 0 & & & &\\
					2 & \mathbb{Z}\big/2 & 0 & \left(\mathbb{Z}\big/2\right)^{2} & 0 & \left(\mathbb{Z}\big/2\right)^{3} & 0 & \left(\mathbb{Z}\big/2\right)^{4} & 0 & & & \\
					1 & \mathbb{Z}\big/2 & 0 & \left(\mathbb{Z}\big/2\right)^{2} & 0 & \left(\mathbb{Z}\big/2\right)^{3} & 0 & \left(\mathbb{Z}\big/2\right)^{4} & 0 & \left(\mathbb{Z}\big/2\right)^{5} & & \\
					0 & \mathbb{Z} & 0 & \mathbb{Z}^{2} & 0 & \mathbb{Z}^{3} & 0 & \mathbb{Z}^{4} & 0 & \mathbb{Z}^{5} & 0 &\\
					\quad\strut & 0 & 1 & 2 & 3 & 4 & 5 & 6 & 7 & 8 & 9 & p\strut \\};
				\draw[-stealth,color=blue,line width=0.3mm] (m-9-10) -- (m-8-8)
				node[midway,above,yshift=-4pt]{\scriptsize \textcolor{blue}{$d^{2}_{8,1}$}};
				\draw[-stealth,color=blue,line width=0.3mm] (m-10-10) -- (m-9-8)
				node[midway,above,yshift=-4pt]{\scriptsize \textcolor{blue}{$d^{2}_{8,0}$}};;
				\draw[thick,<-] (m-1-1.east) -- (m-11-1.east) ;
				\draw[thick,->] (m-11-1.north) -- (m-11-12.north) ;
				\draw [opacity=.4,line width=6mm,line cap=round, color=yellow] 
				(m-2-2.center) -- (m-10-10.center);
				\draw[opacity=.4,line width=6mm,line cap=round, color=orange]
				(m-2-3.center) -- (m-10-11.center);
				\draw[opacity=.4,line width=6mm,line cap=round, color=orange]
				(m-3-2.center) -- (m-10-9.center);
			\end{tikzpicture}
		\end{center}
		\caption{The $E^{2}$-page of Atiyah-Hirzebruch spectral sequence for $\Omega^{Spin}_{8}\left(K_{2}\right)$}
		\label{Figure: E^{2} page of AHSS}
	\end{figure}
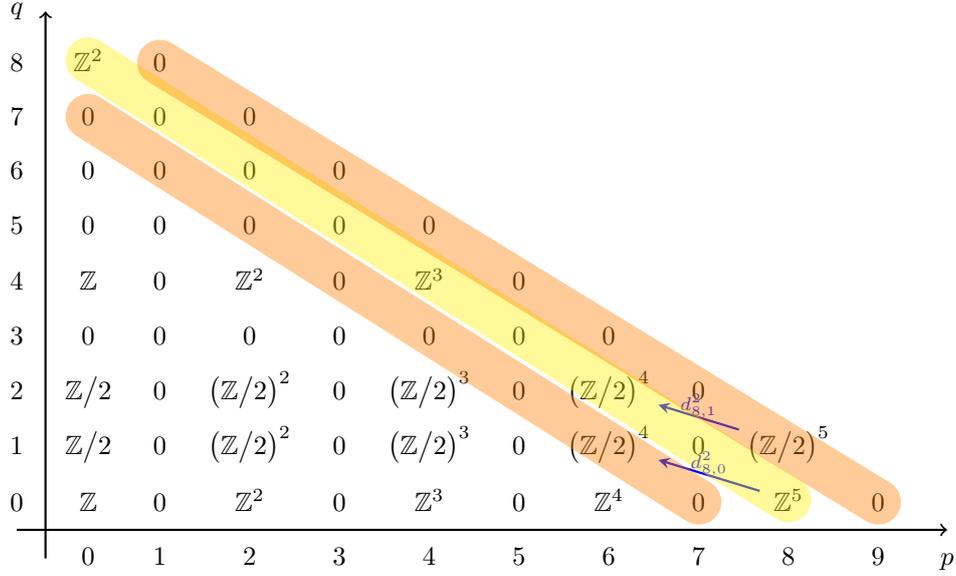
	It can be shown that $E^{3}_{p,q}=E^{\infty}_{p,q}$ when $p+q=8$, and $\bigoplus_{p+q=8}E^{\infty}_{p,q}=\mathbb{Z}^{10}\oplus\mathbb{Z}\big/2$. Hence $\Omega^{Spin}_{8}\left(K_{2}\right)$ is finitely generated and has rank $10$, while in the filtration associated to AHSS there is a non-trivial extension problem.
	
	Now we consider the extension problem. By Pontryagin-Thom construction the bordism group is identified with the homotopy group of corresponding Thom spectrum
	$$\Omega^{Spin}_{8}\left(K_{2}\right)=\pi_{8}\left(MSpin\wedge\left(K_{2}\right)_{+}\right).$$
	Since $\Omega^{Spin}_{8}\left(K_{2}\right)$ contains no odd torsion, it suffices to consider the $2$-local part and by ASS we have
	\begin{equation*}
		E^{s,t}_{2}=\mathrm{Ext}^{s,t}_{\mathcal{A}}\left(H^{*}\left(MSpin\wedge\left(K_{2}\right)_{+};\mathbb{Z}\big/2\right),\mathbb{Z}\big/2\right)\Longrightarrow\left(\Omega^{Spin}_{t-s}\left(K_{2}\right)\right)\sphat_{2}.
	\end{equation*}
	By K\"{u}nneth formula we have
	\begin{equation*}
		H^{*}\left(MSpin\wedge\left(K_{2}\right)_{+};\mathbb{Z}\big/2\right)\cong H^{*}\left(MSpin;\mathbb{Z}\big/2\right)\otimes_{\mathbb{Z}\big/2}H^{*}\left(K_{2};\mathbb{Z}\big/2\right);
	\end{equation*}
	here 
	\begin{compactenum}
		\item $H^{*}\left(K_{2};\mathbb{Z}\big/2\right)\cong\mathbb{Z}\big/2\left[\overline{x},\overline{y}\right]$, $\overline{x},\overline{y}\in H^{2}\left(K_{2};\mathbb{Z}\big/2\right)$;
		\item $H^{*}\left(MSpin;\mathbb{Z}\big/2\right)\cong\left(\mathcal{A}/\mathcal{A}(1)\right)U+\left(\mathcal{A}/\mathcal{A}(1)\right)w_{4}^{2}U+\mathfrak{M}$, $U\in H^{0}\left(MSpin;\mathbb{Z}\big/2\right)$, $w_{4}\in H^{4}\left(MSpin;\mathbb{Z}\big/2\right)$ $\mathfrak{M}\subset H^{\geqslant9}\left(MSpin;\mathbb{Z}\big/2\right)$, $\mathcal{A}$ is the mod $2$ Steenrod algebra and $\mathcal{A}(1)$ is the sub-algebra generated by $\mathrm{Sq}^{1}$ and $\mathrm{Sq}^{2}$
	\end{compactenum}
	To compute $E_{2}^{s,t}$ we need the $\mathcal{A}(1)$-module structure on $H^{*}\left(K_{2};\mathbb{Z}\big/2\right)$. From the action of Steenrod operations
	we obtain an isomorphism of $\mathcal{A}(1)$-modules
	\begin{equation}
		H^{*}\left(K_{2};\mathbb{Z}\big/2\right)\cong L_{1}\oplus\left(\bigoplus_{i\geqslant1}(i+1)\Sigma^{4(i-1)}L_{2}\right)\oplus\left(\bigoplus_{j\geqslant1}j\Sigma^{4j-2}L_{2}\right).
	\end{equation}
	Here $L_{1}$ is a copy of $\mathbb{Z}\big/2$ centered at degree $0$, $L_{2}$ is isomorphic to $\widetilde{H}^{*}\left(\mathbb{C}P^{2};\mathbb{Z}\big/2\right)$ and $\Sigma$ is the operator that shifts degree by $1$.
	Therefore, we have the isomorphism up to degree $8$
	\begin{align*}
		H^{*}\left(MSpin;\mathbb{Z}\big/2\right)\otimes_{\mathbb{Z}\big/2}H^{*}\left(K_{2};\mathbb{Z}\big/2\right)
		\cong\mathcal{A}\otimes_{\mathcal{A}(1)}\left(L_{1}\oplus2L_{2}\oplus\Sigma^{2}L_{2}\oplus3\Sigma^{4}L_{2}\oplus4\Sigma^{6}L_{2}\oplus\Sigma^{8}L_{1}\right).
	\end{align*}
	The Adams charts of $\mathrm{Ext}^{*,*}_{\mathcal{A}}\left(\mathcal{A}\otimes_{\mathcal{A}(1)}L_{i},\mathbb{Z}\big/2\right)\cong\mathrm{Ext}^{*,*}_{\mathcal{A}(1)}\left(L_{i},\mathbb{Z}\big/2\right)$ $(i=1,2)$ are computed in \cite[Section 2.3]{WanWang19ASS}, and up to $t-s\leqslant8$ the $E_{\infty}$-page of the ASS is given as in Figure \ref{Figure: ASS, spin case}, from which we read $\left(\Omega^{Spin}_{8}\left(K_{2}\right)\right)\sphat_{2}\cong\left(\mathbb{Z}_{(2)}\right)^{10}$. Hence $\Omega^{Spin}_{8}\left(K_{2}\right)\cong\mathbb{Z}^{10}$.
	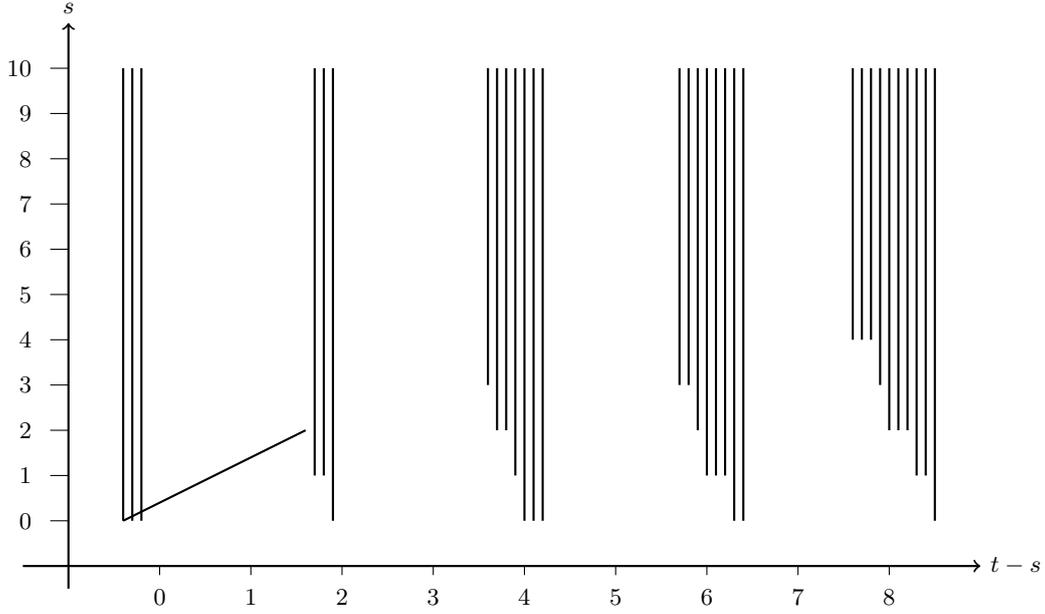
\begin{figure}[H]
		\begin{center}
			\begin{tikzpicture}[xscale=1.2, yscale=0.6]
				\draw[thick,->] (-1.5,-1) -- (9,-1) node[right] {\small $t - s$};
				\draw[thick,->] (-1,-1.5) -- (-1,11) node[above] {\small $s$};
				
				\draw[thick] (-0.4,0) -- (-0.4,10);
				\draw[thick] (-0.4,0) -- (1.6,2);
				\draw[thick] (3.6,3) -- (3.6,10);
				\draw[thick] (7.6,4) -- (7.6,10);
				
				\draw[thick] (-0.3,0) -- (-0.3,10);
				\draw[thick] (1.7,1) -- (1.7,10);
				\draw[thick] (3.7,2) -- (3.7,10);
				\draw[thick] (5.7,3) -- (5.7,10);
				\draw[thick] (7.7,4) -- (7.7,10);
				
				\draw[thick] (-0.2,0) -- (-0.2,10);
				\draw[thick] (1.8,1) -- (1.8,10);
				\draw[thick] (3.8,2) -- (3.8,10);
				\draw[thick] (5.8,3) -- (5.8,10);
				\draw[thick] (7.8,4) -- (7.8,10);
				
				\draw[thick] (1.9,0) -- (1.9,10);
				\draw[thick] (3.9,1) -- (3.9,10);
				\draw[thick] (5.9,2) -- (5.9,10);
				\draw[thick] (7.9,3) -- (7.9,10);
				
				\draw[thick] (4,0) -- (4,10);
				\draw[thick] (6,1) -- (6,10);
				\draw[thick] (8,2) -- (8,10);
				
				\draw[thick] (4.1,0) -- (4.1,10);
				\draw[thick] (6.1,1) -- (6.1,10);
				\draw[thick] (8.1,2) -- (8.1,10);
				
				\draw[thick] (4.2,0) -- (4.2,10);
				\draw[thick] (6.2,1) -- (6.2,10);
				\draw[thick] (8.2,2) -- (8.2,10);
				
				\draw[thick] (6.3,0) -- (6.3,10);
				\draw[thick] (8.3,1) -- (8.3,10);
				
				\draw[thick] (6.4,0) -- (6.4,10);
				\draw[thick] (8.4,1) -- (8.4,10);
				
				\draw[thick] (8.5,0) -- (8.5,10);

				\foreach \y in {0,1,2,3,4,5,6,7,8,9,10} {
					\draw (-1.2,\y) -- (-1,\y);  
					\node[left] at (-1.3, \y) {\small \y}; 
				}
				
				\foreach \x in {0,1,2,3,4,5,6,7,8} {
					\draw (\x,-1.2) -- (\x,-1);  
					\node[below] at (\x, -1.3) {\small \x}; 
				}					
			\end{tikzpicture}
		\end{center}
		\caption{$E_{\infty}$-page of Adams spectral sequence for $\Omega^{Spin}_{8}\left(K_{2}\right)$}
		\label{Figure: ASS, spin case}
	\end{figure}
	
	Now we construct an explicit isomorphism $\Omega^{Spin}_{8}\left(K_{2}\right)\to\mathbb{Z}^{10}$. The following homomorphism
	\begin{align*}
		\Omega^{Spin}_{8}\left(K_{2}\right)&\to\mathbb{Z}^{10},\\
		[M,x,y]&\mapsto \left[\sigma(M),\widehat{A}(M),\frac{1}{2}x^{4},x^{3}y,\frac{1}{2}x^{2}y^{2},xy^{3},\frac{1}{2}y^{4},\frac{\sigma(x,x)}{16},\frac{\sigma(x,y)-(x-y)^2}{8},\frac{\sigma(y,y)}{16}\right]
	\end{align*}
	is well-defined. Its first two components arise from bordism invariants of $\Omega_{8}^{Spin}$, the middle five come from identification of $E^{3}_{8,0}$ and the last three are induced by the signature of certain $4$-dimensional submanifolds.
	It is routine to check Table \ref{Table: epimorphism of spin bordism gp onto Z^10}. Hence the homomorphism we give above
	is an epimorphism. Since its domain and codomain are both isomorphic to $\mathbb{Z}^{10}$, it must be an isomorphism.
	\begin{table}[h]
		\centering
		\caption{Epimorhism $\varphi:\Omega^{Spin}_{8}\left(K_{2}\right)\to\mathbb{Z}^{10}$}
		\label{Table: epimorphism of spin bordism gp onto Z^10}
		\begin{tabular}{ccccccccccc}
			\toprule[1.5pt]
			$\Omega^{Spin}_{8}\left(K_{2}\right)$ & $\sigma$ & $\widehat{A}$ & $\frac{1}{2}x^{4}$ & $x^{3}y$ & $\frac{1}{2}x^{2}y^{2}$ & $xy^{3}$ & $\frac{1}{2}y^{4}$ & $\frac{\sigma(x,x)}{16}$ & $\frac{\sigma(x,y)-(x-y)^{2}}{8}$ & $\frac{\sigma(y,y)}{16}$\\
			\midrule[1pt]
			$\mathbb{H}P^{2}$ & 1 & 0 & 0 & 0 & 0 & 0 & 0 & 0 & 0 & 0\\
			$Bott$ & 224 & -1 & 0 & 0 & 0 & 0 & 0 & 0 & 0 & 0\\
			$[V(2),\widehat{z_{5}},0]$ & 2 & 0 & 1 & 0 & 0 & 0 & 0 & 0 & 0 & 0\\
			$\left[\mathbb{C}P^{3}\times\mathbb{C}P^{1},z_{3},z_{1}\right]$ & 0 & 0 & 0 & 1 & 0 & 0 & 0 & 0 & 0 & 0\\
			$[V(2),\widehat{z_{5}},\widehat{z_{5}}]$ & 2 & 0 & 1 & 2 & 1 & 2 & 1 & 0 & 0 & 0\\
			$\left[\mathbb{C}P^{3}\times\mathbb{C}P^{1},z_{1},z_{3}\right]$ & 0 & 0 & 0 & 0 & 0 & 1 & 0 & 0 & 0 & 0\\
			$[V(2),0,\widehat{z_{5}}]$ & 2 & 0 & 0 & 0 & 0 & 0 & 1 & 0 & 0 & 0\\
			$\left[\left(S^{2}\right)^{4},\Delta,0\right]$ & 0 & 0 & 12 & 0 & 0 & 0 & 0 & -1 & 0 & 0\\
			$\left[\left(S^{2}\right)^{4},x_{1}+x_{2},x_{3}+x_{4}\right]$ & 0 & 0 & 0 & 0 & 2 & 0 & 0 & 0 & 1 & 0\\
			$\left[\left(S^{2}\right)^{4},0,\Delta\right]$ & 0 & 0 & 0 & 0 & 0 & 0 & 12 & 0 & 0 & -1\\
			\bottomrule[1.5pt]
		\end{tabular}
	\end{table}
	
	Finally, it is routine to verify Table \ref{Table: char nums of a basis for spin bordism gp}.
\end{pf}

\begin{rmk}\label{Remark: bordism gp of K_1}
	In \cite{KS88,Kreck1991SomeNH} to develop the homeomorphism and diffeomorphism invariants for manifolds $M_{k,l}$ and $N_{k,l}$, the complete information of bordism groups $\Omega^{TopSpin}_{8}\left(K_{1}\right)$, $\Omega^{TopSpin}_{8}\left(K_{1};\gamma^{1}\right)$, $\Omega^{Spin}_{8}\left(K_{1}\right)$ and $\Omega^{Spin}_{8}\left(K_{1};\gamma^{1}\right)$ is not given explicitly. Instead the authors apply Hirzebruch integrality theorem, compare the topological and smooth bordism groups and give a detailed analysis on the filtrations associated to the AHSS. They directly compute the lattices generated by specific characteristic numbers evaluating on the subgroup of $8$-dimensional bordism group represented by manifolds with vanishing signature. 
	
	We can imitate the proof of Theorem \ref{Theorem: results of spin bordism gp} and deduce that $\Omega^{Spin}_{8}\left(K_{1}\right)\cong\Omega^{Spin}_{8}\left(K_{1};\gamma^{1}\right)\cong\mathbb{Z}^{4}$. To obtain $\Omega^{Spin}_{8}\left(K_{1}\right)$ it suffices to apply AHSS since there is no nontrivial extension. While in the calculation of $\Omega^{Spin}_{8}\left(K_{1};\gamma^{1}\right)$ there is a possibly nontrivial extension problem in AHSS and we need to apply ASS to determine the explicit isomorphism type. Meanwhile the epimorphisms exhibited in Tables \ref{Table: 8th spin bordism group of K_1} and \ref{Table: 8th nonspin bordism group of K_1} establish explicit isomorphisms from $\Omega^{Spin}_{8}\left(K_{1}\right)$ and $\Omega^{Spin}_{8}\left(K_{1};\gamma^{1}\right)$ to $\mathbb{Z}^{4}$. 
	\begin{table}[h]
		\centering
		\begin{minipage}[h]{0.37875\textwidth}
			\caption{Epimorphism $\Omega^{Spin}_{8}\left(K_{1}\right)\to\mathbb{Z}^{4}$}
			\label{Table: 8th spin bordism group of K_1}
			\begin{tabular}{ccccc}
				\toprule[1.5pt]
				$\Omega^{Spin}_{8}\left(K_{1}\right)$ & $\sigma$ & $\widehat{A}$ & $\frac{1}{2}x^{4}$ & $\frac{\sigma(x,x)}{16}$ \\
				\midrule[1pt]
				$\mathbb{H}P^{2}$ & 1 & 0 & 0 & 0\\
				$Bott$ & 224 & -1 & 0 & 0\\
				$\left[V(2),\widehat{z_{5}}\right]$ & 2 & 0 & 1 & 0\\
				$\left[\left(S^{2}\right)^{4},\Delta\right]$ & 0 & 0 & 12 & -1\\
				\bottomrule[1.5pt]
			\end{tabular}
		\end{minipage}
		\hfil
		\begin{minipage}[h]{0.45\textwidth}
			\centering
			\caption{Epimorphism $\Omega^{Spin}_{8}\left(K_{1};\gamma^{1}\right)\to\mathbb{Z}^{4}$}
			\label{Table: 8th nonspin bordism group of K_1}
			\begin{tabular}{ccccc}
				\toprule[1.5pt]
				$\Omega^{Spin}_{8}\left(K_{1};\gamma^{1}\right)$ & $\sigma$ & $e^{\frac{x}{2}}{A}$ & $x^{4}$ & $\frac{\sigma(x,x)-x^{2}}{8}$ \\
				\midrule[1pt]
				$\mathbb{H}P^{2}$ & 1 & 0 & 0 & 0\\
				$Bott$ & 224 & -1 & 0 & 0\\
				$\left[\mathbb{C}P^{4},z_{4}\right]$ & 1 & 0 & 1 & 0\\
				$\left[\left(\mathbb{C}P^{2}\right)^{2},\Delta'\right]$ & 1 & 0 & 6 & -1\\
				\bottomrule[1.5pt]
			\end{tabular}
		\end{minipage}
	\end{table}
\end{rmk}

\subsection{$s$-invariant for $\mathcal{E}_{1}^{+}$-manifolds}\label{Section:  s inv for E_1^+ mfds}
First we extend spin manifolds in $\mathcal{M}_{1}$ to the following broader class of manifolds.
\begin{defn}[$\mathcal{E}_{1}^{+}$-manifold]
	A spin $7$-manifold $M$ is called an \textbf{$\mathcal{E}_{1}^{+}$-manifold}, if it is simply connected, $H_{2}(M)\cong\mathbb{Z}^{2}$ and $H^{4}(M)$ is a finite group generated by $\overline{p_{1}}(M)$ and pairwise products of elements in $H^{2}(M)$. 
	
	A \textbf{polarization} of the $\mathcal{E}_{1}^{+}$-manifold $M$ is an ordered basis $(x,y)$ of $H^{2}(M)$, and the tripple $(M;x,y)$ is called a \textbf{polarized $\mathcal{E}_{1}^{+}$-manifold}. Polarized $\mathcal{E}_{1}^{+}$-manifolds $(M;x,y)$ and $\left(\overline{M};\overline{x},\overline{y}\right)$ are \textbf{polarized diffeomorphic} if there is a diffeomorphism $M\xrightarrow{f}\overline{M}$ such that $f^{*}\overline{x}=x$ and $f^{*}\overline{y}=y$, and we say $(M;x,y)\xrightarrow{f}\left(\overline{M};\overline{x},\overline{y}\right)$ is a \textbf{polarized diffeomorphism}.
\end{defn}
\begin{rmk}
	By definition any spin manifold in $\mathcal{M}_{1}$ is an $\mathcal{E}^{+}_{1}$-manifold.
	Polarized $\mathcal{E}_{1}^{+}$-manifolds are also called $F_{2}$-manifolds in \cite{HepworthPhdThesis}. It is clear that polarized diffeomorphism is an equivalence relation among polarized $\mathcal{E}^{+}_{1}$-manifolds. And if $(M;x,y)$ is a polarized $\mathcal{E}^{+}_{1}$-manifold and $\overline{M}\xrightarrow{f}M$ is a diffeomorphism, then we obtain a new polarized $\mathcal{E}^{+}_{1}$-manifold $\left(\overline{M};f^{*}x,f^{*}y\right)$.
\end{rmk}

Now we apply modified surgery theory to deduce the $s$-invariant of $\mathcal{E}_{1}^{+}$-manifolds. We refer to \cite{SurgeryAndDuality,Kreck2018OnTC} for basic definitions of modified surgery theory. First we shall determine their normal $2$-type and normal $2$-structures. They are given as in Figure \ref{normal 2 form & type of E_{1}^{+} mfd}.
\begin{figure}[h]
	\begin{equation*}
		\xymatrix{
			& BSpin\times K_{2}\ar[d]^{\mathrm{pr}_{1}}\\
			& BSpin\ar[d]^{B\phi}\\
			M\ar[r]^{\nu}\ar[ur]^{\widetilde{\nu}}\ar[uur]^{\left(\widetilde{\nu},x,y\right)} & BO
		}
	\end{equation*}
	\caption{Normal $2$-type and normal $2$-structure of $\mathcal{E}_{1}^{+}$-manifold $M$}
	\label{normal 2 form & type of E_{1}^{+} mfd}
\end{figure}
Here $Spin\xrightarrow{\phi}O$ is the canonical map, $B\phi$ is the induced map between classifying spaces, $\mathrm{pr}_{1}$ is projection onto the first component, $M\xrightarrow{\nu}BO$ is the classification map of stable normal bundle, $M\xrightarrow{\widetilde{\nu}}BSpin$ is a spin structure of stable normal bundle of $M$ and $M\xrightarrow{(x,y)}K_{2}$ is a $2$-equivalence. Since $M$ is simply connected and oriented with $H_{2}(M)\cong\mathbb{Z}^{2}$, $\widetilde{\nu}$ is unique and $(x,y)$ is equivalent to an ordered basis of $H^{2}(M)$.

Denote 
$(\mathcal{B},\mathcal{F})=\left(BSpin\times K_{2},B\phi\circ\mathrm{pr}_{1}\right)$
and let $\varphi=\left(\widetilde{\nu},x,y\right):M\to\mathcal{B}$ be a normal $2$-structure of $\mathcal{E}_{1}^{+}$-manifold $M$. We will not distinguish the notation $(M;x,y)$ and $(M,\varphi)$ while referring to a polarized $\mathcal{E}_{1}^{+}$-manifold. Then
$H^{4}\left(\mathcal{B};\mathbb{Q}\right)=\mathbb{Q}\left\{p_{1},\eta_{1}^{2},\eta_{1}\eta_{2},\eta_{2}^{2}\right\}.$
Since $H^{4}(M)$ is finite we have $H^{4}(M;\mathbb{Q})=0$ and
$$\ker H^{4}(\varphi;\mathbb{Q})=H^{4}\left(\mathcal{B};\mathbb{Q}\right)=\mathbb{Q}\left\{p_{1},\eta_{1}^{2},\eta_{1}\eta_{2},\eta_{2}^{2}\right\}.$$
Set $\kappa_{1}=p_{1}$, $\kappa_{2}=\eta_{1}^{2}$, $\kappa_{3}=\eta_{1}\eta_{2}$, $\kappa_{4}=\eta_{2}^{2}$. Suppose $(X,\psi)=(X;x,y)$ is a closed $8$-dimensional $(\mathcal{B},\mathcal{F})$-manifold, or a compact $8$-dimensional $(\mathcal{B},\mathcal{F})$-manifold whose boundary is a polarized $\mathcal{E}_{1}^{+}$-manifold. Set
$$S_{i,j}(X,\psi)=\left<\left(\psi^{*}\kappa_{i}\right)\cup\left(\psi^{*}\kappa_{j}\right),[X,\partial X]\right>\in\mathbb{Q},\ 1\leqslant i\leqslant j\leqslant4;$$
since $\psi^{*}\kappa_{1}=-p_{1}(X)$, $\psi^{*}\kappa_{2}=x^{2}$, $\psi^{*}\kappa_{3}=xy$, $\psi^{*}\kappa_{4}=y^{2}$ we also have the following expressions:
\begin{alignat*}{4}
	\left(p_{1}\right)^{2}\left(X;x,y\right)&:=&S_{1,1}\left(X;x,y\right)&=\left<\left(p_{1}(X)\right)^{2},[X,\partial X]\right>,\\
	\left(p_{1}x^{2}\right)\left(X;x,y\right)&:=&S_{1,2}\left(X;x,y\right)&=\left<p_{1}(X)x^{2},[X,\partial X]\right>,\\
	\left(p_{1}xy\right)\left(X;x,y\right)&:=&S_{1,3}\left(X;x,y\right)&=\left<p_{1}(X)xy,[X,\partial X]\right>,\\
	\left(p_{1}y^{2}\right)\left(X;x,y\right)&:=&S_{1,4}\left(X;x,y\right)&=\left<p_{1}(X)y^{2},[X,\partial X]\right>,\\
	\left(x^{2}\right)^{2}\left(X;x,y\right)&:=&S_{2,2}\left(X;x,y\right)&=\left<\left(x^{2}\right)^{2},[X,\partial X]\right>,\\
	\left(x^{2}\cdot xy\right)\left(X;x,y\right)&:=&S_{2,3}\left(X;x,y\right)&=\left<x^{2}\cdot xy,[X,\partial X]\right>,\\
	\left(x^{2}\cdot y^{2}\right)\left(X;x,y\right)&:=&S_{2,4}\left(X;x,y\right)&=\left<x^{2}\cdot y^{2},[X,\partial X]\right>,\\
	\left(xy\right)^{2}\left(X;x,y\right)&:=&S_{3,3}\left(X;x,y\right)&=\left<\left(xy\right)^{2},[X,\partial X]\right>,\\
	\left(xy\cdot y^{2}\right)\left(X;x,y\right)&:=&S_{3,4}\left(X;x,y\right)&=\left<xy\cdot y^{2},[X,\partial X]\right>,\\
	\left(y^{2}\right)^{2}\left(X;x,y\right)&:=&S_{4,4}\left(X;x,y\right)&=\left<\left(y^{2}\right)^{2},[X,\partial X]\right>.
\end{alignat*}
We evaluate $S=\left(S_{i,j}\right)_{1\leqslant i\leqslant j\leqslant4}$ on a basis $e$ of $\Omega_{8}^{Spin}\left(K_{2}\right)_{0}$, and according to Table \ref{Table: char nums of a basis for spin bordism gp} we obtain Table \ref{Table: matrix S(e) for E^{+}_{1} mfds}. Here $\Omega_{8}^{Spin}\left(K_{2}\right)_{0}\subset\Omega_{8}^{Spin}\left(K_{2}\right)$ is the subgroup of manifolds with vanishing signature.
\begin{table}[h]
	\centering
	\caption{Matrix $S(e)$ associated to $\mathcal{E}_{1}^{+}$-manifolds, $S:\Omega^{Spin}_{8}\left(K_{2}\right)_{0}\to\mathbb{Q}^{10}$}
	\label{Table: matrix S(e) for E^{+}_{1} mfds}
	\begin{tabular}{ccccccccccc}
		\toprule[1.5pt]
		\scriptsize$\Omega^{Spin}_{8}\left(K_{2}\right)_{0}$ & \scriptsize$p_{1}^{2}$ & \scriptsize$x^{2}p_{1}$ & \scriptsize$xyp_{1}$ & \scriptsize$y^{2}p_{1}$ & \scriptsize$\left(x^{2}\right)^{2}$ & \scriptsize$x^{2}\cdot xy$ & \scriptsize$x^{2}\cdot y^{2}$ & \scriptsize$(xy)^{2}$ & \scriptsize$xy\cdot y^{2}$ & \scriptsize$\left(y^{2}\right)^{2}$ \\
		\midrule[1pt]
		\scriptsize$Bott-224\mathbb{H}P^{2}$ & -896 & 0 & 0 & 0 & 0 & 0 & 0 & 0 & 0 & 0\\
		\scriptsize$\left[V(2),\widehat{z_{5}},0\right]-2\mathbb{H}P^{2}$ & 0 & 4 & 0 & 0 & 2 & 0 & 0 & 0 & 0 & 0\\
		\scriptsize$\left[\mathbb{C}P^{3}\times\mathbb{C}P^{1},z_{3},z_{1}\right]$ & 0 & 0 & 4 & 0 & 0 & 1 & 0 & 0 & 0 & 0\\
		\scriptsize$\left[V(2),\widehat{z_{5}},\widehat{z_{5}}\right]-2\mathbb{H}P^{2}$ & 0 & 4 & 4 & 4 & 2 & 2 & 2 & 2 & 2 & 2\\
		\scriptsize$\left[\mathbb{C}P^{3}\times\mathbb{C}P^{1},z_{1},z_{3}\right]$ & 0 & 0 & 4 & 0 & 0 & 0 & 0 & 0 & 1 & 0\\
		\scriptsize$[V(2),0,\widehat{z_{5}}]-2\mathbb{H}P^{2}$ & 0 & 0 & 0 & 4 & 0 & 0 & 0 & 0 & 0 & 2\\
		\scriptsize$\left[\left(S^{2}\right)^{4},\Delta,0\right]$ & 0 & 0 & 0 & 0 & 24 & 0 & 0 & 0 & 0 & 0\\
		\scriptsize$\left[\left(S^{2}\right)^{4},x_{1}+x_{2},x_{3}+x_{4}\right]$ & 0 & 0 & 0 & 0 & 0 & 0 & 4 & 4 & 0 & 0\\
		\scriptsize$\left[\left(S^{2}\right)^{4},0,\Delta\right]$ & 0 & 0 & 0 & 0 & 0 & 0 & 0 & 0 & 0 & 24\\
		\bottomrule[1.5pt]
	\end{tabular}
\end{table}
Then we derive an equivalent canonical form of $S(e)$ by invertible row transformations over $\mathbb{Z}$ and invertible column transformations over $\mathbb{Q}$, omit zero rows, record the expressions of column titles and obtain Table \ref{Table: E_{1}^{+} manifold canonical form under equiv, nondeg block}. 
\begin{table}[h]
	\centering
	\caption{Nonzero rows of generalized equivalent canonical form of matrix $S(e)$ associated to $\mathcal{E}_{1}^{+}$-manifolds}
	\label{Table: E_{1}^{+} manifold canonical form under equiv, nondeg block}
	\begin{tabular}{cccccccccc}
		\toprule[1.5pt]
		$S_{1}$ & $S_{2}$ &$S_{3}$ &$S_{4}$ &$S_{5}$ &$S_{6}$ &$S_{7}$ &$S_{8}$ &$S_{9}$ &$S_{10}$ \\
		\midrule[1pt]
		896 & 0 & 0 & 0 & 0 & 0 & 0 & 0 & 0 & 0 \\
		0 & 48 & 0 & 0 & 0 & 0 & 0 & 0 & 0 & 0\\
		0 & 0 & 24 & 0 & 0 & 0 & 0 & 0 & 0 & 0\\
		0 & 0 & 0 & 48 & 0 & 0 & 0 & 0 & 0 & 0\\
		0 & 0 & 0 & 0 & 2 & 0 & 0 & 0 & 0 & 0\\
		0 & 0 & 0 & 0 & 0 & 1 & 0 & 0 & 0 & 0\\
		0 & 0 & 0 & 0 & 0 & 0 & 2 & 0 & 0 & 0\\
		0 & 0 & 0 & 0 & 0 & 0 & 0 & 1 & 0 & 0\\
		0 & 0 & 0 & 0 & 0 & 0 & 0 & 0 & 2 & 0\\
		\bottomrule[1.5pt]
	\end{tabular}
\end{table}
Its row vectors generate a integral lattice $L$ and its column titles are given by
\begin{align*}
	S_{1}&=p_{1}^{2},\\
	S_{2}&=x^{2}p_{1}-2\left(x^{2}\right)^{2},\\
	S_{3}&=xyp_{1}-4x^{2}\cdot xy+6x^{2}\cdot y^{2}-4xy\cdot y^{2},\\
	S_{4}&=y^{2}p_{1}-2\left(y^{2}\right)^{2},\\
	S_{5}&=\left(x^{2}\right)^{2},\\
	S_{6}&=x^{2}\cdot xy,\\
	S_{7}&=x^{2}\cdot y^{2},\\
	S_{8}&=xy\cdot y^{2},\\
	S_{9}&=\left(y^{2}\right)^{2},\\
	S_{10}&=x^{2}\cdot y^{2}-(xy)^{2}.
\end{align*}

Hence if $\left(V;\widehat{x},\widehat{y}\right)$ is a coboundary of the polarized $\mathcal{E}_{1}^{+}$-manifold $(M;x,y)$ with vanishing signature, we set
\begin{align*}
	s_{i}(M;x,y)&:=S_{i}\left(V;\widehat{x},\widehat{y}\right)\mathrm{mod}\ d_{i}\mathbb{Z}\in\mathbb{Q}\big/d_{i}\mathbb{Z}\\
	&=\frac{1}{d_{i}}S_{i}\left(V;\widehat{x},\widehat{y}\right)\mathrm{mod}\ \mathbb{Z}\in\mathbb{Q}\big/\mathbb{Z},\ 1\leqslant i\leqslant9;\\
	S_{10}(M;x,y)&:=S_{10}\left(V;\widehat{x},\widehat{y}\right)\in\mathbb{Q}
\end{align*}
where $d_{i}$ is the unique nonzero entry in the $i$-th row of Table \ref{Table: E_{1}^{+} manifold canonical form under equiv, nondeg block}.
\begin{thm}[Polarized diffeomorphism invariant for polarized $\mathcal{E}_{1}^{+}$-manifolds]\label{Theorem: Polarized diffeom inv for polarized E_{1}^{+}-mfds}
	Invariants $s_{i}(M;x,y)$ $(1\leqslant i\leqslant9)$ and $S_{10}(M;x,y)$ are well-defined for a polarized $\mathcal{E}_{1}^{+}$-manifold $(M;x,y)$. Two $\mathcal{E}_{1}^{+}$-manifolds $M$ and $\overline{M}$ are diffeomorphic if and only if there are normal $2$-structures $M\xrightarrow{\varphi}\mathcal{B}$, $\overline{M}\xrightarrow{\overline{\varphi}}\mathcal{B}$ and a tangential isomorphism compactible to the normal $2$-structures, such that the polarized $\mathcal{E}_{1}^{+}$-manifolds $\left(M,\varphi\right)$ and $\left(\overline{M},\overline{\varphi}\right)$ have the same invariants
	$$s=\left(s_{1},\cdots,s_{9},S_{10}\right)\in\left(\mathbb{Q}\big/\mathbb{Z}\right)^{9}\times\mathbb{Q}.$$
\end{thm}

\begin{pf}
	First we show that these invariants are well-defined. We shall prove that different choices of $(\mathcal{B},\mathcal{F})$-coboundaries for the same polarized $\mathcal{E}_{1}^{+}$-manifold $(M;x,y)$ lead the same value of $s(M;x,y)$. Let $\left(W_{i};x_{i},y_{i}\right)$ $(i=1,2)$ be two $(\mathcal{B},\mathcal{F})$-coboundaries of polarized $\mathcal{E}_{1}^{+}$-manifold $(M;x,y)$. Form the closed oriented $8$-manifold $X:=W_{1}\cup_{M}\left(-W_{2}\right)$, and by Mayer-Vietoris sequence there are unique elements $\widetilde{x}$, $\widetilde{y}\in H^{2}(X)$ such that $\left.\widetilde{x}\right|_{W_{i}}=x_{i}$ and $\left.\widetilde{y}\right|_{W_{i}}=y_{i}$. Then $\left(X;\widetilde{x},\widetilde{y}\right)$ is a closed $8$-dimensional $(\mathcal{B},\mathcal{F})$-manifold and 
	$$S_{i}\left(X;\widetilde{x},\widetilde{y}\right)=S_{i}\left(W_{1};x_{1},y_{1}\right)-S_{i}\left(W_{2};x_{2},y_{2}\right),1\leqslant i\leqslant10.$$
	Therefore, different $(\mathcal{B},\mathcal{F})$-coboundaries of the same polarized $\mathcal{E}_{1}^{+}$-manifold lead the same value of $s$-invariant, and the $s$-invariant is well-defined.
	
	Next we show the necessity.
	Suppose $\mathcal{E}_{1}^{+}$-manifolds $M$ and $\overline{M}$ are diffeomorphic and $M\xrightarrow{f}\overline{M}$ is a diffeomorphism. Then $H^{*}\left(\overline{M}\right)\xrightarrow{f^{*}}H^{*}(M)$ is a tangential isomorphism, namely it induces isomorphisms on the second and fourth cohomology groups, preserves the cup product structure $H^{2}\otimes H^{2}\xrightarrow{\cdot\cup\cdot}H^{4}$ and the characteristic classes $w_{2}$, $p_{1}$. Let $(x,y)$ be a polarization of $M$. Then $\left(\overline{x},\overline{y}\right):=f^{*}(x,y)$ is a polarization of $\overline{M}$ and $f^{*}$ is clearly compactible to polarizations $(x,y)$, $\left(\overline{x},\overline{y}\right)$. Let $\left(V;x',y'\right)$ be a $(\mathcal{B},\mathcal{F})$-coboundary of $(M;x,y)$ with vanishing signature. By h-cobordism theorem there is a compact $8$-manifold $W$ with boundary $\partial W=(-M)\sqcup\overline{M}$ and a diffeomorphism $W\xrightarrow{F}M\times I$ such that $\left.F\right|_{M}=\mathrm{id}_{M}$ and $\left.F\right|_{\overline{M}}=f$. Then $\left(W;F^{*}\mathrm{pr}_{1}^{*}x,F^{*}\mathrm{pr}_{1}^{*}y\right)$ is an h-cobordism of $(\mathcal{B},\mathcal{F})$-manifolds $(M;x,y)$ and $\left(\overline{M};\overline{x},\overline{y}\right)$. Form $\overline{V}=V\cup_{M} W$, and by Mayer-Vietoris sequence there are unique classes $\overline{x'}$, $\overline{y'}\in H^{2}\left(\overline{V}\right)$ such that $\left.\overline{x'}\right|_{V}=\overline{x}$, $\left.\overline{x'}\right|_{W}=F^{*}\mathrm{pr}_{1}^{*}x$ and $\left.\overline{y'}\right|_{V}=\overline{y}$, $\left.\overline{y'}\right|_{W}=F^{*}\mathrm{pr}_{1}^{*}y$. Then $\left(\overline{V};\overline{x'},\overline{y'}\right)$ is a coboundary of $\left(\overline{M};\overline{x},\overline{y}\right)$ with vanishing signature, and it is straightforward to compute that $S_{i}\left(V;x',y'\right)=S_{i}\left(\overline{V};\overline{x'},\overline{y'}\right)$ for all $1\leqslant i\leqslant10$. Hence $s(M;x,y)=s\left(\overline{M};\overline{x},\overline{y}\right)$.
	
	Finally we show the suffiency. Suppose there are normal $2$-structures $M\xrightarrow{\varphi}\mathcal{B}$, $\overline{M}\xrightarrow{\overline{\varphi}}\mathcal{B}$ and a tangential isomorphism compactible to the normal $2$-structures such that $s\left(M,\varphi\right)=s\left(\overline{M},\overline{\varphi}\right)$. Let $\left(V,\widetilde{\varphi}\right)$ and $\left(\overline{V},\widetilde{\overline{\varphi}}\right)$ be the $(\mathcal{B},\mathcal{F})$-coboundaries of $\left(M,\varphi\right)$ and $\left(\overline{M},\overline{\varphi}\right)$ with vanishing signature respectively. Then it follows from $s\left(M,\varphi\right)=s\left(\overline{M},\overline{\varphi}\right)$ that $S\left(V,\widetilde{\varphi}\right)-S\left(\overline{V},\widetilde{\overline{\varphi}}\right)\in S\left(\Omega_{8}^{Spin}\left(K_{2}\right)_{0}\right)$, namely there is a closed $8$-dimensional $(\mathcal{B},\mathcal{F})$-manifold $(X;\psi)$ with vanishing signature, such that $S\left(V,\widetilde{\varphi}\right)-S\left(\overline{V},\widetilde{\overline{\varphi}}\right)=S(X,\psi)$. Then $(W,\Psi):=\left(V,\widetilde{\varphi}\right)\#\left(-\left(\overline{V},\widetilde{\overline{\varphi}}\right)\right)\#\left(-(X,\psi)\right)$ is a $(\mathcal{B},\mathcal{F})$-bordism between $\left(M,\varphi\right)$ and $\left(\overline{M},\overline{\varphi}\right)$ such that its signature vanishes and $\left<\left(\Psi^{*}\kappa\right)\cup\left(\Psi^{*}\kappa'\right),[W,\partial W]\right>=0$ for all $\kappa,\kappa'\in\ker H^{4}\left(\varphi\sqcup\overline{\varphi};\mathbb{Q}\right)$. Hence it follows from \cite[Theorem 6]{SurgeryAndDuality} that $M$ and $\overline{M}$ are diffeomorphic.
\end{pf}

\begin{rmk}
	If we intend to deduce $s$-invariants for non-spin manifolds in $\mathcal{M}_{1}$ or manifolds in $\mathcal{M} _{i}$, $2\leqslant i\leqslant5$ and classify them, we shall apply \cite[Theorem 5(a)]{Kreck2018OnTC} to deduce sufficient and necessary condition for two manifolds being diffeomorphic, and we need information of $\Omega^{Spin}_{8}\left(K_{2};\mathrm{pr}_{1}^{*}\gamma^{1}\right)$ to derive $s$-invariants for non-spin manifolds.
\end{rmk}

We end this section by remarks and properties on components of $s$-invariants.

\begin{rmk}\label{Remark: s1 of E_{1}^{+} mfd is ori-pres diff inv}
	Invariant $s_{1}$ does not depend on the choice of polarizations and is a diffeomorphism invariant of $\mathcal{E}_{1}^{+}$-manifolds. It can be viewed as an extension of Eells-Kuiper $\mu$-invariant (\cite{EKinv}).
\end{rmk}

\begin{rmk}
	By \cite[Proposition 3.3.2]{HepworthPhdThesis} invariant $S_{10}$ is related to Massey triple product $\langle\cdot,\cdot,\cdot\rangle$ of $M$ by
	$$S_{10}(M;x,y)=-\left<\left<x,x,y\right>y,[M]\right>.$$
	When $(X,\psi)$ is a closed $8$-dimensional $(\mathcal{B},\mathcal{F})$-manifold we have $S_{10}(X,\psi)=0$, since products of degree $2$ cohomology classes of a closed manifold are commutative. However, if $(X,\psi)$ is a compact $8$-dimensional $(\mathcal{B},\mathcal{F})$-manifold whose boundary $(M,\varphi)$ is an $\mathcal{E}_{1}^{+}$-manifold, $S_{10}(M,\varphi)=S_{10}(X,\psi)=0$ fails in general.
	
	The value of $S_{10}$ does not depend on choices of polarizations and is a diffeomorphism invariant of $\mathcal{E}_{1}^{+}$-manifolds. See Proposition \ref{Proposition: S_{10} does not depend on polarization} for more details.
\end{rmk}

\begin{rmk}
	Next we consider invariants $s_{2}\sim s_{9}$.
	
	\begin{compactenum} 
		\item By \cite[Proposition 3.3.2]{HepworthPhdThesis} invariants $s_{5}\sim s_{9}$ are related to the linking form $lk_{M}$ of $M$ by
		\begin{align*}
			2s_{5}(M;x,y)&=-lk_{M}\left(x^{2},x^{2}\right),\\
			s_{6}(M;x,y)&=-lk_{M}\left(x^{2},xy\right),\\
			2s_{7}(M;x,y)&=-lk_{M}\left(x^{2},y^{2}\right),\\
			s_{8}(M;x,y)&=-lk_{M}\left(xy,y^{2}\right),\\
			2s_{9}(M;x,y)&=-lk_{M}\left(y^{2},y^{2}\right);
		\end{align*} 
		\item Invariants $s_{2}\sim s_{9}$ can be assembled together as a new invariant that does not depend on choices of polarization and is a diffeomorphism invariant of $\mathcal{E}_{1}^{+}$-manifolds. Then we can express the $s$-invariant as a \textbf{diffeomorphism invariant} of $\mathcal{E}_{1}^{+}$-manifolds.
		See Proposition \ref{Proposition: inv s2 to s9 polarization free} for more details.
	\end{compactenum}
\end{rmk}	

\begin{prop}\label{Proposition: S_{10} does not depend on polarization}
	The invariant $S_{10}$ of a polarized $\mathcal{E}_{1}^{+}$-manifold does not depend on polarizations and is a diffeomorphism invariant of $\mathcal{E}_{1}^{+}$-manifolds.
\end{prop}

\begin{pf}
	Let $M$ be an $\mathcal{E}_{1}^{+}$-manifold and let $(x,y)$, $\left(x',y'\right)$ be polarizations of $M$. Let $\left(V;\widehat{x},\widehat{y}\right)$ be a $(\mathcal{B},\mathcal{F})$-coboundary of $(M;x,y)$ with vanishing signature. 
	Since $\left(x',y'\right)$ and $(x,y)$ are ordered bases of $H^{2}(M)$, there exists $A$, $B$, $C$, $D\in\mathbb{Z}$ such that $x'=Ax+By$, $y'=Cx+Dy$ and $AD-BC=\pm1$. Then $\left(V;A\widehat{x}+B\widehat{y},C\widehat{x}+D\widehat{y}\right)$ is a $(\mathcal{B},\mathcal{F})$-coboundary of $(M;Ax+By,Cx+Dy)$ with vanishing signature. 
	
	We shall show that $S_{10}(M;Ax+By,Cx+Dy)=S_{10}(M;x,y)$, namely $S_{10}\left(V;A\widehat{x}+B\widehat{y},C\widehat{x}+D\widehat{y}\right)=S_{10}\left(V;\widehat{x},\widehat{y}\right)$.
	Let $V\xrightarrow{j}(V,M)$ be the inclusion. Recall that $H^{4}(M;\mathbb{Q})=H^{3}(M;\mathbb{Q})=0$. From the long exact sequence of rational cohomology groups associated to the pair $(V,M)$ we have a linear isomorphism $H^{4}(V,M;\mathbb{Q})\xrightarrow{j^{*}}H^{4}(V;\mathbb{Q})$. Let $\Lambda$ be its inverse. Then
	\begin{align*}
		&S_{10}\left(V;\widehat{x},\widehat{y}\right)
		=\left\langle\left(\Lambda\left(\widehat{x}^{2}\right)\right)\left(\Lambda\left(\widehat{y}^{2}\right)\right)-\left(\Lambda\left(\widehat{x}\widehat{y}\right)\right)^{2},[V,M]\right\rangle\\
		&S_{10}\left(V;A\widehat{x}+B\widehat{y},C\widehat{x}+D\widehat{y}\right)\\
		=&\left\langle\left(\Lambda\left(\left(A\widehat{x}+B\widehat{y}\right)^{2}\right)\right)\left(\Lambda\left(\left(C\widehat{x}+D\widehat{y}\right)^{2}\right)\right)-\left(\Lambda\left(\left(A\widehat{x}+B\widehat{y}\right)\left(C\widehat{x}+D\widehat{y}\right)\right)\right)^{2},[V,M]\right\rangle\\
		=&(AD-BC)^{2}\left\langle\left(\Lambda\left(\widehat{x}^{2}\right)\right)\left(\Lambda\left(\widehat{y}^{2}\right)\right)-\left(\Lambda\left(\widehat{x}\widehat{y}\right)\right)^{2},[V,M]\right\rangle
		=S_{10}\left(V;\widehat{x},\widehat{y}\right).
	\end{align*}
	Therefore, given an $\mathcal{E}_{1}^{+}$-manifold $M$, the value $S_{10}(M;x,y)$ does not depend on che choice of polarizations $(x,y)$ and is a diffeomorphism invariant of the $\mathcal{E}_{1}^{+}$-manifold $M$ itself.
\end{pf}

\begin{prop}\label{Proposition: inv s2 to s9 polarization free}
	We can assemble invariants $s_{2}\sim s_{9}$ as one invariant 
	$$\varsigma\in\left(\mathbb{Q}^{8}\big/ L_{0}\right)\Big/GL(2,\mathbb{Z})$$
	that does not depend on polarizations and is a diffeomorphism invariant of $\mathcal{E}^{+}_{1}$-manifolds, and the $s$-invariant of polarized $\mathcal{E}^{+}_{1}$-manifolds can be reinterpreted as a diffeomorphism invariant of $\mathcal{E}^{+}_{1}$-manifolds:
	$$\left(s_{1},\varsigma,S_{10}\right)\in\left(\mathbb{Q}\big/896\mathbb{Z}\right)\times\left(\left(\mathbb{Q}^{8}\big/ L_{0}\right)\Big/GL(2,\mathbb{Z})\right)\times\mathbb{Q}.$$
\end{prop}

\begin{pf}
	Recall that given a polarized $\mathcal{E}^{+}_{1}$-manifold $(M;x,y)$ its $s$-invariant is 
	$$s(M;x,y)=S\left(W;\widehat{x},\widehat{y}\right)+L\in\mathbb{Q}^{10}\big/L,$$
	where $\left(W;\widehat{x},\widehat{y}\right)$ is a coboundary of $(M;x,y)$ with signature zero and $L$ is the integral lattice generated by the row vectors of matrix $S(e)$ (see Table \ref{Table: E_{1}^{+} manifold canonical form under equiv, nondeg block}). Recall also that $s_{1}$ and $S_{10}$ do not depend on polarizations. Set
	\begin{align*}
		\mathbb{Q}^{10}&=\mathbb{Q}\times\mathbb{Q}^{8}\times\mathbb{Q},\\ \check{S}\left(V;z,w\right)&:=\left(z^{2}p_{1},zwp_{1},w^{2}p_{1},\left(z^{2}\right)^{2},z^{2}\cdot zw,z^{2}\cdot w^{2},zw\cdot w^{2},\left(w^{2}\right)^{2}\right)\left(V;x,y\right),\\
		L_{0}&:=\check{S}\left(\Omega_{8}^{Spin}\left(K_{2}\right)_{0}\right),\\
		\varsigma(M;x,y)&:=\check{S}\left(W;\widehat{x},\widehat{y}\right)+L_{0}\in\mathbb{Q}^{8}\big/L_{0}.
	\end{align*} 
	Here $\left(V;z,w\right)$ is any compact $8$-dimensional $(\mathcal{B},\mathcal{F})$-manifold and $\check{S}$ is viewed as a row vector. Then $L_{0}$ is an integral lattice in $\mathbb{Q}^{8}$ spanned by the row vectors of
	$$\widetilde{L_{0}}=\begin{pmatrix}
		4 & 0 & 0 & 2 & 0 & 0 & 0 & 0\\
		0 & 4 & 0 & 0 & 1 & 0 & 0 & 0\\
		0 & -12 & 0 & 0 & 0 & 2 & 0 & 0\\
		0 & 4 & 0 & 0 & 0 & 0 & 1 & 0\\
		0 & 0 & 4 & 0 & 0 & 0 & 0 & 2\\
		0 & 0 & 0 & 24 & 0 & 0 & 0 & 0\\
		0 & 0 & 0 & 0 & 0 & 4 & 0 & 0\\
		0 & 0 & 0 & 0 & 0 & 0 & 0 & 24
	\end{pmatrix}$$
	(this can be read from Table \ref{Table: matrix S(e) for E^{+}_{1} mfds} and we applied some invertible row transformations over $\mathbb{Z}$) and $s$-invariant can be rephrased as
	$$s(M;x,y)=\left(s_{1}(M),\varsigma(M;x,y\right),S_{10}(M))\in \mathbb{Q}\big/896\mathbb{Z}\times \mathbb{Q}^{8}\big/ L_{0}\times \mathbb{Q}$$
	
	We focus on $\check{S}$ and $\varsigma$. Given a compact $8$-dimensional $(\mathcal{B},\mathcal{F})$-manifold $\left(V;z,w\right)$ and $g=\begin{pmatrix}
		A & C \\ B & D
	\end{pmatrix}\in GL(2,\mathbb{Z})$, we have $\check{S}\left(V;(z,w)g\right)=\check{S}\left(V;z,w\right)R(g)$, where
	$$R(g)=
	\scriptsize\begin{pmatrix}
		A^2 & A C & C^2 & 0 & 0 & 0 & 0 & 0 \\
		2 A B & A D+B C & 2 C D & 0 & 0 & 0 & 0 & 0 \\
		B^2 & B D & D^2 & 0 & 0 & 0 & 0 & 0 \\
		0 & 0 & 0 & A^4 & A^3 C & A^2 C^2 & A C^3 & C^4 \\
		0 & 0 & 0 & 4 A^3 B & A^3 D+3 A^2 B C & 2 A^2 C D+2 A B C^2 & 3 A C^2 D+B C^3 & 4 C^3 D \\
		0 & 0 & 0 & 6 A^2 B^2 & 3 A^2 B D+3 A B^2 C & A^2 D^2+4 A B C D+B^2 C^2 & 3 A C D^2+3 B C^2 D & 6 C^2 D^2 \\
		0 & 0 & 0 & 4 A B^3 & 3 A B^2 D+B^3 C & 2 A B D^2+2 B^2 C D & A D^3+3 B C D^2 & 4 C D^3 \\
		0 & 0 & 0 & B^4 & B^3 D & B^2 D^2 & B D^3 & D^4
	\end{pmatrix}.
	$$
	This gives a right action of $GL(2,\mathbb{Z})$ on $\mathbb{Q}^{8}$. Moreover, lattice $L_{0}$ is invariant under the $GL(2,\mathbb{Z})$-action. This follows from direct computation of the matrix $\widetilde{L_{0}}R(g)\widetilde{L_{0}}^{-1}$: after modulo transparent integer terms we have
	\begin{align*}
		&\widetilde{L_{0}}R(g)\widetilde{L_{0}}^{-1}\\
		\equiv&\scriptsize\begin{pmatrix}
			0 & 0 & 4\left(f_{3}(A)C+Af_{3}(C)\right)-AC & 0 & 0 & f_{4}(A) & f_{2}^{-}(AC)-2\left(f_{3}(A)C+Af_{3}(C)\right) & f_{4}(C) \\
			0 & 0 & 2f_{3}(A)D+2Bf_{3}(C) & 0 & 0 & f_{3}(A)B & g_{2}(A,C)(D-B) -f_{3}(A)D-Bf_{3}(C)  & f_{3}(C)D \\
			0 & 0 & 0 & 0 & 0 & f_{2}^{+}(AB) & f_{2}^{-}(AD)+f_{2}^{-}(BC) & f_{2}^{+}(CD)\\
			0 & 0 & 2f_{3}(B)C+2Af_{3}(D) & 0 & 0 & Af_{3}(B) & g_{2}(B,D)(C-A) -f_{3}(B)C-Af_{3}(D) & Cf_{3}(D) \\
			0 & 0 & 4\left(f_{3}(B)D+Bf_{3}(D)\right)-BD  & 0 & 0 & f_{4}(B) & f_{2}^{-}(BD)-2\left(f_{3}(B)D+Bf_{3}(D)\right) & f_{4}(D) \\
			0 & 0 & 0 & 0 & 0 & 0 & 0 & 0 \\
			0 & 0 & 0 & 0 & 0 & 0 & 0 & 0 \\
			0 & 0 & 0 & 0 & 0 & 0 & 0 & 0
		\end{pmatrix}\ \mathrm{mod}\ 1,
	\end{align*}
	where $f_{2}^{\pm}(n):=\frac{n^{2}\pm n}{2}$, $g_{2}(m,n):=\frac{m^{2}n-mn^{2}}{2}$, $f_{3}(n):=\frac{n^{3}-n}{6}$ and $f_{4}(n):=\frac{n^{4}-n^{2}}{12}$
	are polynomials that take integral values whenever $m$ and $n$ are integers.
	
	Therefore, the right $GL(2,\mathbb{Z})$-action on $\mathbb{Q}^{8}$ induces an action on $\mathbb{Q}^{8}\big/ L_{0}$, and the invariant
	$$\varsigma(M):=\varsigma(M;x,y)GL(2,\mathbb{Z})\in\left(\mathbb{Q}^{8}\big/ L_{0}\right)\Big/GL(2,\mathbb{Z})$$
	is well-defined. We conclude that $\varsigma(M)\in\left(\mathbb{Q}^{8}\big/ L_{0}\right)\Big/GL(2,\mathbb{Z})$ does not depend on the choices of polarizations and is a diffeomorphism invariant of $M$.
\end{pf}
	
	\section{Partial classification of manifolds in $\mathcal{M}_{1}$}\label{Section: 7 Partial classification of manifolds in M_1}

In this section we give a partial classification of spin manifolds in $\mathcal{M}_{1}$. We derived the polarized diffeomorphism invariants for $\mathcal{E}^{+}_{1}$-manifolds in Section \ref{Section:  s inv for E_1^+ mfds}, and now we compute the $s$-invariants of certain spin manifolds in $\mathcal{M}_{1}$ in Section \ref{Section: compute s inv of some spin mfds in M_1}. Then in Section \ref{Section: partial classification of spin mfds in M_1} we compare the $s$-invariants, give a partial classification of spin manifolds in $\mathcal{M}_{1}$ and prove Theorem \ref{Theorem: Partial classification of mfds in M_1} for spin manifolds.

\subsection{Computations of $s$-invariants for some spin manifolds in $\mathcal{M}_{1}$}\label{Section: compute s inv of some spin mfds in M_1}

In this section we compute $s$-invariants for certain spin manifolds in $\mathcal{M}_{1}$.

\begin{lem}[$s$-invariants of spin $M_{m,n,l}$ with $\gcd(n,l)=1$]\label{Lemma: s-inv of spin M_{m,n,l} with gcd(n,l)=1}
	Assume integer $m$ and non-zero integers $n$, $l$ satisfy $\left(m,n,l\right)\equiv(0,1,0)\ \mathrm{mod}\ 2$ and $\gcd(n,l)=1$. Set $z=-q\beta+p\gamma$, where integers $p$ and $q$ satisfy $pn+ql=1$. For $\begin{pmatrix}
		A & C \\ B & D
	\end{pmatrix}\in GL(2,\mathbb{Z})$ the $s$-invariant of $\left(M_{m,n,l};Ap^{*}z+Bp^{*}\alpha,Cp^{*}z+Dp^{*}\alpha\right)$ is given by
		\begin{eqnarray}\label{Equation: s-inv of spin M_{m,n,l}, (n,l)=1}
		\small
		{\left\{
		\begin{alignedat}{2}
			s_{1}&=\frac{3m\left(n^{2}-1\right)\left(n^{2}+3\right)}{896n^{2}}+\frac{\left(l^{2}+4\right)^{2}}{896l}-\frac{\mathrm{sgn}(l)}{224}\in\mathbb{Q}\big/\mathbb{Z},\\
			s_{2}&=\frac{A^{2}}{48}\left(\frac{q^{2}m\left(n^{2}-3\right)}{n}+\frac{p^{2}\left(l^{2}+4\right)}{l}\right)-\frac{AB}{24}\cdot\frac{q\left(n^{2}+3\right)}{n}-\frac{A^{4}}{24}\left(\frac{p^{4}}{l}-\frac{q^{4}m}{n^{2}}\right)+\frac{A^{3}B}{6}\cdot\frac{q^{3}}{n}\in\mathbb{Q}\big/\mathbb{Z},\\
			s_{3}&=\frac{AC}{24}\left(\frac{q^{2}m\left(n^{2}-3\right)}{n}+\frac{p^{2}\left(l^{2}+4\right)}{l}\right)-\frac{AD+BC}{24}\cdot\frac{q\left(n^{2}+3\right)}{n}\\
			&\ -\frac{AC\left(2A^{2}-3AC+2C^{2}\right)}{12}\left(\frac{p^{4}}{l}-\frac{q^{4}m}{n^{2}}\right)+\left(\frac{A^{3}D+BC^{3}}{6}+\frac{AC(A-C)(B-D)}{2}\right)\cdot\frac{q^{3}}{n}\in\mathbb{Q}\big/\mathbb{Z},\\
			s_{4}&=\frac{C^{2}}{48}\left(\frac{q^{2}m\left(n^{2}-3\right)}{n}+\frac{p^{2}\left(l^{2}+4\right)}{l}\right)-\frac{CD}{24}\cdot\frac{q\left(n^{2}+3\right)}{n}-\frac{C^{4}}{24}\left(\frac{p^{4}}{l}-\frac{q^{4}m}{n^{2}}\right)+\frac{C^{3}D}{6}\cdot\frac{q^{3}}{n}\in\mathbb{Q}\big/\mathbb{Z},\\
			s_{5}&=\frac{A^{4}}{2}\left(\frac{p^{4}}{l}-\frac{q^{4}m}{n^{2}}\right)-2A^{3}B\cdot\frac{q^{3}}{n}\in\mathbb{Q}\big/\mathbb{Z},\\
			s_{6}&=A^{3}C\left(\frac{p^{4}}{l}-\frac{q^{4}m}{n^{2}}\right)-A^{2}\left(AD+3BC\right)\cdot\frac{q^{3}}{n}\in\mathbb{Q}\big/\mathbb{Z},\\
			s_{7}&=\frac{A^{2}C^{2}}{2}\left(\frac{p^{4}}{l}-\frac{q^{4}m}{n^{2}}\right)-AC(AD+BC)\cdot\frac{q^{3}}{n}\in\mathbb{Q}\big/\mathbb{Z},\\
			s_{8}&=AC^{3}\left(\frac{p^{4}}{l}-\frac{q^{4}m}{n^{2}}\right)-C^{2}\left(3AD+BC\right)\cdot\frac{q^{3}}{n}\in\mathbb{Q}\big/\mathbb{Z},\\
			s_{9}&=\frac{C^{4}}{2}\left(\frac{p^{4}}{l}-\frac{q^{4}m}{n^{2}}\right)-2C^{3}D\cdot\frac{q^{3}}{n}\in\mathbb{Q}\big/\mathbb{Z},\\
			S_{10}&=0\in\mathbb{Q}.\\
		\end{alignedat}
		\right.}
	\end{eqnarray}
\end{lem}

\begin{lem}[$s$-invariants of spin $M_{m,n,l}$ with $\gcd(m,n)=1$]\label{Lemma: s-inv of spin M_{m,n,l} with gcd(m,n)=1}
	Assume non-zero integers $m$, $n$ and $l$ satisfy $\left(m,n,l\right)\equiv(0,1,0)\ \mathrm{mod}\ 2$ and $\gcd(m,n)=1$. Set $z=-v\alpha+u\beta$, where integers $u$ and $v$ satisfy $um+vn=1$. For $\begin{pmatrix}
		A & C \\ B & D
	\end{pmatrix}\in GL(2,\mathbb{Z})$ the $s$-invariant of $\left(M_{m,n,l};Ap^{*}z+Bp^{*}\gamma,Cp^{*}z+Dp^{*}\gamma\right)$ is given by
	\begin{eqnarray}\label{Equation: s-inv of spin M_{m,n,l}, (m,n)=1}
		\small
		{\left\{
		\begin{alignedat}{2}
			s_{1}&=\frac{3m\left(n^{2}-1\right)\left(n^{2}+3\right)}{896n^{2}}+\frac{\left(l^{2}+4\right)^{2}}{896l}-\frac{\mathrm{sgn}(l)}{224}\in\mathbb{Q}\big/\mathbb{Z},\\
			s_{2}&=\frac{A^{2}}{48}\left(u-3uvn-\frac{3u}{n^{2}}-\frac{3uv}{n}\right)+\frac{B^{2}}{48}\cdot\frac{l^{2}+4}{l}+\frac{A^{4}}{24}\left(\frac{u^{3}}{n^{2}}+\frac{3u^{3}v}{n}\right)-\frac{B^{4}}{24l}\in\mathbb{Q}\big/\mathbb{Z},\\
			s_{3}&=\frac{AC}{24}\left(u-3uvn-\frac{3u}{n^{2}}-\frac{3uv}{n}\right)+\frac{BD}{24}\cdot\frac{l^{2}+4}{l}\\
			&\ +\frac{2A^{3}C-3A^{2}C^{2}+2AC^{3}}{12}\left(\frac{u^{3}}{n^{2}}+\frac{3u^{3}v}{n}\right)-\frac{2B^{3}D-3B^{2}D^{2}+2BD^{3}}{12l} \in\mathbb{Q}\big/\mathbb{Z},\\
			s_{4}&=\frac{C^{2}}{48}\left(u-3uvn-\frac{3u}{n^{2}}-\frac{3uv}{n}\right)+\frac{D^{2}}{48}\cdot\frac{l^{2}+4}{l}+\frac{C^{4}}{24}\left(\frac{u^{3}}{n^{2}}+\frac{3u^{3}v}{n}\right)-\frac{D^{4}}{24l}\in\mathbb{Q}\big/\mathbb{Z},\\
			s_{5}&=-\frac{A^{4}}{2}\left(\frac{u^{3}}{n^{2}}+\frac{3u^{3}v}{n}\right)+\frac{B^{4}}{2l}\in\mathbb{Q}\big/\mathbb{Z},\\
			s_{6}&=-A^{3}C\left(\frac{u^{3}}{n^{2}}+\frac{3u^{3}v}{n}\right)+\frac{B^{3}D}{l}\in\mathbb{Q}\big/\mathbb{Z},\\
			s_{7}&=-\frac{A^{2}C^{2}}{2}\left(\frac{u^{3}}{n^{2}}+\frac{3u^{3}v}{n}\right)+\frac{B^{2}D^{2}}{2l}\in\mathbb{Q}\big/\mathbb{Z},\\
			s_{8}&=-AC^{3}\left(\frac{u^{3}}{n^{2}}+\frac{3u^{3}v}{n}\right)+\frac{BD^{3}}{l}\in\mathbb{Q}\big/\mathbb{Z},\\
			s_{9}&=-\frac{C^{4}}{2}\left(\frac{u^{3}}{n^{2}}+\frac{3u^{3}v}{n}\right)+\frac{D^{4}}{2l}\in\mathbb{Q}\big/\mathbb{Z},\\
			S_{10}&=0\in\mathbb{Q}.\\
		\end{alignedat}
		\right.}
	\end{eqnarray}
\end{lem}

\begin{rmk}\label{Remark: M_1^+ mfd, s1 inv}
	In formulae \eqref{Equation: s-inv of spin M_{m,n,l}, (n,l)=1} and \eqref{Equation: s-inv of spin M_{m,n,l}, (m,n)=1} invariant $s_{1}$ share the same expression. We will see from proof that despite the coprime condition required in Lemmata \ref{Lemma: s-inv of spin M_{m,n,l} with gcd(n,l)=1} and \ref{Lemma: s-inv of spin M_{m,n,l} with gcd(m,n)=1}, expression
	\begin{align}\label{Equation: M_1^+ mfd, s1 inv}
		s_{1}=\frac{3m\left(n^{2}-1\right)\left(n^{2}+3\right)}{896n^{2}}+\frac{\left(l^{2}+4\right)^{2}}{896l}-\frac{\mathrm{sgn}(l)}{224}\in\mathbb{Q}\big/\mathbb{Z}
	\end{align}
	is the common formula of $s_{1}$-invariant for any spin manifold in $\mathcal{M}_{1}$.
\end{rmk}

\begin{pf}[of Lemma \ref{Lemma: s-inv of spin M_{m,n,l} with gcd(n,l)=1}]
	Let $V=V_{m,n,l}$ be the associated disc bundle, and $\left(V;A\pi^{*}z+B\pi^{*}\alpha,C\pi^{*}z+D\pi^{*}\alpha\right)$ is a natural $(\mathcal{B},\mathcal{F})$-coboundary of $\left(M;Ap^{*}z+Bp^{*}\alpha,Cp^{*}z+Dp^{*}\alpha\right)$. By Proposition \ref{Proposition: signature of (D,E)} the signature of $V$ is
	$\sigma\left(V,M\right)=\mathrm{sgn}(l)=\pm1$,
	hence our final choice of $(\mathcal{B},\mathcal{F})$-coboundary with vanishing signature is
	$$\left(V;A\pi^{*}z+B\pi^{*}\alpha,C\pi^{*}z+D\pi^{*}\alpha\right)\#\left(-\mathrm{sgn}(l)\mathbb{H}P^{2}\right).$$
	
	Next we compute the $S$-invariant for this coboundary. Note that
	\begin{align*}
		&S_{i}\left(\left(V;A\pi^{*}z+B\pi^{*}\alpha,C\pi^{*}z+D\pi^{*}\alpha\right)\#\left(-\mathrm{sgn}(l)\mathbb{H}P^{2}\right)\right)\\
		=&S_{i}\left(V;A\pi^{*}z+B\pi^{*}\alpha,C\pi^{*}z+D\pi^{*}\alpha\right)-\mathrm{sgn}(l)S_{i}\left({H}P^{2}\right),\ 1\leqslant i\leqslant10.
	\end{align*}
	Hence we compute $S_{i}\left(V;A\pi^{*}z+B\pi^{*}\alpha,C\pi^{*}z+D\pi^{*}\alpha\right)$ and $S_{i}\left({H}P^{2}\right)$ separately. It is easy to compute that
	\begin{align*}
		S_{1}\left(\mathbb{H}P^{2}\right)=4;\quad S_{i}\left(\mathbb{H}P^{2}\right)=0,\ 2\leqslant i\leqslant10.
	\end{align*}
	We first compute the cohomology classes of $\left(V;x,y\right)=\left(V;A\pi^{*}z+B\pi^{*}\alpha,C\pi^{*}z+D\pi^{*}\alpha\right)$:
	\begin{align*}
		p_{1}
		&=\pi^{*}\left(e\cup\left(\frac{m\left(n^{2}-3\right)}{n^{2}}\alpha+\frac{n^{2}+3}{n}\beta+\frac{l^{2}+4}{l}\gamma\right)\right),\\
		x^{2}
		&=\pi^{*}\left(e\cup\left(-\left(\frac{A^{2}q^{2}m}{n^{2}}+\frac{2ABq}{n}\right)\alpha+\frac{A^{2}q^{2}}{n}\beta+\frac{A^{2}p^{2}}{l}\gamma\right)\right),\\
		xy
		&=\pi^{*}\left(e\cup\left(-\left(\frac{ACq^{2}m}{n^{2}}+\frac{(AD+BC)q}{n}\right)\alpha+\frac{ACq^{2}}{n}\beta+\frac{A^{2}p^{2}}{l}\gamma\right)\right),\\
		y^{2}
		&=\pi^{*}\left(e\cup\left(-\left(\frac{C^{2}q^{2}m}{n^{2}}+\frac{2CDq}{n}\right)\alpha+\frac{C^{2}q^{2}}{n}\beta+\frac{A^{2}p^{2}}{l}\gamma\right)\right).
	\end{align*}
	Hence the monomial $8$-dimensional characteristic numbers of $\left(V;x,y\right)$ are
	\begin{eqnarray}\label{Equation: 8d monomial char num of M_{m,n,l}, (n,l)=1}
		\small
		{\left\{
		\begin{alignedat}{2}
			p_{1}^{2}&=\frac{3m\left(n^{2}-1\right)\left(n^{2}+3\right)}{n^{2}}+\frac{\left(l^{2}+4\right)^{2}}{l},\\
			x^{2}\cdot p_{1}&=A^{2}\left(\frac{q^{2}m\left(n^{2}-3\right)}{n}+\frac{p^{2}\left(l^{2}+4\right)}{l}\right)-\frac{2ABq\left(n^{2}+3\right)}{n},\\
			xy\cdot p_{1}&=AC\left(\frac{q^{2}m\left(n^{2}-3\right)}{n}+\frac{p^{2}\left(l^{2}+4\right)}{l}\right)-\frac{(AD+BC)q\left(n^{2}+3\right)}{n},\\
			y^{2}\cdot p_{1}&=C^{2}\left(\frac{q^{2}m\left(n^{2}-3\right)}{n}+\frac{p^{2}\left(l^{2}+4\right)}{l}\right)-\frac{2CDq\left(n^{2}+3\right)}{n},\\
			\left(x^{2}\right)^{2}&=A^{4}\left(\frac{p^{4}}{l}-\frac{q^{4}m}{n^{2}}\right)-4A^{3}B\cdot\frac{q^{3}}{n},\\
			x^{2}\cdot xy&=A^{3}C\left(\frac{p^{4}}{l}-\frac{q^{4}m}{n^{2}}\right)-A^{2}\left(AD+3BC\right)\cdot\frac{q^{3}}{n},\\
			x^{2}\cdot y^{2}=\left(xy\right)^{2}
			&=A^{2}C^{2}\left(\frac{p^{4}}{l}-\frac{q^{4}m}{n^{2}}\right)-2AC(AD+BC)\cdot\frac{q^{3}}{n},\\
			\left(x^{2}\right)^{2}&=AC^{3}\left(\frac{p^{4}}{l}-\frac{q^{4}m}{n^{2}}\right)-C^{2}\left(3AD+BC\right)\cdot\frac{q^{3}}{n},\\
			\left(x^{2}\right)^{2}&=C^{4}\left(\frac{p^{4}}{l}-\frac{q^{4}m}{n^{2}}\right)-4C^{3}D\cdot\frac{q^{3}}{n}.
		\end{alignedat}
		\right.}
	\end{eqnarray}

	The $S$-invariants of $\left(V;A\pi^{*}z+B\pi^{*}\alpha,C\pi^{*}z+D\pi^{*}\alpha\right)$ are rational linear combinations of these monomial $8$-dimensional characteristic numbers (see Theorem \ref{Theorem: Polarized diffeom inv for polarized E_{1}^{+}-mfds} and the paragraphs before). Combining the expressions
	\begin{align*}
		s_{i}\left(M;Ap^{*}z+Bp^{*}\gamma,Cp^{*}z+Dp^{*}\alpha\right)&=\frac{1}{d_{i}}\left(S_{i}\left(V;A\pi^{*}z+B\pi^{*}\alpha,C\pi^{*}z+D\pi^{*}\alpha\right)-S_{i}\left(\mathbb{H}P^{2}\right)\right)\in\mathbb{Q}\big/\mathbb{Z},\ 1\leqslant i\leqslant9\\
		S_{10}\left(M;Ap^{*}z+Bp^{*}\gamma,Cp^{*}z+Dp^{*}\alpha\right)&=S_{10}\left(V;A\pi^{*}z+B\pi^{*}\alpha,C\pi^{*}z+D\pi^{*}\alpha\right)-S_{10}\left(\mathbb{H}P^{2}\right)\in\mathbb{Q},
	\end{align*}
	we obtain formula \eqref{Equation: s-inv of spin M_{m,n,l}, (n,l)=1}.
\end{pf}

\begin{pf}[of Lemma \ref{Lemma: s-inv of spin M_{m,n,l} with gcd(m,n)=1}]
	The proof is parallel to that of Lemma \ref{Lemma: s-inv of spin M_{m,n,l} with gcd(n,l)=1}. Let $V=V_{m,n,l}$ be the associated disc bundle. A natural $(\mathcal{B},\mathcal{F})$-coboundary of $\left(M;Ap^{*}z+Bp^{*}\gamma,Cp^{*}z+Dp^{*}\gamma\right)$ is $\left(V;A\pi^{*}z+B\pi^{*}\gamma,C\pi^{*}z+D\pi^{*}\gamma\right)$. By Proposition \ref{Proposition: signature of (D,E)} we have $\sigma(V,M)=\mathrm{sgn}(l)=\pm1$. Hence $\left(V;A\pi^{*}z+B\pi^{*}\gamma,C\pi^{*}z+D\pi^{*}\gamma\right)\#\left(-\mathrm{sgn}(l)\mathbb{H}P^{2}\right)$ is a $(\mathcal{B},\mathcal{F})$-coboundary with vanishing signature, and it suffices to compute $S_{i}\left(V;A\pi^{*}z+B\pi^{*}\gamma,C\pi^{*}z+D\pi^{*}\gamma\right)$.
	We first compute the cohomology classes of $\left(V;x,y\right)=\left(V;A\pi^{*}z+B\pi^{*}\gamma,C\pi^{*}z+D\pi^{*}\gamma\right)$:
	\begin{align*}
		p_{1}
		&=\pi^{*}\left(e\cup\left(\frac{m\left(n^{2}-3\right)}{n^{2}}\alpha+\frac{n^{2}+3}{n}\beta+\frac{l^{2}+4}{l}\gamma\right)\right),\\
		x^{2}
		&=\pi^{*}\left(e\cup\left(-\left(\frac{A^{2}u}{n^{2}}+\frac{A^{2}uv}{n}\right)\alpha+\frac{A^{2}u^{2}}{n}\beta+\frac{B^{2}}{l}\gamma\right)\right),\\
		xy
		&=\pi^{*}\left(e\cup\left(-\left(\frac{ACu}{n^{2}}+\frac{ACuv}{n}\right)\alpha+\frac{ACu^{2}}{n}\beta+\frac{BD}{l}\gamma\right)\right),\\
		y^{2}
		&=\pi^{*}\left(e\cup\left(-\left(\frac{C^{2}u}{n^{2}}+\frac{C^{2}uv}{n}\right)\alpha+\frac{C^{2}u^{2}}{n}\beta+\frac{D^{2}}{l}\gamma\right)\right).
	\end{align*}
	Hence the monomial $8$-dimensional characteristic numbers of $\left(V;x,y\right)$ are
	\begin{eqnarray}\label{Equation: 8d monomial char num of M_{m,n,l}, (m,n)=1}
		\small
		{\left\{
			\begin{alignedat}{2}
				p_{1}^{2}&=\frac{3m\left(n^{2}-1\right)\left(n^{2}+3\right)}{n^{2}}+\frac{\left(l^{2}+4\right)^{2}}{l},\\
				x^{2}\cdot p_{1}&=A^{2}\left(u-3uvn-\frac{3u}{n^{2}}-\frac{3uv}{n}\right)+B^{2}\cdot\frac{l^{2}+4}{l},\\
				xy\cdot p_{1}&=AC\left(u-3uvn-\frac{3u}{n^{2}}-\frac{3uv}{n}\right)+BD\cdot\frac{l^{2}+4}{l},\\
				y^{2}\cdot p_{1}&=C^{2}\left(u-3uvn-\frac{3u}{n^{2}}-\frac{3uv}{n}\right)+D^{2}\cdot\frac{l^{2}+4}{l},\\
				\left(x^{2}\right)^{2}&=-A^{4}u^{3}\left(\frac{1}{n^{2}}+\frac{3v}{n}\right)+\frac{B^{4}}{l},\\
				x^{2}\cdot xy&=-A^{3}Cu^{3}\left(\frac{1}{n^{2}}+\frac{3v}{n}\right)+\frac{B^{3}D}{l},\\
				x^{2}\cdot y^{2}=\left(xy\right)^{2}
				&=-A^{2}C^{2}u^{3}\left(\frac{1}{n^{2}}+\frac{3v}{n}\right)+\frac{B^{2}D^{2}}{l},\\
				\left(x^{2}\right)^{2}&=-AC^{3}u^{3}\left(\frac{1}{n^{2}}+\frac{3v}{n}\right)+\frac{BD^{3}}{l},\\
				\left(x^{2}\right)^{2}&=-C^{4}u^{3}\left(\frac{1}{n^{2}}+\frac{3v}{n}\right)+\frac{D^{4}}{l}.
			\end{alignedat}
			\right.}
	\end{eqnarray}
	
	The $S$-invariants of $\left(V;A\pi^{*}z+B\pi^{*}\gamma,C\pi^{*}z+D\pi^{*}\gamma\right)$ are rational linear combinations of these monomial $8$-dimensional characteristic numbers. Combining the expressions
	\begin{align*}
		s_{i}\left(M;Ap^{*}z+Bp^{*}\gamma,Cp^{*}z+Dp^{*}\gamma\right)&=\frac{1}{d_{i}}\left(S_{i}\left(V;A\pi^{*}z+B\pi^{*}\gamma,C\pi^{*}z+D\pi^{*}\gamma\right)-S_{i}\left(\mathbb{H}P^{2}\right)\right)\in\mathbb{Q}\big/\mathbb{Z},\ 1\leqslant i\leqslant9\\
		S_{10}\left(M;Ap^{*}z+Bp^{*}\gamma,Cp^{*}z+Dp^{*}\gamma\right)&=S_{10}\left(V;A\pi^{*}z+B\pi^{*}\gamma,C\pi^{*}z+D\pi^{*}\gamma\right)-S_{10}\left(\mathbb{H}P^{2}\right)\in\mathbb{Q},
	\end{align*}
	we obtain formula \eqref{Equation: s-inv of spin M_{m,n,l}, (m,n)=1}.
\end{pf}

\subsection{Partial classification of spin manifolds in $\mathcal{M}_{1}$}\label{Section: partial classification of spin mfds in M_1}

In this section we apply the results in Sections \ref{Section: coh ring & char classes of mfds in M_1} and \ref{Section: compute s inv of some spin mfds in M_1} to give a partial classification of spin manifolds in $\mathcal{M}_{1}$, proving Theorem \ref{Theorem: Partial classification of mfds in M_1} for spin manifolds.

By Proposition \ref{Proposition: M1 mfds, first classification comparing H^4} the $\mathcal{M}_{1}$-manifold pairs $\left(M_{m,n,l},M_{\overline{m},\overline{n},\overline{l}}\right)$ that we wish to distinguish satisfy $\overline{n}=\pm n$ and $\overline{l}=\pm l$. We will see that in these cases $\overline{l}=-l$ cannot happen. In this section we first compare their invariants $s_{1}$ and obtain Lemma \ref{Lemma: M_1 spin mfds, n bar=pm n, l bar =-l,sufficient condition}. Then we prove each statement of Theorem \ref{Theorem: Partial classification of mfds in M_1} for spin manifolds.

\begin{lem}\label{Lemma: M_1 spin mfds, n bar=pm n, l bar =-l,sufficient condition}
	Suppose we have integers $m$, $n$, $l$ and $\overline{m}$, $\overline{n}$ such that $\overline{n}=\pm n\neq0$, $l\neq0$, $\gcd(m,n,l)=\gcd\left(\overline{m},\overline{n},l\right)=1$ and $(m,n,l)\equiv\left(\overline{m},\overline{n},l\right)\equiv(0,1,0)\ \mathrm{mod}\ 2$.
	Then $M=M_{m,n,l}$ and $\overline{M}=M_{\overline{m},\overline{n},l}$ are diffeomorphic only if 
		\begin{align}\label{Equation: M_1^+ mfds, n bar=pm n, l bar=l, sol}
			\overline{m}\equiv m\ \mathrm{mod}\ 2^{\lambda_{2}(n)}7^{\lambda_{7}(n)}\left(\frac{n}{3^{\mu_{3}(n)}}\right)^{2}.
		\end{align} 
\end{lem}

\begin{pf}
	$M$ and $\overline{M}$ are diffeomorphic only if they have the same invariant $s_{1}$. By Remark \ref{Remark: M_1^+ mfd, s1 inv} their $s_{1}$-invariants have parallel expression. Let them be equal. Since $\overline{n}=\pm n$ we have $\overline{n}^{2}=n^{2}$, and we obtain the equation
	\begin{align}\label{Equation: M_1^+ mfd, n bar=pm n, l bar=l, compare s1 inv}
		3\left(\overline{m}-m\right)\left(n^{2}-1\right)\left(n^{2}+3\right)\equiv0\ \mathrm{mod}\ 896n^{2}.
	\end{align}
	We view this as an equation in $\overline{m}$, $m$ and set $n=3^{b}7^{c}n_{0}$, where $b$, $c$ are non-negative integers and $n_{0}$ is coprime to $2$, $3$ and $7$. By Chinese remainder theorem it is equivalent to solve the equation set
	$$
		\left\{
			\begin{alignedat}{2}
				3\left(\overline{m}-m\right)\left(n^{2}-1\right)\left(n^{2}+3\right)&\equiv0\ \mathrm{mod}\ 2^{7},\\
				3\left(\overline{m}-m\right)\left(n^{2}-1\right)\left(n^{2}+3\right)&\equiv0\ \mathrm{mod}\ 3^{2b},\\
				3\left(\overline{m}-m\right)\left(n^{2}-1\right)\left(n^{2}+3\right)&\equiv0\ \mathrm{mod}\ 7^{2c+1},\\
				3\left(\overline{m}-m\right)\left(n^{2}-1\right)\left(n^{2}+3\right)&\equiv0\ \mathrm{mod}\ n_{0}^{2}.
			\end{alignedat}
		\right.
	$$
	
	We begin with the $2$-primary part. By assumption $\overline{m}$ and $m$ are even. Since $n$ is odd we have $\alpha\left(2,n^{2}+3\right)=2$ and $\alpha\left(2,n^{2}-1\right)\geqslant3$. Moreover, $\alpha\left(2,n^{2}-1\right)=3$ if and only if $n\equiv3,5\ \mathrm{mod}\ 8$. Hence the solution depends on the value of $n\ \mathrm{mod}\ 8$:
	\begin{compactenum}
		\item When $n\equiv1,7\ \mathrm{mod}\ 8$ the equation is automatically true.
		\item When $n\equiv3,5\ \mathrm{mod}\ 8$ the equation has the solution
		$\overline{m}\equiv m\ \mathrm{mod}\ 4.$
	\end{compactenum}
	We summarize the above as
	$$\overline{m}\equiv m\ \mathrm{mod}\ 2^{\lambda_{2}(n)}.$$
	
	Next we consider $3$-primary part equation. Its solution depends on the value of $b$:
	\begin{compactenum}
		\item When $b=0$ the equation automatically holds.
		\item When $b\geqslant1$ the equation has the solution
		$\overline{m}\equiv m\ \mathrm{mod}\ 3^{2b-2}.$
	\end{compactenum}
	We summarize the above as
	$$\overline{m}\equiv m\ \mathrm{mod}\ 3^{2\left(b-\mu_{3}(n)\right)}.$$
	
	Then we consider the equation with modulus $7^{2c+1}$. Its solution depends on the value of $c$:
	\begin{compactenum}
		\item When $c=0$ the solution depends further on the value of $n\ \mathrm{mod}\ 7$:
		\begin{compactenum}
			\item When $n\equiv1,2,5,6\ \mathrm{mod}\ 7$ so that $\left(n^{2}-1\right)\left(n^{2}+3\right)\equiv0\ \mathrm{mod}\ 7$, the equation automatically holds.
			\item When $n\equiv3,4\ \mathrm{mod}\ 7$ so that $\left(n^{2}-1\right)\left(n^{2}+3\right)$ is not divided by $7$, the equation has solution
			$\overline{m}\equiv m\ \mathrm{mod}\ 7.$
		\end{compactenum}
		\item When $c\geqslant1$ the equation has solution
		$\overline{m}\equiv m\ \mathrm{mod}\ 7^{2c+1}.$
	\end{compactenum}
	We summarize the above as
	$$\overline{m}\equiv m\ \mathrm{mod}\ 7^{2c+\lambda_{7}(n)}.$$
	
	Finally the equation with modulus $n_{0}^{2}$ has solution
	$\overline{m}\equiv m\ \mathrm{mod}\ n_{0}^{2}.$
	In summary, when $\overline{l}=l$ manifolds $M$ and $\overline{M}$ are diffeomorphic only if formula \eqref{Equation: M_1^+ mfds, n bar=pm n, l bar=l, sol} is true.
\end{pf}

\begin{pf}[of Theorem \ref{Theorem: Partial classification of mfds in M_1}, Statement 1 for spin manifolds]
	Let $n$, $l$ and $\overline{n}$, $\overline{l}$ be non-zero integers with $\gcd(n,l)=\gcd\left(\overline{n},\overline{l}\right)=1$ and $(n,l)\equiv\left(\overline{n},\overline{l}\right)\equiv(1,0)\ \mathrm{mod}\ 2$. We are determining when $M=M_{0,n,l}$ and $\overline{M}=M_{0,\overline{n},\overline{l}}$ are diffeomorphic. By Proposition \ref{Proposition: M1 mfds, first classification comparing H^4}, Statement 1 they are diffeomorphic only if $\overline{n}=\pm n$ and $\overline{l}=\pm l$. Hence it remains to determine the signs.
	
	We first sketch how to determine the signs from $s$-invariants. Suppose ${M}\xrightarrow{\varphi}\overline{M}$ is a diffeomorphism. Then $\varphi$ induces a polarized diffeomorphism
	$$\left({M};A{p}^{*}{z}+B{p}^{*}\alpha,C{p}^{*}{z}+D{p}^{z}\alpha\right)\xrightarrow{\varphi}\left(\overline{M};\overline{p}^{*}\overline{z},\overline{p}^{*}\alpha\right)$$
	for some $g=\begin{pmatrix}
		A & C \\ B & D
	\end{pmatrix}\in GL(2,\mathbb{Z})$ and a tangential isomorphism $H^{*}\left(\overline{M}\right)\xrightarrow{\varphi^{*}}H^{*}(M)$. Since $\varphi^{*}$ is a tangential isomorphism $g$ satisfies certain equation $\Phi(g)=0$. Then $M$ and $\overline{M}$ are diffeomorphic if and only if equation
	\begin{eqnarray}\label{Equation: Phi(g)=0, s(M;(x,y)g)=s(M bar; x bar, y bar)}
		\left\{
		\begin{alignedat}{2}
			\Phi(g)&=0,\\
			s\left({M};\left(p^{*}z,p^{*}\alpha\right)g\right)&=s\left(\overline{M};\overline{p}^{*}\overline{z},\overline{p}^{*}\alpha\right),
		\end{alignedat}
		\right.
	\end{eqnarray}
	viewed as an equation in $g\in GL(2,\mathbb{Z})$ parametrized by $n$, $l$ and $\overline{n}$, $\overline{l}$, has a solution in $GL(2,\mathbb{Z})$.
	Comparison of invariant $s_{5}$ excludes the possibility $\overline{l}=-l$, and we find $-{A}={D}=1$, ${B}={C}=0$ is a solution for any $\overline{n}=\pm n$, $\overline{l}=l$. Hence we conclude $M_{0,n,l}$ and $M_{0,-n,l}$ are diffeomorphic.
	
	\noindent\textbf{Polarization and tangential isomorphism}
	
	Set $\psi=\varphi^{-1}$ and
	$$\varphi^{*}\begin{pmatrix}
		\overline{p}^{*}\overline{z} & \overline{p}^{*}\alpha
	\end{pmatrix}=\begin{pmatrix}
		{p}^{*}{z} & {p}^{*}\alpha
	\end{pmatrix}\begin{pmatrix}
		{A} & {C} \\ {B} & {D}
	\end{pmatrix},\quad \psi^{*}\begin{pmatrix}
	{p}^{*}{z} & {p}^{*}\alpha
	\end{pmatrix}=\begin{pmatrix}
	\overline{p}^{*}\overline{z} & \overline{p}^{*}\alpha
	\end{pmatrix}\begin{pmatrix}
		{\varepsilon} {D} & -{\varepsilon} {C} \\ -{\varepsilon} {B} & {\varepsilon} {A}
	\end{pmatrix}.$$
	Here $\begin{pmatrix}
		{A} & {C} \\ {B} & {D}
	\end{pmatrix}\in GL(2,\mathbb{Z})$ and ${\varepsilon}={A}{D}-{B}{C}=\pm1$. Then it is straightforward to compute that
	\begin{align*}
		\left\{
		\begin{alignedat}{2}
			\varphi^{*}\left(\left(\overline{p}^{*}\overline{z}\right)^{2}\right)&={A}^{2}\left({p}^{*}{z}\right)^{2}+2{A}{B}\left({p}^{*}{z}\right)\left({p}^{*}\alpha\right),\\
			\varphi^{*}\left(\left(\overline{p}^{*}{z}\right)\left(\overline{p}^{*}\alpha\right)\right)&={A}{C}\left({p}^{*}{z}\right)^{2}+\left({A}{D}+{B}{C}\right)\left({p}^{*}{z}\right)\left({p}^{*}\alpha\right),\\
			\varphi^{*}\left(\left(\overline{p}^{*}\alpha\right)^{2}\right)&={C}^{2}\left({p}^{*}{z}\right)^{2}+2{C}{D}\left({p}^{*}{z}\right)\left({p}^{*}\alpha\right).
		\end{alignedat}
		\right.
	\end{align*}
	Recall that
	\begin{align*}
		H^{4}\left({M}\right)&=\mathbb{Z}\big/|nl|\left\{\left({p}^{*}{z}\right)^{2}\right\}\oplus\mathbb{Z}\big/d\left\{\left({p}^{*}\overline{z}\right)\left({p}^{*}\alpha\right)\right\},\\
		p_{1}\left({M}\right)&=\left(3l^{2}+4n^{2}\right)\left({p}^{*}{z}\right)^{2},
	\end{align*} 
	and we also have the bar counterpart.
	
	Now we analyze the restrictions given by tangential isomorphisms. By definition $\varphi^{*}$ and $\psi^{*}$ should be isomorphisms on the fourth cohomology groups and preserve the first Pontryagin classes.
	We begin with $\left(p^{*}\alpha\right)^{2}$. Since it vanishes we have
	$\varphi^{*}\left(\left(\overline{p}^{*}\alpha\right)^{2}\right)={C}^{2}\left({p}^{*}{z}\right)^{2}+2{C}{D}\left({p}^{*}{z}\right)\left({p}^{*}\alpha\right)=0$,
	namely
	$$
	\left\{
	\begin{alignedat}{2}
		{C}^{2}&\equiv0\ \mathrm{mod}\ |nl|,\\
		2{C}{D}&\equiv0\ \mathrm{mod}\ |n|.
	\end{alignedat}
	\right.
	$$
	Recall that $n$ is odd and coprime to $l$. Also we have ${A}{D}-{B}{C}={\varepsilon}=\pm1$. Hence we obtain
	$$
	\left\{
	\begin{alignedat}{2}
		{C}&\equiv0\ \mathrm{mod}\ |n|,\\
		\gcd\left(n,{D}\right)&=\gcd\left(n,{A}\right)=1;\\
		{C}&\equiv0\ \mathrm{mod}\ |l|.
	\end{alignedat}
	\right.
	$$
	
	The fourth cohomology groups are both isomorphic to $\mathbb{Z}\big/|l|\oplus\left(\mathbb{Z}\big/|n|\right)^{2}$ and $\gcd\left(n,l\right)=1$. Hence we can study the $\mathbb{Z}\big/|l|$-summand and $\left(\mathbb{Z}\big/|n|\right)^{2}$-summand separately. We begin with the $\mathbb{Z}\big/|l|$-summand of fourth cohomology groups. Since
	$n\left({p}^{*}{z}\right)^{2}$, $n\left(\overline{p}^{*}\overline{z}\right)^{2}$ both generate $\mathbb{Z}\big/|l|$-summands and
	$\varphi^{*}\left(n\left(\overline{p}^{*}\overline{z}\right)^{2}\right)={A}^{2}n\left({p}^{*}{z}\right)^{2}$,
	we must have
	$\gcd\left(l,{A}^{2}\right)=1,$
	namely
	$\gcd\left(l,{A}\right)=1.$
	By symmetry we also have
	$\gcd\left(l,{D}\right)=1.$
	Next we consider the $\left(\mathbb{Z}\big/|n|\right)^{2}$-summand of fourth cohomology groups. Note that $\left(l\left(p^{*}z\right)^{2},\left(p^{*}z\right)\left(p^{*}\alpha\right)\right)$ is a $\mathbb{Z}\big/|n|$-basis for the $\left(\mathbb{Z}\big/|n|\right)^{2}$-summand of $H^{4}(M)$, and the bar counterpart is also true. Since
	$$
	\varphi^{*}\begin{pmatrix}
		l\left(\overline{p}^{*}\overline{z}\right)^{2} & \left(\overline{p}^{*}\overline{z}\right)\left(\overline{p}^{*}\alpha\right)
	\end{pmatrix}=\begin{pmatrix}
		l\left({p}^{*}{z}\right)^{2} & \left({p}^{*}{z}\right)\left({p}^{*}\alpha\right)
	\end{pmatrix}\begin{pmatrix}
		{A}^{2} & \frac{{A}{C}}{l} \\
		2{A}{B}|l| & {A}{D}+{B}{C}
	\end{pmatrix},
	$$
	we must have
	$$
	\left\{
	\begin{alignedat}{2}
		{A}{C}&\equiv0\ \mathrm{mod}\ |l|,\\
		\begin{pmatrix}
			{A}^{2} & \frac{{A}{C}}{|l|} \\
			2{A}{B}|l| & {A}{D}+{B}{C}
		\end{pmatrix}&\in GL\left(2,\mathbb{Z}\big/|n|\right).
	\end{alignedat}
	\right.
	$$
	We solve the equation above and obtain
	${C}\equiv0\ \mathrm{mod}\ |l|.$
	
	Finally we consider the first Pontryagin classes. Since $\varphi^{*}p_{1}\left(\overline{M}\right)=p_{1}(M)$ we have
	$$
	\left\{
	\begin{alignedat}{2}
		{A}^{2}\left(3l^{2}+4n^{2}\right)&\equiv3l^{2}+4n^{2}\ \mathrm{mod}\ |nl|,\\
		2{A}{B}\left(3l^{2}+4n^{2}\right)&\equiv0\ \mathrm{mod}\ |n|.
	\end{alignedat}
	\right.
	$$
	We solve the equation above and obtain
	$$
	\left\{
	\begin{alignedat}{2}
		3\left({A}^{2}-1\right)&\equiv0\ \mathrm{mod}\ |n|,\\
		3{B}&\equiv0\ \mathrm{mod}\ |n|;\\
		4\left({A}^{2}-1\right)&\equiv0\ \mathrm{mod}\ |l|.
	\end{alignedat}
	\right.
	$$
	By symmetry we also have
	$$
	\left\{
	\begin{alignedat}{2}
		3\left({D}^{2}-1\right)&\equiv0\ \mathrm{mod}\ |n|,\\
		4\left({D}^{2}-1\right)&\equiv0\ \mathrm{mod}\ |l|.
	\end{alignedat}
	\right.
	$$
	
	In summary, if ${M}\xrightarrow{\varphi}\overline{M}$ is a diffeomorphism with 
	$$\varphi^{*}\begin{pmatrix}
		\overline{p}^{*}\overline{z} & \overline{p}^{*}\alpha
	\end{pmatrix}=\begin{pmatrix}
		{p}^{*}{z} & {p}^{*}\alpha
	\end{pmatrix}g,\quad g=\begin{pmatrix}
		{A} & {C} \\ {B} & {D}
	\end{pmatrix}\in GL(2,\mathbb{Z}),$$
	then equation $\Phi(g)=0$ reads
	\begin{align}
		\left\{
		\begin{alignedat}{2}
			{C}&\equiv0\ \mathrm{mod}\ |n|,\\
			\gcd\left(n,{A}\right)&=\gcd\left(n,{D}\right)=1,\\
			3\left({A}^{2}-1\right)&\equiv 3{B}\equiv3\left({D}^{2}-1\right)\equiv0\ \mathrm{mod}\ |n|;
		\end{alignedat}
		\right.\quad 
		\left\{
		\begin{alignedat}{2}
			{C}&\equiv0\ \mathrm{mod}\ |l|,\\
			\gcd(l,A)&=\gcd(l,D)=1,\\
			4\left({A}^{2}-1\right)&\equiv4\left({D}^{2}-1\right)\equiv0\ \mathrm{mod}\ |l|
		\end{alignedat}
		\right.
	\end{align}
	
	\noindent\textbf{Compare invariant $s_{5}$ and deduce $\overline{l}=l$}
	
	We set $m=0$ in formula \eqref{Equation: s-inv of spin M_{m,n,l}, (n,l)=1} and obtain the $s$-invariant of $\left(M;Ap^{*}z+Bp^{*}\alpha,Cp^{*}z+Dp^{*}\alpha\right)$:
	\begin{eqnarray}\label{Equation: s-inv of spin M_{0,n,l}, (n,l)=1}
		\small\left\{
		\begin{alignedat}{2}
			s_{1}&=\frac{\left(l^{2}+4\right)^{2}}{896l}-\frac{\mathrm{sgn}(l)}{224}\in\mathbb{Q}\big/\mathbb{Z},\\
			s_{2}&=\frac{A^{2}}{48}\cdot\frac{p^{2}\left(l^{2}+4\right)}{l}-\frac{AB}{24}\cdot\frac{q\left(n^{2}+3\right)}{n}-\frac{A^{4}}{24}\cdot\frac{p^{4}}{l}+\frac{A^{3}B}{6}\cdot\frac{q^{3}}{n}\in\mathbb{Q}\big/\mathbb{Z},\\
			s_{3}&=\frac{AC}{24}\cdot\frac{p^{2}\left(l^{2}+4\right)}{l}-\frac{AD+BC}{24}\cdot\frac{q\left(n^{2}+3\right)}{n}\\
			&\ -\frac{AC\left(2A^{2}-3AC+2C^{2}\right)}{12}\cdot\frac{p^{4}}{l}+\left(\frac{A^{3}D+BC^{3}}{6}+\frac{AC(A-C)(B-D)}{2}\right)\cdot\frac{q^{3}}{n}\in\mathbb{Q}\big/\mathbb{Z},\\
			s_{4}&=\frac{C^{2}}{48}\cdot\frac{p^{2}\left(l^{2}+4\right)}{l}-\frac{CD}{24}\cdot\frac{q\left(n^{2}+3\right)}{n}-\frac{C^{4}}{24}\cdot\frac{p^{4}}{l}+\frac{C^{3}D}{6}\cdot\frac{q^{3}}{n}\in\mathbb{Q}\big/\mathbb{Z},\\
			s_{5}&=\frac{A^{4}}{2}\cdot\frac{p^{4}}{l}-2A^{3}B\cdot\frac{q^{3}}{n}\in\mathbb{Q}\big/\mathbb{Z},\\
			s_{6}&=A^{3}C\cdot\frac{p^{4}}{l}-A^{2}\left(AD+3BC\right)\cdot\frac{q^{3}}{n}\in\mathbb{Q}\big/\mathbb{Z},\\
			s_{7}&=\frac{A^{2}C^{2}}{2}\cdot\frac{p^{4}}{l}-AC(AD+BC)\cdot\frac{q^{3}}{n}\in\mathbb{Q}\big/\mathbb{Z},\\
			s_{8}&=AC^{3}\cdot\frac{p^{4}}{l}-C^{2}\left(3AD+BC\right)\cdot\frac{q^{3}}{n}\in\mathbb{Q}\big/\mathbb{Z},\\
			s_{9}&=\frac{C^{4}}{2}\cdot\frac{p^{4}}{l}-2C^{3}D\cdot\frac{q^{3}}{n}\in\mathbb{Q}\big/\mathbb{Z},\\
			S_{10}&=0\in\mathbb{Q}.\\
		\end{alignedat}
		\right.
	\end{eqnarray}
	
	Now we compare invariant $s_{5}$ of $\left(\overline{M};A\overline{p}^{*}\overline{z}+B\overline{p}^{*}\alpha,C\overline{p}^{*}\overline{z}+D\overline{p}^{*}\alpha\right)$ and $\left(M;p^{*}z,p^{*}\alpha\right)$.
	Recall that $\overline{n}=\pm n$ and $\overline{l}=\pm l$. Hence in B\'{e}zout's identities $pn+qn=1$ and $\overline{p}\overline{n}+\overline{q}\overline{l}=1$ we may assume $\overline{p}=\mathrm{sgn}(n)\mathrm{sgn}\left(\overline{n}\right)p$ and $\overline{q}=\mathrm{sgn}(l)\mathrm{sgn}\left(\overline{l}\right)q$, since different values of $p$ and $q$ lead the same value of $p^{*}z$ and $s$-invariants. Now we solve
	$$A^{4}{p}^{4}|n|\mathrm{sgn}\left(\overline{l}\right)-4A^{3}B\mathrm{sgn}(n)|l|\equiv p^{4}|n|\mathrm{sgn}(l)\ \mathrm{mod}\ 2|nl|.$$
	We consider the equation with modulus $2|l|$. Recall that $n$ is odd and coprime to $l$, thus $2l$. Meanwhile it follows from $pn+ql=1$ that $p$ is also odd and coprime to $l$, thus $2l$. Hence we obtain
	$$A^{4}\equiv\mathrm{sgn}(l)\mathrm{sgn}\left(\overline{l}\right)\ \mathrm{mod}\ 2|l|.$$
	Suppose $\overline{l}=-l$. Then the equation above reads
	${A}^{4}+1\equiv0\ \mathrm{mod}\ 2|l|.$
	Since $l$ is even and $\gcd\left(l,{A}\right)=1$, we have $A$ is odd, ${A}^{4}+1\equiv0\ \mathrm{mod}\ 4$ and ${A}^{2}-1\equiv0\ \mathrm{mod}\ 4$. Then
	$
	2=\left({A}^{4}+1\right)-\left({A}^{2}+1\right)\left({A}^{2}-1\right)
	$
	is also divided by $4$. This is a contradiction, and $M_{0,n,l}$, $M_{0,\overline{n},\overline{l}}$ are diffeomorphic only if $\overline{l}=l$.
	
	\noindent\textbf{Conclude $M_{0,n,l}$ and $M_{0,-n,l}$ are diffeomorphic}
	
	Next we show that $M=M_{0,n,l}$ and $\overline{M}=M_{0,\overline{n},l}$ are diffeomorphic when $\overline{n}=-n$. In B\'{e}zout's identities $pn+ql=1$ and $\overline{p}\overline{n}+\overline{q}\overline{l}=1$ we set $\overline{p}=-p$ and $\overline{q}=q$. Then $z=-q\beta+p\gamma$ and $\overline{z}=-q\beta-p\gamma$. Consider the homomorphism $H^{2}\left(\overline{M}\right)\xrightarrow{F}H^{2}(M)$ defined by
	\begin{align*}
		F\begin{pmatrix}
			\overline{p}^{*}\overline{z} & \overline{p}^{*}\alpha
		\end{pmatrix}=\begin{pmatrix}
			p^{*}z & p^{*}\alpha
		\end{pmatrix}\begin{pmatrix}
			A & C \\ B & D
		\end{pmatrix},\ \begin{pmatrix}
			A & C \\ B & D
		\end{pmatrix}=\begin{pmatrix}
			-1 & 0 \\ 0 & 1
		\end{pmatrix},
	\end{align*}
	and it is straightforward to verify that
	\begin{compactenum}
		\item $F\left(w_{2}\left(\overline{M}\right)\right)=w_{2}(M)$; here with a slight abuse of notation we also use $F$ to denotes the mod $2$ reduction of isomorphism $H^{2}\left(\overline{M};\mathbb{Z}\big/2\right)\xrightarrow{F}H^{2}\left(M;\mathbb{Z}\big/2\right)$;
		\item $F$ extends to an isomorphism $H^{4}\left(\overline{M}\right)\xrightarrow{F}H^{4}(M)$ with
		\begin{align*}
			F\begin{pmatrix}
				\left(\overline{p}^{*}\overline{z}\right)^{2} & \left(\overline{p}^{*}\overline{z}\right)\left(\overline{p}^{*}\alpha\right)
			\end{pmatrix}=\begin{pmatrix}
				\left(p^{*}z\right)^{2} & \left(p^{*}z\right)\left(p^{*}\alpha\right)
			\end{pmatrix}\begin{pmatrix}
				1 & 0 \\ 0 & -1
			\end{pmatrix};
		\end{align*}
		\item $F\left(p_{1}\left(\overline{M}\right)\right)=p_{1}(M)$.
	\end{compactenum}
	Hence $F$ defines a tangential isomorphism that is compatible with the polarizations $\left(-p^{*}z,p^{*}\alpha\right)$ and $\left(\overline{p}^{*}\overline{z},\overline{p}^{*}\alpha\right)$.
	
	Finally, we apply formula \eqref{Equation: s-inv of spin M_{0,n,l}, (n,l)=1} to $\left(M;-p^{*}z,p^{*}\alpha\right)$ and $\left(\overline{M};\overline{p}^{*}\overline{z},\overline{p}^{*}\alpha\right)$. It is straightforward to see that
	$s\left(M;-p^{*}z,p^{*}\alpha\right)=s\left(\overline{M};\overline{p}^{*}\overline{z},\overline{p}^{*}\alpha\right)$.
	Therefore, $M=M_{0,n,l}$ and $\overline{M}=M_{0,-n,l}$ are diffeomorphic, and there is a diffeomorphism $M\xrightarrow{f}\overline{M}$ such that its induced homomorphism on cohomology ring is exactly $F$. 
	
	In summary, given non-zero integers $\overline{n}$, $n$ and $\overline{l}$, $l$ such that $\gcd(n,l)=\gcd\left(\overline{n},\overline{l}\right)=1$ and $(n,l)\equiv\left(\overline{n},\overline{l}\right)\equiv(1,0)\ \mathrm{mod}\ 2$, manifolds $M_{0,n,l}$ and $M_{0,\overline{n},\overline{l}}$ are diffeomorphic if and only if $\overline{n}=\pm n$ and $\overline{l}=l$. This proves the spin case of Theorem \ref{Theorem: Partial classification of mfds in M_1}, Statement 1.
\end{pf}

\begin{rmk}
	In the case $\overline{m}=m=0$ we can also deduce $\overline{l}=l$ by comparing invariant $s_{1}$, which would be faster. This approach is insufficient to exclude the possibility $\overline{l}=-l$ when $\overline{m}$ or $m$ is non-zero.
\end{rmk}

\begin{pf}[of Theorem \ref{Theorem: Partial classification of mfds in M_1}, Statement 2]
	Let $m$, $n$, $l$ and $\overline{n}$, $\overline{l}$ be non-zero integers with $\gcd(m,n,l)=\gcd\left(\overline{n},\overline{l}\right)=1$ and $(m,n,l)\equiv\left(0,\overline{n},\overline{l}\right)\equiv(0,1,0)\ \mathrm{mod}\ 2$. We are determining when $M=M_{m,n,l}$ and $\overline{M}=M_{0,\overline{n},\overline{l}}$ are diffeomorphic. By the second statement of Proposition \ref{Proposition: M1 mfds, first classification comparing H^4} they are diffeomorphic only if $n$ divides $m$, $\overline{n}=\pm n$ and $\overline{l}=\pm l$. 
	
	Let us first sketch how to further distinguish such pairs of manifolds by $s$-invariants. Suppose ${M}\xrightarrow{\varphi}\overline{M}$ is a diffeomorphism. According to proof of Proposition \ref{Proposition: M1 mfds, first classification comparing H^4}, Statement 2 there is a polarized diffeomorphism
	$$\left({M};\widetilde{A}{p}^{*}{z}+\widetilde{B}{p}^{*}\alpha,\widetilde{C}{p}^{*}{z}+\widetilde{D}{p}^{*}\alpha\right)\xrightarrow{\varphi}\left(\overline{M};\overline{p}^{*}\overline{z},\overline{p}^{*}\alpha\right).$$
	Here $g=\begin{pmatrix}
		\widetilde{A} & \widetilde{C} \\ \widetilde{B} & \widetilde{D}
	\end{pmatrix}\in GL(2,\mathbb{Z})$, and since $\varphi^{*}$ is a tangential isomorphism $g$ satisfies certain equation
	$\Phi(g)=0.$
	Then $M$ and $\overline{M}$ are diffeomorphic if and only if equation \eqref{Equation: Phi(g)=0, s(M;(x,y)g)=s(M bar; x bar, y bar)},
	viewed as an equation in $g\in GL(2,\mathbb{Z})$ parametrized by $m$, $n$, $l$ and $\overline{n}$, $\overline{l}$ has a solution in $GL(2,\mathbb{Z})$. 
	
	Invariant $S_{10}$ vanishes for both manifolds, and it remains to compare invariants $s_{1}\sim s_{9}$ and combine the known result (Formula \eqref{Equation: M_1^+ mfds, m bar=0, m neq0, tangential isom}) from tangential isomorphism. We recall our conventions on parameters before solving the equations. We have assumed $\overline{p}=\mathrm{sgn}(n)\mathrm{sgn}\left(\overline{n}\right)p$ and $\overline{q}=\mathrm{sgn}(l)\mathrm{sgn}\left(\overline{l}\right)q=q$, hence in the expression of $s$-invariant for $\left(\overline{M};\overline{p}^{*}\overline{z},\overline{p}^{*}\alpha\right)$ no bar sign occurs. And always remember in this case $d=|n|=\left|\overline{n}\right|$.
	
	First we compare invariant $s_{5}$ and exclude the possibility that $\overline{l}=-l$. The result of comparing invariant $s_{1}$ is given in Lemma \ref{Lemma: M_1 spin mfds, n bar=pm n, l bar =-l,sufficient condition}.
	By formula \eqref{Equation: M_1^+ mfds, m bar=0, m neq0, tangential isom} we see that invariants $s_{7}\sim s_{9}$ are already equal. From the comparison of invariants $s_{6}$ we obtain further restrictions on $g$. Then we compare $s_{2}\sim s_{4}$, which can be reduced to solve congruence equations with modulus $48$ or $24$. We may assume a priori that $\alpha(2,q)\gg1$ so that $2$-primary parts of equations comparing $s_{2}\sim s_{4}$ automatically holds. The $3$-primary part of equations comparing $s_{3}$ and $s_{4}$ are already true, and the $3$-primary part of equation comparing $s_{2}$ refines the result of formula \eqref{Equation: M_1^+ mfds, m bar=0, m neq0, tangential isom}. Then we complete the proof of Theorem \ref{Theorem: Partial classification of mfds in M_1}, Statemen 2 for spin manifolds.
	
	When we solve equations from $s$-invariants, we would require further divisibility conditions on parameters $p$, $q$ satisfying $pn+ql=1$ and also on the bar counterpart, so that equation could be simplified. This might change value of $p^{*}z$ and thus value of $s$-invariants, while we can still deduce classification from solving the equation. The reason is as follows.
	From the proof of Proposition \ref{Proposition: M1 mfds, H^4, coh ring & char classes}, Statement 3 we see that altering $p$, $q$ induces a base change
	$$
		\begin{pmatrix}
			p^{*}z_{\lambda} & p^{*}\alpha
		\end{pmatrix}=\begin{pmatrix}
			p^{*}z & p^{*}\alpha
		\end{pmatrix}T_{\lambda},\ T_{\lambda}=\begin{pmatrix}
			1 & 0 \\ -\lambda m & 1
		\end{pmatrix}
	$$
	and thus the change polarization
	$$
		\begin{pmatrix}
			p^{*}z_{\lambda} & p^{*}\alpha
		\end{pmatrix}g=\begin{pmatrix}
			p^{*}z & p^{*}\alpha
		\end{pmatrix}g_{\lambda},\ g_{\lambda}=gT_{\lambda}.
	$$
	Hence for any integer $\lambda$, solvabilities of Equation \eqref{Equation: Phi(g)=0, s(M;(x,y)g)=s(M bar; x bar, y bar)} in $g$ and $g_{\lambda}$ are equivalent.
	
	\noindent\textbf{Polarization and tangential isomorphism}
	
	To compute $s$-invariants we shall first equip the manifolds with polarizations, examine tangential isomorphisms and then compare the $s$-invariants of polarized manifolds. We continue to use the notations in the proof of the second statement of Proposition \ref{Proposition: M1 mfds, first classification comparing H^4}. Let ${M}\xrightarrow{\varphi}\overline{M}$ be a diffeomorphism and $\psi=\varphi^{-1}$. Set
	$$\varphi^{*}\begin{pmatrix}
		\overline{p}^{*}\overline{z} & \overline{p}^{*}\alpha
	\end{pmatrix}=\begin{pmatrix}
	p^{*}z & p^{*}\alpha
	\end{pmatrix}\begin{pmatrix}
		\widetilde{A} & \widetilde{C} \\ \widetilde{B} & \widetilde{D}
	\end{pmatrix},\quad \psi^{*}\begin{pmatrix}
	p^{*}z & p^{*}\alpha
	\end{pmatrix}=\begin{pmatrix}
	\overline{p}^{*}\overline{z} & \overline{p}^{*}\alpha
	\end{pmatrix}\begin{pmatrix}
	\widetilde{\varepsilon} \widetilde{D} & -\widetilde{\varepsilon} \widetilde{C} \\ -\widetilde{\varepsilon} \widetilde{B} & \widetilde{\varepsilon} \widetilde{A}
	\end{pmatrix}.$$
	Here $\begin{pmatrix}
		\widetilde{A} & \widetilde{C} \\ \widetilde{B} & \widetilde{D}
	\end{pmatrix}\in GL(2,\mathbb{Z})$ and $\widetilde{\varepsilon}=\widetilde{A}\widetilde{D}-\widetilde{B}\widetilde{C}=\pm1$. Then it is straightforward to compute that
	\begin{align*}
		\left\{
			\begin{alignedat}{2}
				\varphi^{*}\left(\left(\overline{p}^{*}\overline{z}\right)^{2}\right)&=\widetilde{A}^{2}\left({p}^{*}{z}\right)^{2}+2\widetilde{A}\widetilde{B}\left({p}^{*}{z}\right)\left({p}^{*}\alpha\right),\\
				\varphi^{*}\left(\left(\overline{p}^{*}\overline{z}\right)\left(\overline{p}^{*}\alpha\right)\right)&=\widetilde{A}\widetilde{C}\left({p}^{*}{z}\right)^{2}+\left(\widetilde{A}\widetilde{D}+\widetilde{B}\widetilde{C}\right)\left({p}^{*}{z}\right)\left({p}^{*}\alpha\right),\\
				\varphi^{*}\left(\left(\overline{p}^{*}\alpha\right)^{2}\right)&=\widetilde{C}^{2}\left({p}^{*}{z}\right)^{2}+2\widetilde{C}\widetilde{D}\left({p}^{*}{z}\right)\left({p}^{*}\alpha\right).
			\end{alignedat}
		\right.
	\end{align*}
	Recall that
	\begin{align*}
		H^{4}\left(\overline{M}\right)&=\mathbb{Z}\big/d|l|\left\{\left(\overline{p}^{*}\overline{z}\right)^{2}\right\}\oplus\mathbb{Z}\big/d\left\{\left(\overline{p}^{*}\overline{z}\right)\left(\overline{p}^{*}\alpha\right)\right\},\\
		p_{1}\left(\overline{M}\right)&=\left(3|l|^{2}+4d^{2}\right)\left(\overline{p}^{*}\overline{z}\right)^{2};\\
		H^{4}(M)&=\mathrm{Z}\big/d\left\{p^{*}\omega\right\}\oplus\mathbb{Z}\big/d\left\{p^{*}\rho\right\}\oplus\mathbb{Z}\big/|l|\left\{\left(p^{*}\gamma\right)^{2}\right\};
	\end{align*} 
	and under the assumption $u'$ is divided by $d$ we have a neater relation between generator sets $\left\{p^{*}\omega,p^{*}\rho,\left(p^{*}\gamma\right)^{2}\right\}$ and $\left\{\left(p^{*}z\right)^{2},\left(p^{*}z\right)\left(p^{*}\alpha\right)\right\}$. Since $n_{1}=\pm1$ in the relation $u'm_{1}+v'n_{1}=1$ we can assume further that $u'=0$ and $v'=n_{1}$ so that
	\begin{align*}
		\left\{
		\begin{alignedat}{2}
			\left(p^{*}z\right)^{2}&=q^{2}n_{1}\left(p^{*}\omega+m_{1}n_{1}p^{*}\rho\right)+p^{2}\left(p^{*}\gamma\right)^{2},\\
			\left(p^{*}z\right)\left(p^{*}\alpha\right)&=qn_{1}p^{*}\rho,\\
			\left(p^{*}\alpha\right)^{2}&=0.
		\end{alignedat}
		\right.
	\end{align*}
	In particular, we also have
	\begin{align*}
		H^{4}\left({M}\right)&=\mathbb{Z}\big/d|l|\left\{\left({p}^{*}{z}\right)^{2}\right\}\oplus\mathbb{Z}\big/d\left\{\left({p}^{*}{z}\right)\left({p}^{*}\alpha\right)\right\},\\
		p_{1}(M)&=\left(3|l|^{2}+4d^{2}\right)\left(p^{*}z\right)^{2}.
	\end{align*}
	
	The analysis of polarization and tangential isomorphism is exactly the same as the case where $\overline{m}=m=0$ (see the Proof of Theorem \ref{Theorem: Partial classification of mfds in M_1}, Statement 1 for spin manifolds), and it suffices to replace $A$, $B$, $C$, $D$ and $\varepsilon$ by their tilde counterpart and write $d=|n|=\left|\overline{n}\right|$. The results are similar: if $\overline{M}\xrightarrow{\varphi}M$ is a diffeomorphism with 
	$$\varphi^{*}\begin{pmatrix}
		\overline{p}^{*}\overline{z} & \overline{p}^{*}\alpha
	\end{pmatrix}=\begin{pmatrix}
		p^{*}z & p^{*}\alpha
	\end{pmatrix}g,\quad g=\begin{pmatrix}
	\widetilde{A} & \widetilde{C} \\ \widetilde{B} & \widetilde{D}
	\end{pmatrix}\in GL(2,\mathbb{Z}),$$
	then the equation $\Phi(g)=0$ reads
	\begin{align}\label{Equation: M_1^+ mfds, m bar=0, m neq0, tangential isom}
		\left\{
			\begin{alignedat}{2}
				\widetilde{C}&\equiv0\ \mathrm{mod}\ d,\\
				\gcd\left(d,\widetilde{A}\right)&=\gcd\left(d,\widetilde{D}\right)=1,\\
				3\left(\widetilde{A}^{2}-1\right)&\equiv 3\widetilde{B}\equiv3\left(\widetilde{D}^{2}-1\right)\equiv0\ \mathrm{mod}\ d;
			\end{alignedat}
		\right.\quad 
		\left\{
			\begin{alignedat}{2}
				\widetilde{C}&\equiv0\ \mathrm{mod}\ |l|,\\
				\gcd\left(l,\widetilde{A}\right)&=\gcd\left(l,\widetilde{D}\right)=1,\\
				4\left(\widetilde{A}^{2}-1\right)&\equiv4\left(\widetilde{D}^{2}-1\right)\equiv0\ \mathrm{mod}\ |l|
			\end{alignedat}
		\right.
	\end{align}
	
	\noindent\textbf{Compare invariant $s_{5}$ and exclude the case $\overline{l}=-l$}
	
	Now we compare invariant $s_{5}$ and show that the case $\overline{l}=-l$ will never happen. We read the invariants $s_{5}$ of $\left(M;\widetilde{A}p^{*}z+\widetilde{B}p^{*}\alpha,\widetilde{C}p^{*}z+\widetilde{D}p^{*}\alpha\right)$ and $\left(\overline{M};\overline{p}^{*}\overline{z},\overline{p}^{*}\alpha\right)$, let them be equal, and we obtain:
	\begin{align}\label{Equation: M_1^+ mfds M_{m,n,l} & M_{0,n bar, l bar}, inv s_5}
		\widetilde{A}^{4}p^{4}n^{2}\mathrm{sgn}(l)-\widetilde{A}^{4}q^{4}m|l|-4\widetilde{A}^{3}\widetilde{B}q^{3}n|l|\equiv\overline{p}^{4}n^{2}\mathrm{sgn}\left(\overline{l}\right)\ \mathrm{mod}\ 2n^{2}|l|.
	\end{align}
	The equation above with modulus $2|l|$ is simplified as
	$\widetilde{A}^{4}\equiv\mathrm{sgn}\left(\overline{l}\right)\mathrm{sgn}(l)\ \mathrm{mod}\ 2|l|$,
	which must imply $\overline{l}=l$ from the proof of Theorem \ref{Theorem: Partial classification of mfds in M_1}, Statement 1 for spin manifolds. Hence given non-zero integers $m$, $n$, $l$ and $\overline{n}$, $\overline{l}$ with $\gcd(m,n,l)=\gcd\left(\overline{n},\overline{l}\right)=1$ and $(m,n,l)\equiv\left(0,\overline{n},\overline{l}\right)\equiv(0,1,0)\ \mathrm{mod}\ 2$, manifolds $M=M_{m,n,l}$ and $\overline{M}=M_{0,\overline{n},\overline{l}}$ are diffeomorphic only if $n$ divides $m$, $\overline{n}=\pm n$ and $\overline{l}=l$. 
	
	We continue to solve equation \eqref{Equation: M_1^+ mfds M_{m,n,l} & M_{0,n bar, l bar}, inv s_5}. Now we have
	\begin{align*}
		\widetilde{A}^{4}p^{4}n^{2}\mathrm{sgn}(l)-\widetilde{A}^{4}q^{4}m|l|-4\widetilde{A}^{3}\widetilde{B}q^{3}n|l|\equiv{p}^{4}n^{2}\mathrm{sgn}\left({l}\right)\ \mathrm{mod}\ 2n^{2}|l|.
	\end{align*}
	Since $\gcd\left(n^{2},2l\right)=1$, by Chinese remainder theorem it is equivalent to solve the congruence equation set
	\begin{align}\label{Equation: M_1^+ mfds M_{m,n,l} & M_{0,n bar, l}, inv s_5}
		\left\{
		\begin{alignedat}{2}
			\widetilde{A}^{4}p^{4}n^{2}\mathrm{sgn}(l)-\widetilde{A}^{4}q^{4}m|l|-4\widetilde{A}^{3}\widetilde{B}q^{3}n|l|&\equiv{p}^{4}n^{2}\mathrm{sgn}\left({l}\right)\ \mathrm{mod}\ n^{2},\\
			\widetilde{A}^{4}p^{4}n^{2}\mathrm{sgn}(l)-\widetilde{A}^{4}q^{4}m|l|-4\widetilde{A}^{3}\widetilde{B}q^{3}n|l|&\equiv{p}^{4}n^{2}\mathrm{sgn}\left({l}\right)\ \mathrm{mod}\ 2|l|.
		\end{alignedat}
		\right.
	\end{align}
	Recall that $\gcd\left(d,\widetilde{A}\right)=1$, $n$ divides $m$ and $pn+ql=1$. Recall also from formula \eqref{Equation: M_1^+ mfds, m bar=0, m neq0, tangential isom} that $3\widetilde{B}\equiv0\ \mathrm{mod}\ d$. Hence from the equation with modulus $n^{2}$ we obtain
	\begin{align}\label{Equation: M_1^+ mfds M_{m,n,l} & M_{0,n bar, l}, inv s_5 sol, d-primary}
		\widetilde{A}\cdot\frac{m}{n}+\widetilde{B}l\equiv0\ \mathrm{mod}\ d.
	\end{align}
	Then we consider the equation with modulus $2|l|$. In this we obtain
	$\widetilde{A}^{4}-1\equiv0\ \mathrm{mod}\ 2|l|.$
	Meanwhile, by formula \eqref{Equation: M_1^+ mfds, m bar=0, m neq0, tangential isom} we have 
	$4\left(\widetilde{A}^{2}-1\right)\equiv0\ \mathrm{mod}\ |l|.$
	Let us denote $l=2^{a}l_{0}$, where $a$ is a positive integer and $l_{0}$ is odd. We begin with the equations with modulus $\left|l_{0}\right|$. We have $\widetilde{A}^{2}-1\equiv0\ \mathrm{mod}\ \left|l_{0}\right|$, and $\widetilde{A}^{2}+1\equiv0\ \mathrm{mod}\ \left|l_{0}\right|$ if and only if $l_{0}$ divides $2$, namely $l_{0}=\pm1$ as it is odd. Next we consider the $2$-primary part. This implies that $\widetilde{A}^{2}+1$ is always coprime to $l_{0}$. Recall that $\widetilde{A}$ is odd, hence $\alpha\left(2,\widetilde{A}^{2}+1\right)=1$, and it follows from $\widetilde{A}^{4}-1=\left(\widetilde{A}^{2}+1\right)\left(\widetilde{A}^{2}-1\right)$ and $\gcd\left(l_{0},\widetilde{A}^{2}+1\right)=1$ that
	$\widetilde{A}^{2}-1\equiv0\ \mathrm{mod}\ |l|$.
	This refines the result $4\left(\widetilde{A}^{2}-1\right)\equiv0\ \mathrm{mod}\ |l|$ in formula \eqref{Equation: M_1^+ mfds, m bar=0, m neq0, tangential isom}.
	
	\noindent\textbf{Compare invariants $s_{7}\sim s_{9}$}
	
	Next we consider invariants $s_{7}\sim s_{9}$. As for the polarized $\mathcal{E}_{1}^{+}$-manifold $\left(\overline{M};\overline{p}^{*}\overline{z},\overline{p}^{*}\alpha\right)$, invariants $s_{7}\sim s_{9}$ are all already $0\in\mathbb{Q}\big/\mathbb{Z}$. Meanwhile, invariants $s_{7}\sim s_{9}$ of $\left({M};\widetilde{A}{p}^{*}{z}+\widetilde{B}{p}^{*}\alpha,\widetilde{C}{p}^{*}{z}+\widetilde{D}{p}^{*}\alpha\right)$ are also $0\in\mathbb{Q}\big/\mathbb{Z}$. This follows from $C\equiv0\ \mathrm{mod}\ d$, $C\equiv0\ \mathrm{mod} |l|$ (formula \eqref{Equation: M_1^+ mfds, m bar=0, m neq0, tangential isom}), the assumption that $l$ is even, and the known result $n$ divides $m$. Therefore, invariants $s_{7}\sim s_{9}$ of both polarized manifolds are equal.
	
	\noindent\textbf{Compare invariant $s_{6}$}
	
	Now we compare invariant $s_{6}$. Since $l$ divides $\widetilde{C}$, $n$ divides $\widetilde{C}$, $m$ and $pn+ql=1$, we obtain
	\begin{eqnarray}\label{Equation: M_1^+ mfds M_{m,n,l} & M_{0,n bar, l}, inv s_6}
		\widetilde{A}^{3}\widetilde{D}\equiv1\ \mathrm{mod}\ d.
	\end{eqnarray}
	
	\noindent\textbf{Compare invariant $s_{4}$}
	
	Next we compare invariant $s_{4}$. We have
	\begin{eqnarray}\label{Equation: M_1^+ mfds M_{m,n,l} & M_{0,n bar, l}, inv s_4}
		\left.
		\begin{alignedat}{2}
			&\frac{\widetilde{C}^{2}}{n}\cdot q^{2}m\left(n^{2}-3\right)+\frac{\widetilde{C}^{2}}{l}\cdot p^{2}\left(l^{2}+4\right)-2\cdot\frac{\widetilde{C}}{n}\cdot\widetilde{D}q^{2}\left(n^{2}+3\right)\\
			-&2\cdot\frac{\widetilde{C}^{4}}{l}\cdot p^{4}+2\cdot\frac{\widetilde{C}^{4}}{n^{2}}\cdot q^{4}m+8\cdot\frac{\widetilde{C}^{3}}{n}\cdot\widetilde{D}q^{3}\equiv0\ \mathrm{mod}\ 48.
		\end{alignedat}
		\right.
	\end{eqnarray}
	Since $\widetilde{C}$ is divided by $n$ and $l$, in the left hand side each factor of each term is an integer. By Chinese remainder theorem we solve the $2$-primary part and $3$-primary part of equation \eqref{Equation: M_1^+ mfds M_{m,n,l} & M_{0,n bar, l}, inv s_4} separately.
	
	We begin with the $2$-primary part. Since $pn+ql=1$ and $n$ is odd, we may assume a priori $\alpha(2,q)\gg2$ and the $2$-primary part automatically holds.
	Next we consider the $3$-primary part. Note that for any integer $r$ we have $r^{3}\equiv r\ \mathrm{mod}\ 3$. Using this relation the $3$-primary part of equation \eqref{Equation: M_1^+ mfds M_{m,n,l} & M_{0,n bar, l}, inv s_4} can be reduced as
	$$\left(\frac{\widetilde{C}}{n}\right)^{2}q^{2}m\left(n-n^{2}\right)+\frac{\widetilde{C}}{n}\cdot\widetilde{D}\left(q^{2}-q\right)n^{2}\equiv0\ \mathrm{mod}\ 3.$$
	When $n$ is divided by $3$, this equation automatically holds. And when $n$ is not divided by $3$, from the relation $pn+ql=1$ we may assume a priori that $q$ is divided by $3$ and the equation above also automatically holds.
	
	To summarize, the invariants $s_{4}$ of both polarized $\mathcal{E}_{1}^{+}$-manifolds under consideration are always equal.
	
	\noindent\textbf{Compare invariant $s_{2}$}
	
	Then we move to compare invariant $s_{2}$. We have
	\begin{align}\label{Equation: M_1^+ mfds M_{m,n,l} & M_{0,n bar, l}, inv s_2}
		\left.
			\begin{alignedat}{2}
				&\widetilde{A}^{2}q^{2}\cdot\frac{m}{n}\cdot\left(n^{2}-3\right)+\frac{\widetilde{A}^{2}-1}{l}\cdot p^{2}\left(l^{2}+4\right)-2\widetilde{A}\widetilde{B}qn\\
				-&2\widetilde{A}\cdot\frac{3\widetilde{B}}{n}\cdot q-2\cdot\frac{\widetilde{A}^{4}-1}{l}\cdot p^{4}+2\widetilde{A}^{3}q^{3}\cdot\frac{\widetilde{A}q\cdot\frac{m}{n}+4\widetilde{B}}{n}\equiv0\ \mathrm{mod}\ 48.
			\end{alignedat}
		\right.
	\end{align}
	Since $n$ divides $m$, $3\overline{B}\equiv0\ \mathrm{mod}\ d$ (formula \eqref{Equation: M_1^+ mfds, m bar=0, m neq0, tangential isom}), $\widetilde{A}^{2}-1\equiv0\ \mathrm{mod}\ |l|$  and $\widetilde{A}q\cdot\frac{m}{n}+4\widetilde{B}\equiv0\ \mathrm{mod}\ d$ (result of comparing invariant $s_{5}$), in the left hand side each factor of each term is an integer. By Chinese remainder theorem we solve the $2$-primary part and $3$-primary part of equation \eqref{Equation: M_1^+ mfds M_{m,n,l} & M_{0,n bar, l}, inv s_2} separately.
	
	We begin with the $2$-primary part. We may assume a priori $\alpha(2,q)\gg1$ as before. 
	Since $pn+ql=1$ and $l$ is even, we see that $p$ is odd. Hence $p^{2}-1\equiv0\ \mathrm{mod}\ 8$ and $p^{4}-1\equiv0\ \mathrm{mod}\ 16$, and the $2$-primary part equation can be simplified as
	$$\left(\widetilde{A}^{2}-1\right)l-2\cdot\frac{\widetilde{A}^{2}-1}{l}\cdot\left(\widetilde{A}^{2}-1\right)\equiv0\ \mathrm{mod}\ 16.$$
	Since $\widetilde{A}$ is odd we have $\widetilde{A}^{2}-1\equiv0\ \mathrm{mod}\ 8$, and the equation above automatically holds. Therefore, the $2$-primary part of equation \eqref{Equation: M_1^+ mfds M_{m,n,l} & M_{0,n bar, l}, inv s_2} is automatically true.
	
	Next we consider the $3$-primary part. We are solving the equation
	\begin{eqnarray}\label{Equation: M_1^+ mfds M_{m,n,l} & M_{0,n bar, l}, inv s_2, 3 primary}
		\widetilde{A}^{2}q^{2}\cdot\frac{m}{n}\cdot n^{2}+\frac{\widetilde{A}^{2}-1}{l}\cdot p^{2}\left(l^{2}+1\right)+\widetilde{A}\widetilde{B}qn+\frac{\widetilde{A}^{4}-1}{l}\cdot p^{2}-\widetilde{A}q\cdot\frac{\widetilde{A}q\cdot\frac{m}{n}+\widetilde{B}}{n}\equiv0\ \mathrm{mod}\ 3.
	\end{eqnarray}
	Whether $d$ is divided by $3$ affects the solution to equation \eqref{Equation: M_1^+ mfds M_{m,n,l} & M_{0,n bar, l}, inv s_2, 3 primary}.
	\begin{compactenum}
		\item When $n$ is not divided by $3$, from $pn+ql=1$ we may assume a priori that $q$ is divided by $3$ and thus $p$ is not. Then equation \eqref{Equation: M_1^+ mfds M_{m,n,l} & M_{0,n bar, l}, inv s_2, 3 primary} is reduced as
		$$\frac{\widetilde{A}^{2}-1}{l}\cdot\left(l^{2}+\widetilde{A}^{2}-1\right)\equiv0\ \mathrm{mod}\ 3.$$
		The treatment depends further on divisibility of $l$ by $3$:
		\begin{compactenum}
			\item When $3$ divides $l$ the left hand side can be written as $\frac{\widetilde{A}^{2}-1}{l}\cdot l^{2}+\left(\frac{\widetilde{A}^{2}-1}{l}\right)^{2}\cdot l$ and is divided by $3$.
			\item When $3$ does not divide $l$ we have $l^{2}-1\equiv0\ \mathrm{mod}\ 3$, and we are reduced to solve the equation
			$\left(\widetilde{A}^{2}-1\right)\widetilde{A}^{2}\equiv0\ \mathrm{mod}\ 3$,
			which is automatically true.
		\end{compactenum}
		Therefore, when $n$ is not divided by $3$, equation \eqref{Equation: M_1^+ mfds M_{m,n,l} & M_{0,n bar, l}, inv s_2, 3 primary} automatically holds.
		\item When $n$ is divided by $3$ so that $l$ is not, from $pn+ql=1$ we may assume a priori that $p$ is divided by $3$ and thus $q$ is not. Then equation \eqref{Equation: M_1^+ mfds M_{m,n,l} & M_{0,n bar, l}, inv s_2, 3 primary} is reduced as
		$$\widetilde{A}\cdot\frac{\widetilde{A}q\cdot\frac{m}{n}+\widetilde{B}}{n}\equiv0\ \mathrm{mod}\ 3.$$
		Recall that $n$ divides $\widetilde{C}$ and $\widetilde{A}\widetilde{D}-\widetilde{B}\widetilde{C}=\pm1$. Hence $\widetilde{A}$ is not divided by $3$. Moreover, $\widetilde{A}$ is coprime to $3d$. Since $pn+ql=1$ and by assumption $p$ is divided by $3$, we have $ql\equiv0\ \mathrm{mod}\ 3d$. Hence equation \eqref{Equation: M_1^+ mfds M_{m,n,l} & M_{0,n bar, l}, inv s_2, 3 primary} can be further reduced as
		$$\frac{m}{n}\equiv-\widetilde{A}^{-1}\widetilde{B}l\ \mathrm{mod}\ 3d.$$
		Here with a slight abuse of notation we denote by $\widetilde{A}^{-1}$ the multiplicative inverse of $\widetilde{A}$ modulo $3d$.
		
		Recall from formula \eqref{Equation: M_1^+ mfds, m bar=0, m neq0, tangential isom} that $3\widetilde{B}\equiv0\ \mathrm{mod}\ d$, namely $\widetilde{B}\equiv0\ \mathrm{mod}\ \frac{d}{3}$. Hence we may set $d=3^{b}7^{c}d_{0}$ and $\widetilde{B}=\widetilde{B_{1}}\cdot 3^{b-1}7^{c}d_{0}$, where $b$ is a positive integer, $c$ is a non-negative integer and $d_{0}$ is coprime to both $3$ and $7$. Then we obtain
		$$\frac{m}{n}\equiv-\widetilde{A}^{-1}\widetilde{B}_{1}l\cdot3^{b-1}7^{c}d_{0}\ \mathrm{mod}\ 3^{b+1}7^{c}d_{0},$$
		namely,
		\begin{eqnarray}\label{Equation: M_1^+ mfds M_{m,n,l} & M_{0,n bar, l}, inv s_2, 3 primary, 3|n, sol}
			m\equiv0\ \mathrm{mod}\ 3^{2b-1}7^{2c}d_{0}^{2}.
		\end{eqnarray}
		The modulus is exactly $\frac{n^{2}}{3}$.
	\end{compactenum}
	
	In summary, by comparing invariant $s_{2}$ we obtain
	$$m\equiv0\ \mathrm{mod}\ \frac{n^{2}}{3^{\mu_{3}(b)}}.$$
	This result refines what we deduced from tantengial isomorphism and comparison of invariants $s_{1}$, $s_{4}\sim s_{9}$ and $S_{10}$, and combining all of them we obtain
	$$m\equiv0\ \mathrm{mod}\ \frac{2^{\lambda_{2}(n)}7^{\lambda_{7}(n)}n^{2}}{3^{\mu_{3}(b)}}.$$
	
	\noindent\textbf{Compare invariant $s_{3}$ and conclude}
	
	Finally we compare invariant $s_{3}$. We are solving
	\begin{eqnarray}\label{Equation: M_1^+ mfds M_{m,n,l} & M_{0,n bar, l}, inv s_3}
		\left.
			\begin{alignedat}{2}
				&\widetilde{A}\cdot\frac{\widetilde{C}}{n}\cdot q^{2}mn^{2}+\widetilde{A}\cdot\frac{\widetilde{C}}{l}\cdot p^{2}\left(l^{2}+4\right)-\widetilde{B}\cdot\frac{C}{n}\cdot qn^{2}-4\widetilde{A}\left(\widetilde{A}^{2}+\widetilde{C}^{2}\right)\left(\frac{\widetilde{C}}{l}\cdot p^{4}-\frac{\widetilde{C}}{n}\cdot q^{4}\cdot\frac{m}{n}\right)\\
				+&6\widetilde{A}\widetilde{C}\cdot\frac{\widetilde{C}}{l}\cdot p^{4}+4\widetilde{B}\widetilde{C}^{2}\cdot\frac{\widetilde{C}}{n}\cdot q^{3}-\left(\widetilde{A}\widetilde{D}-1\right)qn+4\cdot\frac{\widetilde{A}^{3}\widetilde{D}-1}{n}\cdot q^{3}-\frac{3\left(\widetilde{A}\widetilde{D}-1\right)q}{n}\equiv0\ \mathrm{mod}\ 24.
			\end{alignedat}
		\right.
	\end{eqnarray}
	According to our previous results, all but the last term of left hand side are integers. Hence the last term of the left hand side should also be an integer, and since $pn+ql=1$ we have
	$$3\left(\widetilde{A}\widetilde{D}-1\right)\equiv0\ \mathrm{mod}\ d.$$
	
	We return to equation \eqref{Equation: M_1^+ mfds M_{m,n,l} & M_{0,n bar, l}, inv s_3}. By Chinese remainder theorem we treat the $2$-primary part and $3$-primary part separately. We begin with the $2$-primary part. We assume a priori that $\alpha(2,q)\gg1$ as before. Since $\widetilde{A}$ and $p$ are odd, we have $\widetilde{A}^{2}\equiv1\ \mathrm{mod}\ 8$, $p^{4}\equiv1\ \mathrm{mod}\ 16$ and the $2$-primary part equation is simplified as
	$$\frac{\widetilde{C}}{l}\cdot \left(\frac{l}{2}\right)^{2}-\left(\frac{\widetilde{C}}{l}\right)^{2}\cdot\frac{l}{2}\equiv0\ \mathrm{mod}\ 2,$$
	which is automatically true.
	
	Next we consider the $3$-primary part of equation \eqref{Equation: M_1^+ mfds M_{m,n,l} & M_{0,n bar, l}, inv s_3}. We are reduced to solve the following equation
	\begin{eqnarray}\label{Equation: M_1^+ mfds M_{m,n,l} & M_{0,n bar, l}, inv s_3, 3-primary}
		\widetilde{A}\cdot\frac{\widetilde{C}}{n}\cdot q^{2}\cdot\frac{m}{n}\cdot\left(n^{2}+n+1\right)-\left(\widetilde{A}\widetilde{D}-1\right)qn-\frac{3\left(\widetilde{A}\widetilde{D}-1\right)}{n}\cdot q+\frac{\widetilde{A}^{3}\widetilde{D}-1}{n}\cdot q\equiv0\ \mathrm{mod}\ 3,
	\end{eqnarray}
	and the treatment is different according to divisibility of $n$ by $3$:
	\begin{compactenum}
		\item When $n$ is not divided by $3$, from $pn+ql=1$ we may assume a priori that $q$ is divided by $3$, and equation \eqref{Equation: M_1^+ mfds M_{m,n,l} & M_{0,n bar, l}, inv s_3, 3-primary} is automatically true.
		\item When $n$ is divided by $3$, we have $q$, $l$ and $\widetilde{A}$ are coprime to $3$. Then equation \eqref{Equation: M_1^+ mfds M_{m,n,l} & M_{0,n bar, l}, inv s_3, 3-primary} can be simplified as
		$$\frac{\widetilde{C}}{n}\cdot\frac{m}{n}\equiv\widetilde{A}^{-1}l\left(\frac{3\left(\widetilde{A}\widetilde{D}-1\right)}{n}-\frac{\widetilde{A}^{3}\widetilde{D}-1}{n}\right)\ \mathrm{mod}\ 3.$$
		Recall that $\widetilde{A}\widetilde{D}-\widetilde{B}\widetilde{C}=\widetilde{\varepsilon}=\pm1$ and $n$ divides $\widetilde{C}$, hence
		\begin{align*}
			\frac{3\left(\widetilde{A}\widetilde{D}-1\right)}{n}-\frac{\widetilde{A}^{3}\widetilde{D}-1}{n}&=\frac{3\left(\widetilde{\varepsilon}-1\right)}{n}+3\widetilde{B}\cdot\frac{\widetilde{C}}{n}-\widetilde{A}^{2}\widetilde{B}\cdot\frac{\widetilde{C}}{n}-\frac{\widetilde{A}^{2}\widetilde{\varepsilon}-1}{n}\\
			&\equiv\frac{3\left(\widetilde{\varepsilon}-1\right)}{n}-\widetilde{B}\cdot\frac{\widetilde{C}}{n}-\frac{\widetilde{A}^{2}\widetilde{\varepsilon}-1}{n}\ \mathrm{mod}\ 3,
		\end{align*}
		where we still have
		$$3\left(\widetilde{\varepsilon}-1\right)\equiv\widetilde{A}^{2}\widetilde{\varepsilon}-1\equiv0\ \mathrm{mod}\ d.$$
		If $\widetilde{\varepsilon}$ can take $-1$, then $d$ divides $6$ and thus must be $3$ as $d$ is odd, and we have $\widetilde{A}^{2}+1\equiv0\ \mathrm{mod}\ 3$. However, $\widetilde{A}$ is coprime to $\widetilde{C}$ and thus to $3$ and we shall have $\widetilde{A}^{2}+1\equiv2\ \mathrm{mod}\ 3$. This contradiction implies that $\widetilde{\varepsilon}$ can only take $1$, and we obtain 
		$$\frac{\widetilde{C}}{n}\cdot\frac{m}{n}\equiv
		-\widetilde{A}^{-1}l\left(\frac{\widetilde{A}^{2}-1}{n}+\widetilde{B}\cdot\frac{\widetilde{C}}{n}\right)\ \mathrm{mod}\ 3.$$
		This equation automatically holds if $\widetilde{C}\equiv\widetilde{A}^{2}-1\equiv0\ \mathrm{mod}\ 3d$.
	\end{compactenum}
	Therefore, we obtain no further restriction to $m$, $n$ and $l$ from the comparison of invariant $s_{3}$.
	
	To summarize,  given non-zero integers $m$, $n$, $l$ and $\overline{n}$, $\overline{l}$ with $\gcd(m,n,l)=\gcd\left(\overline{n},\overline{l}\right)=1$ and $(m,n,l)\equiv\left(0,\overline{n},\overline{l}\right)\equiv(0,1,0)\ \mathrm{mod}\ 2$, manifolds $M_{m,n,l}$ and $M_{0,\overline{n},\overline{l}}$ are diffeomorphic if and only if
	$$
		\left\{
			\begin{alignedat}{2}
				\overline{l}&=l,\\
				\overline{n}&=\pm n,\\
				m&\equiv0\ \mathrm{mod}\ \frac{2^{\lambda_{2}(n)}7^{\lambda_{7}(n)}n^{2}}{3^{\mu_{3}(n)}}.
			\end{alignedat}
		\right.
	$$
	Now the proof of Theorem \ref{Theorem: Partial classification of mfds in M_1}, Statement 2 is completed.
\end{pf}

\begin{pf}[of Theorem \ref{Theorem: Partial classification of mfds in M_1}, Statement 3]
	Let $m$, $n$, $l$ and $\overline{m}$, $\overline{n}$, $\overline{l}$ be non-zero integers such that $\gcd(m,n)=\gcd\left(\overline{m},\overline{n}\right)=1$ and $(m,n,l)\equiv\left(\overline{m},\overline{n},\overline{l}\right)\equiv(0,1,0)\ \mathrm{mod}\ 2$. Now we determine when $M=M_{m,n,l}$ and $\overline{M}=M_{\overline{m},\overline{n},\overline{l}}$ are diffeomorphic. By Proposition \ref{Proposition: M1 mfds, first classification comparing H^4}, Statement 3 they are diffeomorphic only if $\overline{n}=\pm n$ and $\overline{l}=\pm l$. 
	
	Let us first sketch how to further classify such manifolds by comparing the $s$-invariants. The proof is similar to that of Theorem \ref{Theorem: Partial classification of mfds in M_1}, Statement 2 for spin manifolds. Suppose $M\xrightarrow{\varphi}\overline{M}$ is a diffeomorphism. According to the proof of Proposition \ref{Proposition: M1 mfds, first classification comparing H^4}, Statement 3 we have a polarized diffeomorphism
	$$\left({M};{A}{p}^{*}{z}+{B}{p}^{*}\gamma,{C}{p}^{*}{z}+{D}{p}^{*}\gamma\right)\xrightarrow{\varphi}\left(\overline{M};\overline{p}^{*}\overline{z},\overline{p}^{*}{\gamma}\right).$$
	Here $g=\begin{pmatrix}
		{A} & {C} \\ {B} & {D}
	\end{pmatrix}\in GL(2,\mathbb{Z})$. Since $\varphi^{*}$ is a tangential isomorphism $g$ satisfies certain equation $\Phi(g)=0$. Then $M$ and $\overline{M}$ are diffeomorphic if and only if equation \eqref{Equation: Phi(g)=0, s(M;(x,y)g)=s(M bar; x bar, y bar)}, viewed as an equation in $g\in GL(2,\mathbb{Z})$ parametrized by $m$, $n$, $l$ and $\overline{m}$, $\overline{n}$, $\overline{l}$ has a solution in $GL(2,\mathbb{Z})$.
	
	Invariant $S_{10}$ vanishes for both manifolds. Hence it remains to compare invariants $s_{1}\sim s_{9}$ and combine the known result (formula \eqref{Equation: M_1^+ mfds, m bar, m neq0, d bar=d=1, tangential isom}) from tangential isomorphism.
	We begin with invariant $s_{9}$ and exclude the possibility that $\overline{l}=-l$. The result of comparing invariant $s_{1}$ is given in Lemma \ref{Lemma: M_1 spin mfds, n bar=pm n, l bar =-l,sufficient condition}. 
	By formula \eqref{Equation: M_1^+ mfds, m bar, m neq0, d bar=d=1, tangential isom} invariants $s_{6}\sim s_{8}$ are already equal. By invariant $s_{5}$ we relate $\overline{m}^{3}$ and $m^{3}$. Comparison of invariants $s_{4}$ and $s_{3}$ impose no further restriction on $m$, $n$, $l$ and $\overline{m}$, $\overline{n}$. Finally by invariant $s_{2}$ we refine the result at $3$-primary part. Then we complete the proof of Theorem \ref{Theorem: Partial classification of mfds in M_1}, Statement 3.
	
	When we solve equations from $s$-invariants, we would require further divisibility conditions on parameters $u$, $v$ satisfying $um+vn=1$ and also on the bar counterpart, so that equation could be simplified. We have justified the plausibility in proof sketch of Theorem \ref{Theorem: Partial classification of mfds in M_1}, Statement 2.
	
	\noindent\textbf{Polarization and tangential isomorphism}
	
	First we  equip the manifolds with polarizations and examine tangential isomorphisms. We continue to use the notations in the proof of the third statement of Proposition \ref{Proposition: M1 mfds, first classification comparing H^4}. Let $M\xrightarrow{\varphi}\overline{M}$ be a diffeomorphism and $\psi=\varphi^{-1}$. Then $H^{*}\left(\overline{M}\right)\xrightarrow{\varphi^{*}}H^{*}(M)$ and $H^{*}\left(M\right)\xrightarrow{\psi^{*}}H^{*}\left(\overline{M}\right)$ are isomorphisms of cohomology rings and they preserve the first Pontryagin classes. We continue the proof of Proposition \ref{Proposition: M1 mfds, first classification comparing H^4}, Statement 3 and deduce further restrictions on $A$, $B$, $C$, $D$ from tangential isomorphism. Let us begin with the cohomology ring. Set
	\begin{align*}
		\varphi^{*}\begin{pmatrix}
			\overline{p}^{*}\overline{z} & \overline{p}^{*}\gamma
		\end{pmatrix}=\begin{pmatrix}
			p^{*}z & p^{*}\gamma
		\end{pmatrix}\begin{pmatrix}
			A & C \\ B & D
		\end{pmatrix},\quad \psi^{*}\begin{pmatrix}
			p^{*}z & p^{*}\gamma
		\end{pmatrix}=\begin{pmatrix}
			\overline{p}^{*}\overline{z} & \overline{p}^{*}\gamma
		\end{pmatrix}\begin{pmatrix}
			\varepsilon D & -\varepsilon C \\ -\varepsilon B & \varepsilon A
		\end{pmatrix}.
	\end{align*}
	Here $\begin{pmatrix}
		A & C \\ B & D
	\end{pmatrix}\in GL(2,\mathbb{Z})$ and $\varepsilon=AD-BC=\pm1$. We have obtained that $\overline{n_{2}}=n_{2}$ and $\overline{l_{1}}=\pm l_{1}$, hence
	\begin{align*}
		H^{4}(M)&=\mathbb{Z}\big/\delta n_{2}\left\{p^{*}\left(z^{2}\right)\right\}\oplus\mathbb{Z}\big/\delta \left|l_{1}\right|\left\{p^{*}\left(\gamma^{2}\right)\right\},\\
		p_{1}(M)&=3m^{2}\left(p^{*}z\right)^{2}+4\left(p^{*}\gamma\right)^{2};
	\end{align*}
	And we also have the bar counterpart. Recall the action of $\varphi^{*}$ on fourth cohomology groups: 
	\begin{align}\label{Equation: M_1 mfds, (m,n)=(m bar, n bar)=1, isom on H^4}
		\left\{
		\begin{alignedat}{2}
			\varphi^{*}\left(\left(\overline{p}^{*}\overline{z}\right)^{2}\right)&=A^{2}\left(p^{*}z\right)^{2}+B^{2}\left(p^{*}\gamma\right)^{2},\\
			\varphi^{*}\left(\left(\overline{p}^{*}\gamma\right)^{2}\right)&=C^{2}\left(p^{*}z\right)^{2}+D^{2}\left(p^{*}\gamma\right)^{2}.
		\end{alignedat}
		\right.
	\end{align}
	Here we did not exhibit the action on mixing term $\left(\overline{p}^{*}\overline{z}\right)\left(\overline{p}^{*}\gamma\right)$ and $\left({p}^{*}{z}\right)\left({p}^{*}\gamma\right)$. This is because we have applied this before in the proof of Proposition \ref{Proposition: M1 mfds, first classification comparing H^4}, Statement 3, where we deduced the following results:
	\begin{equation*}
		\left\{\begin{alignedat}{2}
			\overline{n_{2}}&=n_{2},\quad\overline{l_{1}}=\pm l_{1},\\
			AC&\equiv CD\equiv0\ \mathrm{mod}\ \delta\overline{n_{2}},\\
			AB&\equiv BD\equiv0\ \mathrm{mod}\ \delta \left|l_1\right|.
		\end{alignedat}\right.
	\end{equation*}
	Recall that $AD-BC=\varepsilon=\pm1$ so that $\gcd(A,D)=1$, hence from the above we further obtain that
	\begin{eqnarray}\label{Equation: M_1 mfds, (m,n)=(m bar, n bar)=1, l|B, n^2 | C}
		B\equiv0\ \mathrm{mod}\ \left|l\right|,\quad C\equiv0\ \mathrm{mod}\  n^{2}.
	\end{eqnarray}
	Formulae \eqref{Equation: M_1 mfds, (m,n)=(m bar, n bar)=1, isom on H^4} and \eqref{Equation: M_1 mfds, (m,n)=(m bar, n bar)=1, l|B, n^2 | C} guarantees that $\varphi^{*}$ induces an isomorphism on the fourth cohomology groups. This can be deduced from the following lemma from basic group theory.
	\begin{lem}\label{Lemma: Z^2 endomorphism descend to homom & isom on coordinate-wise quot}
		Let $\mathbb{Z}^{2}\xrightarrow{\mathcal{L}}\mathbb{Z}^{2}$ be a homomorphism of abelian groups that is given in terms of coordinates by
		$\mathcal{L}(x,y)=(ax+cy,bx+dy)$.
		Let $m$, $n$ be positive integers and denote by $\mathbb{Z}\xrightarrow{\rho_{m}}\mathbb{Z}\big/m$ the mod $m$ quotient map, $\rho_{n}$ defined likewise. Denote $G=\left(\mathbb{Z}\big/m\right)\times\left(\mathbb{Z}\big/n\right)$ and $q=\rho_{m}\times\rho_{n}$.
		\begin{compactenum}
			\item There is a homomorphism $G\xrightarrow{\overline{\mathcal{L}}}G$ such that $q\mathcal{L}=\overline{\mathcal{L}}q$ if and only if
			$cn\equiv0\ \mathrm{mod}\ m$ and $bm\equiv0\ \mathrm{mod}\ n.$
			\item Suppose the homomorphism $\overline{\mathcal{L}}$ exists. Then $\overline{\mathcal{L}}$ is an isomorphism if and only if
			$\gcd\left(ad-bc,mn\right)=1.$
		\end{compactenum}
	\end{lem}
	Let us apply Lemma \ref{Lemma: Z^2 endomorphism descend to homom & isom on coordinate-wise quot} to our case. In view of expression of $\varphi^{*}$ with respect to generators, we define $\mathcal{L}$ in terms of coordinates by 
	$\mathcal{L}(x,y)=\left(A^{2}x+C^{2}y,B^{2}x+D^{2}y\right),$
	and set $q=\rho_{n^{2}}\times\rho_{|l|}.$
	Then it follows from formula \eqref{Equation: M_1 mfds, (m,n)=(m bar, n bar)=1, l|B, n^2 | C} that $C^{2}|l|\equiv0\ \mathrm{mod}\ n^{2}$, $B^{2}n^{2}\equiv0\ \mathrm{mod}\ |l|$ and $A^{2}D^{2}-B^{2}C^{2}=1+2\varepsilon BC$ is coprime to $n^{2}|l|$. Hence $\mathcal{L}$ defines an automorphism $\overline{\mathcal{L}}$ of $G$ by Lemma \ref{Lemma: Z^2 endomorphism descend to homom & isom on coordinate-wise quot}. If we view $\varphi^{*}$ as an endomorphism of the abstract group $\left(\mathbb{Z}\big/n^{2}\right)\times\left(\mathbb{Z}\big/|l|\right)$, then it admits a lifting to an endomorphism of $\mathbb{Z}^{2}$ which satisfies all conditions given in Lemma \ref{Lemma: Z^2 endomorphism descend to homom & isom on coordinate-wise quot}. As a consequence, $\varphi^{*}$ is an isomorphism on fourth cohomology groups.
	
	It remains to study the effect on the first Pontryagin classes. From $\varphi^{*}\left(p_{1}\left(\overline{M}\right)\right)=p_{1}(M)$ we have
	$$\left(3A^{2}\overline{m}^{2}+4C^{2}\right)\left(p^{*}z\right)^{2}+\left(3B^{2}\overline{m}^{2}+4D^{2}\right)\left(p^{*}\gamma\right)^{2}=3m^{2}\left(p^{*}z\right)^{2}+4\left(p^{*}\gamma\right)^{2}.$$
	Recall that $\left(p^{*}\gamma\right)^{2}$ has order $|l|$ and $\left(p^{*}z\right)^{2}$ has order $n^{2}$. Also by formula \eqref{Equation: M_1 mfds, (m,n)=(m bar, n bar)=1, l|B, n^2 | C} we have $l$ divides $B$ and $n^{2}$ divides $C$. Hence we obtain
	$$
		\left\{
			\begin{alignedat}{2}
				3A^{2}\overline{m}^{2}&\equiv3m^{2}\ \mathrm{mod}\ n^{2},\\
				4D^{2}&\equiv4\ \mathrm{mod}\ |l|.
			\end{alignedat}
		\right.
	$$
	By symmetry we also have
	$$
	\left\{
	\begin{alignedat}{2}
		3D^{2}{m}^{2}&\equiv3\overline{m}^{2}\ \mathrm{mod}\ n^{2},\\
		4A^{2}&\equiv4\ \mathrm{mod}\ |l|,
	\end{alignedat}
	\right.
	$$
	which indeed is a direct consequence of $AD-BC=\varepsilon=\pm1$ and formula \eqref{Equation: M_1 mfds, (m,n)=(m bar, n bar)=1, l|B, n^2 | C}.
	
	To summarize, if $M\xrightarrow{\varphi}\overline{M}$ is a diffeomorphism with 
	$$\varphi^{*}\begin{pmatrix}
		\overline{p}^{*}\overline{z} & \overline{p}^{*}{\gamma}
	\end{pmatrix}=\begin{pmatrix}
		p^{*}z & p^{*}\gamma
	\end{pmatrix}g,\quad g=\begin{pmatrix}
		A & C \\ B & D
	\end{pmatrix}\in GL(2,\mathbb{Z}),$$
	then the equation $\Phi(g)=0$ reads
	\begin{eqnarray}\label{Equation: M_1^+ mfds, m bar, m neq0, d bar=d=1, tangential isom}
		\left\{
			\begin{alignedat}{2}
				C&\equiv0\ \mathrm{mod}\ n^{2},\\
				\gcd\left(n^{2},A\right)&=\gcd\left(n^{2},D\right)=1,\\
				3D^{2}m^{2}&\equiv3\overline{m}^{2}\ \mathrm{mod}\ n^{2},\\
				3A^{2}\overline{m}^{2}&\equiv3{m}^{2}\ \mathrm{mod}\ n^{2};
			\end{alignedat}
		\right.
		\quad
		\left\{
		\begin{alignedat}{2}
			B&\equiv0\ \mathrm{mod}\ |l|,\\
			\gcd\left(l,A\right)&=\gcd\left(l,D\right)=1,\\
			4D^{2}&\equiv4\ \mathrm{mod}\ |l|,\\
			4A^{2}&\equiv4\ \mathrm{mod}\ |l|.
		\end{alignedat}
		\right.
	\end{eqnarray}
	
	\noindent\textbf{Compare invariant $s_{9}$ and deduce $\overline{l}=l$}
	
	Now we compare invariant $s_{9}$ and show that $\overline{l}=l$. We read invariants $s_{9}$ of these two polarized manifolds from formula \eqref{Equation: s-inv of spin M_{m,n,l}, (m,n)=1}. Let them be equal and we obtain the following equation:
	\begin{eqnarray}\label{Equation: M_1^+ mfds M_{m,n,l} & M_{m bar,n bar, l bar}, (m,n)=(m bar,n bar)=1, inv s_9}
		-C^{4}u^{3}(1+3vn)|l|+D^{4}n^{2}\mathrm{sgn}(l)\equiv n^{2}\mathrm{sgn}\left(\overline{l}\right)\ \mathrm{mod}\ 2n^{2}|l|.
	\end{eqnarray}
	Recall that $C\equiv0\ \mathrm{mod}\ n^{2}$. Also recall that $um+vn=1$ and $m$ is even. Hence $vn$, $1+3vn$ are both even and we have $-C^{4}u^{3}(1+3vn)|l|\equiv0\ \mathrm{mod}\ 2n^{2}|l|$. Equation \eqref{Equation: M_1^+ mfds M_{m,n,l} & M_{m bar,n bar, l bar}, (m,n)=(m bar,n bar)=1, inv s_9} can be further reduced as
	$D^{4}\equiv\mathrm{sgn}(l)\mathrm{sgn}\left(\overline{l}\right)\ \mathrm{mod}\ 2|l|$.
	Since $l$ is even, we obtain that $\mathrm{sgn}(l)\mathrm{sgn}\left(\overline{l}\right)$ must be $1$ and $\overline{l}=l$. Then we combine the result $4D^{2}\equiv4\ \mathrm{mod}\ |l|$ and obtain $D^{2}\equiv1\ \mathrm{mod}\ |l|$. 
	This has been proved before in the proof of Theorem \ref{Theorem: Partial classification of mfds in M_1}, Statements 1 and 2 for spin manifolds, where we excluded the possibility $\overline{l}=-l$ and deduced $A^{2}\equiv1\ \mathrm{mod} |l|$. 
	
	Therefore, given non-zero integers $m$, $n$, $l$ and $\overline{m}$, $\overline{n}$, $\overline{l}$ with $\gcd(m,n)=\gcd\left(\overline{m},\overline{n}\right)=1$ and $(m,n,l)\equiv\left(\overline{m},\overline{n},\overline{l}\right)\equiv(0,1,0)\ \mathrm{mod}\ 2$, $M=M_{m,n,l}$ and $\overline{M}=M_{\overline{m},\overline{n},\overline{l}}$ are diffeomorphic only if $\overline{n}=\pm n$ and $\overline{l}=l$.
	
	\noindent\textbf{Compare invariants $s_{6}\sim s_{8}$}
	
	Next we consider invariants $s_{6}\sim s_{8}$. As for the polarized $\mathcal{E}_{1}^{+}$-manifold $\left(\overline{M};\overline{p}^{*}\overline{z},\overline{p}^{*}{\gamma}\right)$, invariants $s_{6}\sim s_{8}$ are all already $0\in\mathbb{Q}\big/\mathbb{Z}$. Meanwhile, invariants $s_{6}\sim s_{8}$ of $\left(M;Ap^{*}z+Bp^{*}\gamma,Cp^{*}z+Dp^{*}\gamma\right)$ are also $0\in\mathbb{Q}\big/\mathbb{Z}$ since $B\equiv0\ \mathrm{mod} |l|$, $C\equiv0\ \mathrm{mod}\ n^{2}$ (formula \eqref{Equation: M_1^+ mfds, m bar, m neq0, d bar=d=1, tangential isom}) and $1+3vn$ is even. Therefore, invariants $s_{6}\sim s_{8}$ of both polarized manifolds are equal.

	\noindent\textbf{Compare invariant $s_{5}$}
	
	Now we compare invariant $s_{5}$. Since $um+vn=1$ $m$ is even and $n$ is odd, we may assume a priori that $\alpha(2,u)\gg1$ and $n$ divides $v$. We also assume the bar counterpart. Then we obtain the equation
	$$A^{4}u^{3}\equiv\overline{u}^{3}\ \mathrm{mod}\ n^{2}.$$
	Since $um+\frac{v}{n}\cdot n^{2}=1$ we have $um\equiv1\ \mathrm{mod}\ n^{2}$, and likewise $\overline{u}\overline{m}\equiv1\ \mathrm{mod}\ n^{2}$. Hence we also have
	$$A^{4}\overline{m}^{3}\equiv m^{3}\ \mathrm{mod}\ n^{2}.$$
	By symmetry we also have
	$$D^{4}\overline{u}^{3}\equiv{u}^{3}\ \mathrm{mod}\ n^{2},\quad D^{4}{m}^{3}\equiv \overline{m}^{3}\ \mathrm{mod}\ n^{2}.$$
	
	\noindent\textbf{Compare invariant $s_{4}$}
	
	Next we compare invariant $s_{4}$. We have
	\begin{align}\label{Equation: M_1^+ mfds M_{m,n,l} & M_{m bar,n bar, l bar}, (m,n)=(m bar,n bar)=1, inv s_4}
		\left.
		\begin{alignedat}{2}
			&C^{2}u-3C^{2}u\cdot\frac{v}{n}-3\cdot\frac{C^{2}}{n^{2}}\cdot u-3C^{2}u\cdot\frac{v}{n}+2\cdot\frac{C^{4}}{n^{2}}\cdot u^{3}+6C^{4}u^{3}\cdot\frac{v}{n}\\
			+&\frac{D^{2}-1}{l}\cdot\left(l^{2}+4\right)-2\left(D^{2}+1\right)\cdot\frac{D^{2}-1}{l}\equiv0\ \mathrm{mod}\ 48.
		\end{alignedat}
		\right.
	\end{align}
	By assumption $n$ divides $v$. Meanwhile we have seen that $D^{2}\equiv1\ \mathrm{mod}\ |l|$ and $C\equiv0\ \mathrm{mod}\ n^{2}$. Hence in the left hand side of equation \eqref{Equation: M_1^+ mfds M_{m,n,l} & M_{m bar,n bar, l bar}, (m,n)=(m bar,n bar)=1, inv s_4}, factors of terms are integers. By Chinese remainder theorem we treat the $2$-primary part and $3$-primary part separately, and we begin with the $2$-primary part. We have assumed that $\alpha(2,u)\gg1$. Then the $2$-primary part is reduced as
	$$\frac{D^{2}-1}{l}\cdot\left(2\left(D^{2}-1\right)-l^{2}\right)\equiv0\ \mathrm{mod}\ 16.$$
	Since $AD-BC=\pm1$, $l$ divides $B$ and $l$ is even, we obtain that $D$ is odd and $D^{2}\equiv1\ \mathrm{mod}\ 8$. Hence the $2$-primary part equation can be further reduced as
	$\left(D^{2}-1\right)l\equiv0\ \mathrm{mod}\ 16,$
	which is automatically true.
	
	Next we consider the $3$-primary part. Again we use the fact that $r^{3}\equiv r\ \mathrm{mod}\ 3$ for any integer $r$ and obtain
	$$
		\frac{D^{2}-1}{l}\cdot\left(l^{2}+D^{2}-1\right)\equiv0\ \mathrm{mod}\ 3.
	$$
	The expression on the left hand side is always divided by $3$. See the proof of Theorem \ref{Theorem: Partial classification of mfds in M_1}, Statement 2 for spin manifolds, \textbf{Compare invariant $s_{2}$}, the case $n$ is not divided by $3$.
	
	Therefore, the invariant $s_{4}$ of both polarized $\mathcal{E}_{1}^{+}$-manifolds under consideration are always equal.
	
	\noindent\textbf{Compare invariant $s_{2}$}
	
	We move to compare invariant $s_{2}$. We assume $\alpha(2,u),\alpha\left(2,\overline{u}\right)\gg1$ as before, and we are solving
	\begin{align}\label{Equation: M_1^+ mfds M_{m,n,l} & M_{m bar,n bar, l bar}, (m,n)=(m bar,n bar)=1, inv s_2}
		A^{2}u-\overline{u}+\left(\frac{B}{l}\right)^{2}l\left(l^{2}+4\right)-2\left(\frac{B}{l}\right)^{4}l^{3}+2\cdot\frac{A^{4}u^{3}-\overline{u}^{3}}{n^{2}}-\frac{3\left(A^{2}u-\overline{u}\right)}{n^{2}}\equiv0\ \mathrm{mod}\ 48.
	\end{align}
	We have seen that $B$ is divided by $l$ from tangential isomorphism and that $A^{4}u^{3}\equiv\overline{u}^{3}\ \mathrm{mod}\ n^{2}$ from comparison of invariant $s_{5}$. Hence equation \eqref{Equation: M_1^+ mfds M_{m,n,l} & M_{m bar,n bar, l bar}, (m,n)=(m bar,n bar)=1, inv s_2} implies that $\frac{3\left(A^{2}u-\overline{u}\right)}{n^{2}}$ is an integer, namely
	$
		3\left(A^{2}u-\overline{u}\right)\equiv0\ \mathrm{mod}\ n^{2},
	$
	or equivalently
	\begin{align}\label{Equation: M_1^+ mfds M_{m,n,l} & M_{m bar,n bar, l bar}, (m,n)=(m bar,n bar)=1, inv s_2, first step mod 48 each term should be integer}
		3\left(A^{2}\overline{m}-{m}\right)\equiv0\ \mathrm{mod}\ n^{2}.
	\end{align}
	The treatment depends on divisibility of $n$ by $3$, and we can write them uniformly as
	$$A^{2}\overline{m}\equiv m\ \mathrm{mod}\ \frac{n^{2}}{3^{\mu_{3}(n)}}.$$
	We combine the known result $A^{4}\overline{m}^{3}\equiv m^{3}\ \mathrm{mod}\ n^{2}$ and obtain that
	$
		m^{3}\equiv\left(A^{2}\overline{m}\right)^{2}\overline{m}\equiv m^{2}\overline{m}\ \mathrm{mod}\ \frac{n^{2}}{3^{\mu_{3}(n)}}.
	$
	By assumption $m$ and $n$ are coprime, hence from the above we obtain that $$\overline{m}\equiv m\ \mathrm{mod}\ \frac{n^{2}}{3^{\mu_{3}(n)}}.$$ 
	When $n$ is not divided by $3$ the result above is already contained by formula \eqref{Equation: M_1^+ mfds, n bar=pm n, l bar=l, sol}, the result of comparin invariant $s_{1}$. And when $n$ is divided by $3$ the result above refines the $3$-primary part of formula \eqref{Equation: M_1^+ mfds, n bar=pm n, l bar=l, sol}.
	
	We return to equation \eqref{Equation: M_1^+ mfds M_{m,n,l} & M_{m bar,n bar, l bar}, (m,n)=(m bar,n bar)=1, inv s_2}. Factors of terms on the left hand side are all integer. By Chinese remainder theorem we treat $2$-primary part and $3$-primary part separately. We begin with the $2$-primary part. We have assumed that $\alpha(2,u)\gg1$ and $\alpha\left(2,\overline{u}\right)\gg1$. Since $n$ is odd and $l$ is even, the $2$-primary part can be reduced as
	$$\left(\frac{B}{l}\right)^{2}\cdot\frac{l}{2}\cdot\left(\left(\frac{l}{2}\right)^{2}+1\right)\equiv0\ \mathrm{mod}\ 2.$$
	The equation above automatically holds, since for any integer $r$ we have $r^{2}\equiv r\ \mathrm{mod}\ 2$.
	
	Then we consider the $3$-primary part. We have the following equation
	\begin{align*}
		\left.
		\begin{alignedat}{2}
			A^{2}u-\overline{u}-\frac{A^{4}u^{3}-\overline{u}^{3}}{n^{2}}-\frac{3\left(A^{2}u-\overline{u}\right)}{n^{2}}\equiv0\ \mathrm{mod}\ 3,
		\end{alignedat}
		\right.
	\end{align*}
	or equivalently,
	\begin{align*}
		\left.
		\begin{alignedat}{2}
			\left(A^{2}u-\overline{u}\right)n^{2}-\left(A^{4}u^{3}-\overline{u}^{3}\right)-3\left(A^{2}u-\overline{u}\right)\equiv0\ \mathrm{mod}\ 3n^{2},
		\end{alignedat}
		\right.
	\end{align*}
	We have set $n=3^{b}7^{c}n_{0}$ before. Here it suffices to discuss $3$-primary part separately, and we set $n_{3}=7^{c}n_{0}$ so that $n=3^{b}n_{3}$ with $n_{3}$ coprime to $3$. By Chinese remainder theorem it is equivalent to solve the following equation set
		\begin{align}\label{Equation: M_1^+ mfds M_{m,n,l} & M_{m bar,n bar, l bar}, (m,n)=(m bar,n bar)=1, inv s_2, 3-primary}
			\left\{
			\begin{alignedat}{2}
				\left(A^{2}u-\overline{u}\right)n^{2}-\left(A^{4}u^{3}-\overline{u}^{3}\right)-3\left(A^{2}u-\overline{u}\right)&\equiv0\ \mathrm{mod}\ 3^{2b+1},\\
				\left(A^{2}u-\overline{u}\right)n^{2}-\left(A^{4}u^{3}-\overline{u}^{3}\right)-3\left(A^{2}u-\overline{u}\right)&\equiv0\ \mathrm{mod}\ n_{3}^{2}.
			\end{alignedat}
			\right.
		\end{align}
	We have obtained from previous results that $A^{4}u^{3}\equiv\overline{u}^{3}\ \mathrm{mod}\ 3^{2b}n_{3}^{2}$ and $A^{2}u\equiv\overline{u}\ \mathrm{mod}\ 3^{2b-\mu_{3}(n)}n_{3}^{2}$. Hence the equation with modulus $n_{3}^{2}$ of \eqref{Equation: M_1^+ mfds M_{m,n,l} & M_{m bar,n bar, l bar}, (m,n)=(m bar,n bar)=1, inv s_2, 3-primary} automatically holds, and solution to the equation with modulus $3^{2b+1}$ depends on the value of $b$:
	\begin{compactenum}
		\item When $b=0$ we have
			$$\left(A^{2}u-\overline{u}\right)n^{2}-\left(A^{4}u^{3}-\overline{u}^{3}\right)\equiv0\ \mathrm{mod}\ 3.$$
			The left hand side can be further reduced as $\left(A^{2}u-\overline{u}\right)\left(n^{2}-1\right)\ \mathrm{mod}\ 3$ as $r^{3}\equiv r\ \mathrm{mod}\ 3$. Now $n=n_{3}$ is coprime to $3$ and $n^{2}\equiv1\ \mathrm{mod}\ 3$, hence the equation above automatically holds.
		\item When $b\geqslant1$ we already have $\overline{m}\equiv m$, $\overline{u}\equiv u\ \mathrm{mod}\ 3^{2b-1}$, $A^{2}u\equiv\overline{u}\ \mathrm{mod}\ 3^{2b-1}$ and $A^{4}u^{3}\equiv\overline{u}^{3}\ \mathrm{mod}\ 3^{2b}$. Hence it suffices to solve the equation
		$$\left(A^{4}u^{3}-\overline{u}^{3}\right)+3\left(A^{2}u-\overline{u}\right)\equiv0\ \mathrm{mod}\ 3^{2b+1}.$$
		We set
		$\overline{u}=u+3^{2b-1}\varepsilon_{1}+3^{2b}k_{1}$, $\varepsilon_{1}\in\{\pm1,0\}$, $k_{1}\in\mathbb{Z}$,
		substitute into the following equivalent equation
		$$\left(A^{2}+1\right)\left(A^{2}-1\right)u^{3}+3\left(A^{2}-1\right)u-\left(\overline{u}^{3}-u^{3}\right)-3\left(\overline{u}-u\right)\equiv0\ \mathrm{mod}\ 3^{2b+1},$$
		and it remains to solve
		\begin{align}\label{Equation: M_1^+ mfds M_{m,n,l} & M_{m bar,n bar, l bar}, (m,n)=(m bar,n bar)=1, inv s_2, 3-primary, b>=1}
			\left(A^{2}-1\right)\left(\left(A^{2}+1\right)u^{3}+3u\right)+3^{2b}\varepsilon_{1}u^{2}\equiv0\ \mathrm{mod}\ 3^{2b+1}.
		\end{align}
		First we consider equation \eqref{Equation: M_1^+ mfds M_{m,n,l} & M_{m bar,n bar, l bar}, (m,n)=(m bar,n bar)=1, inv s_2, 3-primary, b>=1} with modulus $3^{2b}$. We have
		$$\left(A^{2}-1\right)\left(\left(A^{2}+1\right)u^{3}+3u\right)\equiv0\ \mathrm{mod}\ 3^{2b}.$$
		It follows from $\overline{u}\equiv u\ \mathrm{mod}\ 3^{2b-1}$ and $A^{2}u\equiv\overline{u}\ \mathrm{mod}\ 3^{2b-1}$ that $\left(A^{2}-1\right)u\equiv0\ \mathrm{mod}\ 3^{2b-1}$. Since $um+vn=1$ and $b\geqslant1$ so that $n$ is divided by $3$, we have $u$ is coprime to $3$, $A^{2}-1\equiv0\ \mathrm{mod}\ 3^{2b-1}$ and $A^{2}+1$ is also coprime to $3$. Then the equation above implies 
		$$A^{2}-1\equiv0\ \mathrm{mod}\ 3^{2b}.$$
		We return to equation \eqref{Equation: M_1^+ mfds M_{m,n,l} & M_{m bar,n bar, l bar}, (m,n)=(m bar,n bar)=1, inv s_2, 3-primary, b>=1}. Recall that $AD-BC=\pm1$ and $C\equiv0\ \mathrm{mod}\ n^{2}$, so that $A$ is coprime to $3$ and $A^{2}+1\equiv-1\ \mathrm{mod}\ 3$. Hence we have
		$$\varepsilon_{1}m\equiv\frac{A^{2}-1}{3^{2b}}\ \mathrm{mod}\ 3.$$
		Recall that $m$ is not divided by $3$. Solving the equation above means that we wish to find $\varepsilon_{1}$ and $m$ such that the equation above has a solution in $A$, $B$, $C$ and $D$ with $AD-BC=\varepsilon=\pm1$. Here $\varepsilon_{1}$ is determined by $\overline{m}$ and $m$. When $\varepsilon_{1}=0$ we require $A^{2}\equiv1\ \mathrm{mod}\ 3^{2b+1}$, and when $\varepsilon_{1}=\pm1$ we just require that $\frac{A^{2}-1}{3^{2b}}\equiv\varepsilon_{1}m\ \mathrm{mod}\ 3$. Hence this equation does not impose any further restriction on $\overline{m}$ and $m$.
	\end{compactenum}
	
	Therefore, by comparing invariant $s_{2}$ we obtain
	$$\overline{m}\equiv m\ \mathrm{mod}\ \frac{n^{2}}{3^{\mu_{3}(n)}},$$
	which refines the $3$-primary part of our previous result. If we combine the results from tangential isomorphism and comparisons of invariants $s_{1}$, $s_{2}$, $s_{4}\sim s_{9}$ and $S_{10}$, we obtain formula \eqref{Equation: M_1^+ mfds, n bar=pm n, l bar=l, sol}.
		
	\noindent\textbf{Compare invariant $s_{3}$}
	
	Finally we compare invariant $s_{3}$. We are solving
	\begin{eqnarray}\label{Equation: M_1^+ mfds M_{m,n,l} & M_{m bar,n bar, l bar}, (m,n)=(m bar,n bar)=1, inv s_3}
		ACu+\frac{B}{l}\cdot D\left(l^{2}+4\right)+4A^{3}\cdot\frac{C}{n^{2}}\cdot u^{3}+4A\cdot\frac{C}{n^{2}}\cdot C^{2}u^{3}-4\cdot\frac{B}{l}\cdot B^{2}D-4\cdot\frac{B}{l}\cdot D^{3}\equiv0\ \mathrm{mod}\ 24.
	\end{eqnarray}
	According to our previous results, factors of terms in the left hand side are integers. By Chinese remainder theorem we treat the $2$-primary part and $3$-primary part separately, and we begin with the $2$-primary part. We assume $\alpha(2,u)\gg1$ as before. We also recall that $l$ is even and divides $B$, so that $B$ is also even. By $AD-BC=\pm1$ we have $D$ is odd, hence $D^{2}-1\equiv0\ \mathrm{mod}\ 8$. Then the $2$-primary part equation is reduced as
	$\frac{B}{l}\cdot\frac{l}{2}\equiv0\ \mathrm{mod}\ 2$.
	The treatment depends on parity of $\frac{l}{2}$:
	\begin{compactenum}
		\item When $\frac{l}{2}$ is odd, we have 
		$\frac{B}{l}\equiv0\ \mathrm{mod}\ 2.$
		\item When $\frac{l}{2}$ is even, the equation is automatically true.
	\end{compactenum}
	Therefore, the $2$-primary part of equation \eqref{Equation: M_1^+ mfds M_{m,n,l} & M_{m bar,n bar, l bar}, (m,n)=(m bar,n bar)=1, inv s_3} imposes no further restriction on $m$, $n$, $l$ and $\overline{m}$, $\overline{n}$.
	
	Then we consider the $3$-primary part. Since $r^{3}\equiv r\ \mathrm{mod}\ 3$ it remains to solve
	\begin{eqnarray}\label{Equation: M_1^+ mfds M_{m,n,l} & M_{m bar,n bar, l bar}, (m,n)=(m bar,n bar)=1, inv s_3, 3-primary}
		A\cdot\frac{C}{n^{2}}\cdot u\left(n^{2}-1\right)\equiv0\ \mathrm{mod}\ 3.
	\end{eqnarray}
	The treatment depends on the divisibility of $n$ by $3$:
	\begin{compactenum}
		\item When $n$ is not divided by $3$, we have $n^{2}\equiv1\ \mathrm{mod}\ 3$ and equation \eqref{Equation: M_1^+ mfds M_{m,n,l} & M_{m bar,n bar, l bar}, (m,n)=(m bar,n bar)=1, inv s_3, 3-primary} is automatically true.
		\item When $n$ is divided by $3$ so that $u$ and $A$ are not, we have
		$\frac{C}{n^{2}}\equiv0\ \mathrm{mod}\ 3$.
	\end{compactenum}
	Again the $3$-primary part of equation \eqref{Equation: M_1^+ mfds M_{m,n,l} & M_{m bar,n bar, l bar}, (m,n)=(m bar,n bar)=1, inv s_3} imposes no further restriction on $m$, $n$, $l$ and $\overline{m}$, $\overline{n}$.
	
	To summarize, given non-zero integers $m$, $n$, $l$ and $\overline{m}$, $\overline{n}$, $\overline{l}$ with $\gcd(m,n)=\gcd\left(\overline{m},\overline{n}\right)=1$ and $(m,n,l)\equiv\left(\overline{m},\overline{n},\overline{l}\right)\equiv(0,1,0)\ \mathrm{mod}\ 2$, manifolds $M_{m,n,l}$ and $M_{\overline{m},\overline{n},\overline{l}}$ are diffeomorphic if and only if
	$$
		\left\{
			\begin{alignedat}{2}
				\overline{l}&=l,\\
				\overline{n}&=\pm n,\\
				\overline{m}&\equiv m\ \mathrm{mod}\ \frac{2^{\lambda_{2}(n)}7^{\lambda_{7}(n)}n^{2}}{3^{\mu_{3}(n)}}.
			\end{alignedat}
		\right.
	$$
	Now the proof of Theorem \ref{Theorem: Partial classification of mfds in M_1}, Statement $3$ is completed.
\end{pf}

	\section{Application to positive Ricci curvature problem}\label{Section: 8 Application to positive Ricci curvature}

In this section we will show that there are spin manifolds $M_{m,n,l}$ in $\mathcal{M}_{1}$ whose spaces and moduli spaces of positive Ricci curvature metrics have infinitely many components. 

First in Section \ref{Section 8.1: Two specific metrics on M_{m,n,l}} we explore two metrics $g_{m,n,l}$ and $h_{m,n,l}$ on $M_{m,n,l}$, such that $\left(M_{m,n,l},g_{m,n,l}\right)$ has positive Ricci curvature (\cite{PBJPWT98PositiveRicciPrinBdl}), that $\left(M_{m,n,l},h_{m,n,l}\right)$ has positive scalar curvature (\cite[Proposition 4.3, Formula (4.4)]{KreckStolz1993}, \cite[9.59 Theorem, 9.70 Proposition]{BesseEinsteinMfd}) and the $s$-invariant of $\left(M_{m,n,l},h_{m,n,l}\right)$ can be computed from \cite[Theorem 3.11, Formula (3.12)]{KreckStolz1993}. We justify that after suitable rescaling $h_{m,n,l}$ along the fiber direction, there is a continuous path in $\mathfrak{R}_{scal>0}\left(M_{m,n,l}\right)$ that joins $g_{m,n,l}$ and $h_{m,n,l}$, which also defines a continuous path in $\mathfrak{M}_{scal>0}\left(M_{m,n,l}\right)$. Next in Section \ref{Section 8.2: Infinitely many path components of space and moduli space of metrics with positive Ricci curvature} we compute $s\left(M_{m,n,l},h_{m,n,l}\right)$ and combine the classification result (Theorem \ref{Theorem: Partial classification of mfds in M_1}) to show that there is a manifold $M_{m,n,l}$ such that $\mathfrak{R}_{scal>0}\left(M_{m,n,l}\right)$ and $\mathfrak{M}_{scal>0}\left(M_{m,n,l}\right)$ both have infinitely many components.

\subsection{Metrics $g_{m,n,l}$ and $h_{m,n,l}$ on $M_{m,n,l}$}\label{Section 8.1: Two specific metrics on M_{m,n,l}}

We begin with the construction of metrics $g_{m,n,l}$ and $h_{m,n,l}$. Recall that $N=\left(\mathbb{C}P^{1}\times\mathbb{C}P^{2}\right)\#\mathbb{C}P^{3}$ admits a core metric $g_{N}$ (see the two paragraphs after Problem \ref{Problem: Classify M_{m,n,l}}), and in particular it has positive Ricci curvature and thus positive scalar curvature. Recall also that $M_{m,n,l}$ is simply connected and the bundle projection $M_{m,n,l}\xrightarrow{p}N$ can be viewed as a principal $S^{1}$-bundle over $N$ with Euler class $m\alpha+n\beta+l\gamma$. Fix a principal connection $\mathcal{H}\subset\tau_{M_{m,n,l}}$, and let $\theta\in\Omega^{1}\left(M_{m,n,l};\mathfrak{u}(1)\right)\cong\Omega^{1}\left(M_{m,n,l}\right)$ be the corresponding principal connection $1$-form. Then under local coordinates $(t,x)=\left(t,x_{1},\cdots,x_{n}\right)$ we have $\theta=dt+A(x)$. Here $t$ is the coordinate corresponding to the fiber, $x=\left(x_{1},\cdots,x_{n}\right)$ corresponds to the base space and $A(x)\in\Omega^{1}(N)$ and the Euler class of the bundle is $m\alpha+n\beta+l\gamma=\left[\frac{F}{2\pi}\right]$ with $F=dA$.

By \cite[0.1 Theorem]{PBJPWT98PositiveRicciPrinBdl} $M_{m,n,l}$ admits an $S^{1}$-invariant metric $g_{m,n,l}$ with positive Ricci curvature, such that the bundle projection $\left(M_{m,n,l},g_{m,n,l}\right)\xrightarrow{p}\left(N,g_{N}\right)$ is a Riemmanian submersion with horizontal distribution $\mathcal{H}$. Under local coordinates we have
\begin{eqnarray}\label{Equation: metric g_{m,n,l}}
	g_{m,n,l}=p^{*}g_{N}+\varphi(x)^{2}\left(dt+A(x)\right)^{2}.
\end{eqnarray}
Here $\varphi(x)=e^{f(x)}$ for certain delicately chosen smooth function $f$ on $N$ and is always positive (see \cite{PBJPWT98PositiveRicciPrinBdl} for more details). Then it can be deduced from \cite[9.37 Corollary]{BesseEinsteinMfd} that the scalar curvature of $\left(M_{m,n,l},g_{m,n,l}\right)$ and $\left(N,g_{N}\right)$ are related by
\begin{eqnarray}\label{Equation: scalar curvature of g_{m,n,l}}
	scal\left(M_{m,n,l},g_{m,n,l}\right)=p^{*}scal\left(N,g_{N}\right)-\frac{\varphi^{2}|F|^{2}}{4}-\frac{2\Delta\varphi}{\varphi},
\end{eqnarray}
where $\Delta=\Delta_{N}$ is the Laplacian operator on $\left(N,g_{N}\right)$.

By \cite[Proposition 4.3]{KreckStolz1993} (see also \cite[9.59 Theorem]{BesseEinsteinMfd}) $M_{m,n,l}$ admits the unique metric $h_{m,n,l}$ such that $\left(M_{m,n,l},h_{m,n,l}\right)\xrightarrow{p}\left(N,g_{N}\right)$ is a Riemmanian submersion with totally geodesic fibers isometric to $\left(S^{1}, g_{std}(r)\right)$ and horizontal distribution again $\mathcal{H}$. Here $g_{std}(r)$ is the standard Euclidean metric on $S^{1}$ so that $\left(S^{1}, g_{std}(r)\right)$ can be isometrically embedded into the $2$-dimensional Euclidean space $\mathbb{E}^{2}$ as the circle centered at the origin with radius $r>0$. Under local coordinates this metric can be expressed as
\begin{eqnarray}\label{Equation: metric h_{m,n,l}}
	h_{m,n,l}=p^{*}g_{N}+r^{2}\left(dt+A(x)\right)^{2}
\end{eqnarray}
and the scalar curvature is
\begin{eqnarray}\label{Equation: scalar curvature of h_{m,n,l}}
	scal\left(M_{m,n,l},g_{m,n,l}\right)=p^{*}scal\left(N,g_{N}\right)-\frac{r^{2}|F|^{2}}{4}.
\end{eqnarray}

We have seen that the Ricci curvature of $\left(N,g_{N}\right)$, denoted by $Ric\left(N,g_{N}\right)$, is posivite. Hence $scal\left(N,g_{N}\right)$ is also positive and we may assume $r$ is small enough so that $h_{m,n,l}$ has positive scalar curvature (\cite[Formula (4.4)]{KreckStolz1993}, see also \cite[9.70 Proposition]{BesseEinsteinMfd}). Meanwhile, by construction $Ric\left(M_{m,n,l},g_{m,n,l}\right)>0$ and thus $scal\left(M_{m,n,l},g_{m,n,l}\right)>0$. Metrics $g_{m,n,l}$ and $h_{m,n,l}$ are in general distinct, and we will show that they lie in the same path component of $\mathfrak{R}_{scal>0}\left(M_{m,n,l}\right)$ and thus also in $\mathfrak{M}_{scal>0}\left(M_{m,n,l}\right)$. We will prove that after requiring $r$ to be smaller if necessary, there is a family of metrics on $M_{m,n,l}$ that joins $g_{0}=g_{m,n,l}$ and $g_{1}=h_{m,n,l}$ among which each has positive scalar curvature.

\begin{lem}
	Suppose further that $0<r<\min_{x\in N}\varphi(x)$. Then there is a path in $\mathfrak{R}_{scal>0}\left(M_{m,n,l}\right)$ joining $g_{m,n,l}$ and $h_{m,n,l}$.
\end{lem}

\begin{pf}
Let $\lambda\in[0,1]$. Let $\left\{g_{\lambda}:\lambda\in[0,1]\right\}$ be the family of metrics on $M_{m,n,l}$ such that under local coordinates
\begin{align*}
	g_{\lambda}&=p^{*}g_{N}+\left(\lambda r+(1-\lambda)\varphi\right)^{2}\left(dt+A(x)\right)^{2},\\
	scal\left(M_{m,n,l},g_{\lambda}\right)&=p^{*}scal\left(N,g_{N}\right)-\frac{\left(\lambda r+(1-\lambda)\varphi\right)^{2}|F|^{2}}{4}-\frac{2(1-\lambda)\Delta\varphi}{\lambda r+(1-\lambda)\varphi}.
\end{align*}
We will prove $scal\left(M_{m,n,l},g_{\lambda}\right)>0$ holds for any $\lambda\in[0,1]$. Set 
$$\rho(\lambda)=\frac{\left(\lambda r+(1-\lambda)\varphi\right)^{2}|F|^{2}}{4}+\frac{2(1-\lambda)\Delta\varphi}{\lambda r+(1-\lambda)\varphi},\ \lambda\in[0,1],$$
and we shall show $\rho(\lambda)<p^{*}scal\left(N,g_{N}\right)$ for any $\lambda\in[0,1]$. 

Since $scal\left(M_{m,n,l},g_{m,n,l}\right)>0$ and $scal\left(M_{m,n,l},h_{m,n,l}\right)>0$, we have $\rho(0)<p^{*}scal\left(N,g_{N}\right)$ and $\rho(1)<p^{*}scal\left(N,g_{N}\right)$. Differentiate with respect to $\lambda$ and we have
\begin{align*}
	\rho'(\lambda)&=-\frac{(\varphi-(\varphi-r)\lambda)(\varphi-r)|F|^{2}}{2}-\frac{2r\Delta\varphi}{\left(\varphi-(\varphi-r)\lambda\right)^{2}},\\
	\rho''(\lambda)&=\frac{(\varphi-r)^{2}|F|^{2}}{2}-\frac{4r(\varphi-r)\Delta\varphi}{(\varphi-(\varphi-r)\lambda)^{3}}.
\end{align*}
When $F(x)\neq0$ the critical point of $\rho$ is 
$$\lambda_{0}=-\sqrt[3]{\frac{4r\Delta\varphi}{(\varphi-r)|F|^{2}}}.$$
Therefore:
\begin{compactenum}
	\item If $\Delta\varphi>0$, then $\lambda_{0}<0$ and $\rho'<0$ on $[0,1]$, $\rho$ strictly decreases on $[0,1]$ and $\rho(\lambda)\leqslant\rho(0)<p^{*}scal\left(N,g_{N}\right)$ for any $\lambda\in[0,1]$.
	\item If $\Delta\varphi\leqslant0$, then $\rho''\geqslant0$ on $[0,1]$, $\rho$ is concave up on $[0,1]$ and $\rho(\lambda)\leqslant\min\left\{\rho(0),\rho(1)\right\}<p^{*}scal\left(N,g_{N}\right)$ for any $\lambda\in[0,1]$.
\end{compactenum}
And when $F(x)=0$ at some point $x\in N$ we can argue similarly and obtain again that $\rho(\lambda)<p^{*}scal\left(N,g_{N}\right)$ for any $\lambda\in[0,1]$.

In conclusion, $\left\{g_{\lambda}\right\}_{0\leqslant\lambda\leqslant1}$ defines a path in $\mathfrak{R}_{scal>0}\left(M_{m,n,l}\right)$ that connects $g_{0}=g_{m,n,l}$ and $g_{1}=h_{m,n,l}$
\end{pf}

\begin{rmk}
	The argument above applies for all integers $m$, $n$ and $l$ that are not simultaneously zero. The construction $h_{m,n,l}$  (\cite[Proposition 4.3]{KreckStolz1993}, see also \cite[9.59 Theorem]{BesseEinsteinMfd}) imposes no topological requirement, and the construction $g_{m,n,l}$  (\cite{PBJPWT98PositiveRicciPrinBdl}) only requires that $M_{m,n,l}$ has finite fundamental group.
\end{rmk}

\subsection{Infinitely many components of moduli space}\label{Section 8.2: Infinitely many path components of space and moduli space of metrics with positive Ricci curvature}

Now we assume that $M_{m,n,l}\in\mathcal{M}_{1}$ is spin, namely $n$, $ l\neq0$ and $(m,n,l)\equiv(0,1,0)\ \mathrm{mod}\ 2$. Then $M_{m,n,l}$ has vanishing rational Pontryagin classes (Propositions \ref{Proposition: 4th cohomology gps of E} and \ref{Proposition: char classes of E}), and if $g$ is a metric on $M_{m,n,l}$ with positive scalar curvature the $s$-invariant $s\left(M_{m,n,l},g\right)\in\mathbb{Q}$ is defined. We will compute the $s$-invariant of $s\left(M_{m,n,l},h_{m,n,l}\right)$ and the result is as follows:
\begin{lem}\label{Lemma: s-inv of Riem mfd (M_{m,n,l},h_{m,n,l})}
	Assume $n$, $ l\neq0$ and $(m,n,l)\equiv(0,1,0)\ \mathrm{mod}\ 2$. Let the metric $h_{m,n,l}$ on $M_{m,n,l}$ be given as in Section \ref{Section 8.1: Two specific metrics on M_{m,n,l}}. The $s$-invariant of $\left(M_{m,n,l},h_{m,n,l}\right)$ is given by
	\begin{eqnarray}\label{Equation: s inv of (M_{m,n,l},h_{m,n,l})}
		s\left(M_{m,n,l},h_{m,n,l}\right)=-\frac{3m\left(n^{2}+3\right)\left(n^{2}-1\right)}{896n^{2}}-\frac{\left(l^{2}+4\right)^{2}}{896l}+\frac{\mathrm{sgn}(l)}{224}\in\mathbb{Q}.
	\end{eqnarray}
\end{lem}

First we prove Lemma \ref{Lemma: s-inv of Riem mfd (M_{m,n,l},h_{m,n,l})}. Then we combine the classification result (Theorem \ref{Theorem: Partial classification of mfds in M_1}) and obtain that there is a manifold $M_{m,n,l}$ such that $\mathfrak{M}_{Ric>0}\left(M_{m,n,l}\right)$ and $\mathfrak{R}_{Ric>0}\left(M_{m,n,l}\right)$ both have infinitely many path components, completing the proof of Theorem \ref{Theorem: mfd whose space & moduli space of Ric+ metrics have infinitely many components}.

\begin{pf}[of Lemma \ref{Lemma: s-inv of Riem mfd (M_{m,n,l},h_{m,n,l})}]
	Since $N$ is non-spin and the principal $S^{1}$-bundle $M_{m,n,l}\to N$ has Euler class $e=m\alpha+n\beta+l\gamma$, we shall apply the second formula of \cite[Theorem 3.11, Formula (3.12)]{KreckStolz1993}:
	\begin{align*}
		s\left(M_{m,n,l},h_{m,n,l}\right)&=-\frac{\left\langle p_{1}(N)\check{p_{1}}(N)+2p_{1}(N)e+e^{3},[N]\right\rangle}{896}+\frac{\sigma\left(H^{2}(N;\mathbb{Q}),(\cdot\cup\cdot\cup e)\cap[N]\right)}{224}.
	\end{align*}
	Here $H^{2}(N;\mathbb{Q})\xrightarrow{\cdot\cup e}H^{4}(N;\mathbb{Q})$ is an isomorphism and there is a unique class $\check{p_{1}}(N)\in H^{2}(N;\mathbb{Q})$ such that $\check{p_{1}}(N)e=p_{1}(N)\in H^{4}(N;\mathbb{Q})$. Since $p_{1}(N)=3\beta^{2}+4\gamma^{2}$ and $e=m\alpha+n\beta+l\gamma$, we have
	$\check{p_{1}}(N)=-\frac{3m}{n^{2}}\alpha+\frac{3}{n}\beta+\frac{4}{l}\gamma,$
	and it is routine to compute that
	\begin{align*}
		\left\langle p_{1}(N)\check{p_{1}}(N)+2p_{1}(N)e+e^{3},[N]\right\rangle=\frac{3m\left(n^{2}+3\right)\left(n^{2}-1\right)}{n^{2}}+\frac{\left(l^{2}+4\right)^{2}}{l}.
	\end{align*}
	The signature term is computed by
	\begin{align*}
		\sigma\left(H^{2}(N;\mathbb{Q}),(\cdot\cup\cdot\cup e)\cap[N]\right)&=\sigma\left(\begin{pmatrix}
			0 & n & 0\\
			n & m & 0\\
			0 & 0 & l
		\end{pmatrix}\right)=\mathrm{sgn}(l).
	\end{align*}
	Combine these two results and we obtain formula \eqref{Equation: s inv of (M_{m,n,l},h_{m,n,l})}.
\end{pf}

\begin{pf}[of Theorem \ref{Theorem: mfd whose space & moduli space of Ric+ metrics have infinitely many components}]
We have seen that $\mathfrak{R}_{scal>0}\left(M_{m,n,l}\right)$ is not empty. Then the $s$-invariant of Riemannian manifolds with positive scalar curvature induces two maps
\begin{alignat*}{4}
	s&:\pi_{0}\left(\mathfrak{R}_{scal>0}\left(M_{m,n,l}\right)\right)&\to&\mathbb{Q},\\
	|s|&:\pi_{0}\left(\mathfrak{M}_{scal>0}\left(M_{m,n,l}\right)\right)&\to&\mathbb{Q}_{\geqslant0}.
\end{alignat*}
Since $h_{m,n,l}$ and $g_{m,n,l}$ are joined by a path in $\mathfrak{R}_{scal>0}\left(M_{m,n,l}\right)$, we have 
\begin{align*}
	s\left(M_{m,n,l},h_{m,n,l}\right)&=s\left(M_{m,n,l},g_{m,n,l}\right)\\
	&=-\frac{3m\left(n^{2}+3\right)\left(n^{2}-1\right)}{896n^{2}}-\frac{\left(l^{2}+4\right)^{2}}{896l}+\frac{\mathrm{sgn}(l)}{224}\in\mathbb{Q}.
\end{align*} 
Notice that $g_{m,n,l}$ has positive Ricci curvature, and if two metrics on $M_{m,n,l}$ which admit positive Ricci curvature lie in distinct path components of $\mathfrak{R}_{scal>0}\left(M_{m,n,l}\right)$ $\left(\text{resp. }\mathfrak{M}_{scal>0}\left(M_{m,n,l}\right)\right)$, then they must lie in in distinct path components of $\mathfrak{R}_{Ric>0}\left(M_{m,n,l}\right)$ $\left(\text{resp. }\mathfrak{M}_{Ric>0}\left(M_{m,n,l}\right)\right)$ as well.

Now we fix an odd number $n\geqslant3$, a positive even number $l$ and an even number $m_{0}>\frac{n^{2}\left(\left(l^{2}+4\right)^{2}-4l\right)}{3\left(n^{2}+3\right)\left(n^{2}-1\right)l}$ such that $\gcd\left(m_{0},n,l\right)=1$. Set $M=M_{m_{0},n,l}$ and $\Delta_{0}=\frac{2^{\lambda_{2}(n)}7^{\lambda_{7}(n)}n^{2}}{3^{\mu_{3}(n)}}$. Then $M_{m_{0}+k\Delta_{0},n,l}$ is diffeomorphic to $M$ for any positive integer $k$, $\left\{g_{m_{0}+k\Delta_{0},n,l}:k\in\mathbb{N}\right\}$ is an infinite family of positive Ricci curvature metrics on $M$ and $\left\{s\left(M_{m_{0}+k\Delta_{0},n,l},g_{m_{0}+k\Delta_{0},n,l}\right):k\in\mathbb{N}\right\}$ is a strictly decreasing sequence of negative rational numbers. Hence $\mathfrak{R}_{scal>0}\left(M\right)$ and $\mathfrak{M}_{scal>0}\left(M\right)$ both have infinitely many path components for such $M=M_{m_{0},n,l}$.
\end{pf}
	
	\bibliography{refs}
	
\end{document}